\documentclass[reqno,12pt]{amsart}
\theoremstyle{plain}    
\newtheorem{thm}{Theorem}[section]
\theoremstyle{plain}    
\newtheorem{cor}[thm]{Corollary} 
\newtheorem{lemma}[thm]{Lemma} 
\newtheorem{prop}[thm]{Proposition}
\theoremstyle{remark}
\newtheorem{defi}[thm]{Definition}

\newtheorem{example}[thm]{Example}

\newtheorem{remark}[thm]{Remark}

\usepackage{amscd,amssymb,comment,eepic,euscript,graphics}

\newcommand\Afr{{\mathfrak A}}
\newcommand\alg{{\operatorname{alg}}}
\newcommand\Bc{{\mathcal{B}}}
\newcommand\Bchk{{\check B}}
\newcommand\Bt{{\widetilde B}}
\newcommand\bt{{\tilde b}}
\newcommand\card{{\operatorname{card}}}
\newcommand\cinv{^{\langle-1\rangle}}
\newcommand\contr{{\downharpoonleft}}
\newcommand\Cpx{{\mathbf C}}
\newcommand\Dc{{\mathcal{D}}}
\newcommand\Dco{\overline{\Dc}}
\newcommand\diag{\text{\rm diag}}
\newcommand\Ec{{\mathcal{E}}}
\newcommand\eps{\epsilon}
\newcommand\Fh{{\widehat F}}
\newcommand\Gc{{\mathcal{G}}}
\newcommand\Ic{{\mathcal{I}}}
\newcommand\id{{\operatorname{id}}}
\newcommand\IP{{\operatorname{IP}}}
\newcommand\Lc{{\mathcal{L}}}
\newcommand\ldeg{\operatorname{ldeg}}
\newcommand\lenc{\le_{\operatorname{nc}}}
\newcommand\mfs{{\operatorname{mfs}}}
\newcommand\Mul{{\operatorname{Mul}}}
\newcommand\Nats{{\mathbf N}}
\newcommand\NC{\operatorname{NC}}
\newcommand\NCL{\operatorname{NCL}}
\newcommand\NCLo{\overline{\NCL}}
\newcommand\otdt{\otimes\cdots\otimes}
\newcommand\Phit{{\widetilde\Phi}}
\newcommand\piar{{\overset{\to}\pi}}
\newcommand\pichk{{\check\pi}}
\newcommand\pihat{{\hat\pi}}
\newcommand\pit{{\tilde\pi}}
\newcommand\polar{{\operatorname{polar}}}
\newcommand\Qc{{\mathcal{Q}}}
\newcommand\qt{{\tilde q}}
\newcommand\Rt{{\widetilde R}}
\newcommand\Sc{{\mathcal{S}}}
\newcommand\sigmachk{{\check\sigma}}
\newcommand\sigmahat{{\hat\sigma}}
\newcommand\sigmat{{\tilde\sigma}}
\newcommand\St{{\widetilde S}}
\newcommand\stackdots{\raisebox{0.2\height}[0.8\height]{\vdots}} 
\newcommand\Sym{{\operatorname{Sym}}}
\newcommand\SymMul{{\operatorname{SymMul}}}
\newcommand\tauar{{\overset{\to}\tau}}

\newcommand\tauhat{{\hat\tau}}
\newcommand\taut{{\tilde\tau}}
\newcommand\tdt{\times\cdots\times}
\newcommand\Tt{{\widetilde T}}

\begin{document}

\pagestyle{myheadings}

 \title{Multilinear function series and transforms in free probability theory}

 \author{Kenneth J.\ Dykema}

\address{Mathematisches Institut, Westf\"alische Wilhelms--Universit\"at M\"uns\-ter,
Ein\-stein\-str.\ 62, 48149 M\"unster, Germany;
{\rm permanent address:} Department of Mathematics, Texas A\&M University,
College Station, TX 77843-3368, USA}
\email{kdykema@math.tamu.edu}

 \thanks{Supported in part by NSF grant DMS--0300336 and by the Alexander von Humboldt Foundation.}

 \date{\today}

 \begin{abstract}
The algebra $\Mul[[B]]$ of formal multilinear function series over an algebra $B$
and its quotient $\SymMul[[B]]$ are introduced,
as well as corresponding operations of formal composition.
In the setting of $\Mul[[B]]$, the unsymmetrized R-- and T--transforms
of random variables in $B$--valued noncommutative probability spaces
are introduced.
These satisfy properties analogous to the usual R-- and T--transforms, (the latter being just the
reciprocal of the S--transform), but describe all moments of a random variable,
not only the symmetric moments.
The partially ordered set of noncrossing linked partitions is introduced and is used
to prove properties of the unsymmetrized T--transform.
 \end{abstract}

 \maketitle


\markboth{\sc multilinear function series}{\sc multilinear function series}

\section{Introduction and summary of results}

Recall that,
given a field $K$, the algebra $K[[X]]$
of formal power series in one variable over $K$ contains the algebra $K[X]$
of polynomials in one variable over $K$.
Also composition is defined for formal power series with zero constant term,
extending the usual notion of composition in $K[X]$.
If we take $K$ to be the field of complex numbers, $\Cpx$, then the algebra of germs
of analytic functions at zero is embedded in $\Cpx[[X]]$ by
associating to each analytic function its power series expansion at zero.
Also, for such series, composition in $\Cpx[[X]]$ corresponds to composition of analytic
functions.

Let $B$ be a Banach algebra.
If $F$ is a $B$--valued analytic function defined in a neighborhood of $0$ in $B$,
then $F$ has a series expansion of the form
\begin{equation}\label{eq:Fanalytic}
F(b)=F(0)+\sum_{k=1}^\infty F_k(b,\ldots,b),
\end{equation}
where $F_k$ is a symmetric multilinear function from the $k$--fold product $B\tdt B$
to $B$, (see~\cite{HP}).

Given an arbitrary algebra $B$ over a field $K$, we define the algebra $\Mul[[B]]$
of formal multilinear function series.
An element of $\Mul[[B]]$ is a sequence
\[
\alpha=(\alpha_0,\alpha_1,\alpha_2,\ldots)
\]
where $\alpha_0$ belongs to the unitization of $B$ and where for every $k\ge1$,
$\alpha_k$ is a multilinear function from the $k$--fold product $B\tdt B$ to $B$.
Then $\Mul[[B]]$ is a vector space over $K$ with the obvious operations.
Given $\alpha,\beta\in\Mul[[B]]$, we define the formal product $\alpha\beta$ by setting
\[
(\alpha\beta)_n(b_1,\ldots,b_n)=\sum_{k=0}^n\alpha_k(b_1,\ldots,b_k)\beta_{n-k}(b_{k+1},\ldots,b_n).
\]
This makes $\Mul[[B]]$ an algebra, and the unitization of $B$ is embedded in $\Mul[[B]]$ as
a unital subalgebra by
\[
b\mapsto(b,0,0,\ldots).
\]
We denote by $1=(1,0,0,\ldots)$
the identity element of $\Mul[[B]]$ with respect to formal product.

If $\beta_0=0$, then we define the formal composition $\alpha\circ\beta$
by setting $(\alpha\circ\beta)_0=\alpha_0$ and, for $n\ge1$,
\begin{multline*}
(\alpha\circ\beta)_n(b_1,\ldots,b_n)= \\
=\sum_{k=1}^n\sum_{\substack{p_1,\ldots,p_k\ge1 \\ p_1+\cdots+p_k=n}}\begin{aligned}[t]
\alpha_k(\beta_{p_1}(b_1,\ldots,b_{p_1}),\beta_{p_2}(b_{q_2+1},\ldots,b_{q_2+p_2}),\ldots&,\\
\beta_{p_k}(b_{q_k+1},\ldots,b_{q_k+p_k}&)),\end{aligned}
\end{multline*}
where $q_j=p_1+\cdots+p_{j-1}$.
Then
\begin{equation}\label{eq:I}
I:=(0,\id_B,0,0,\ldots)
\end{equation}
is a left and right identity element for formal composition;
we will also show that formal composition is associative and behaves as expected
with formal products and sums.

The name ``formal multilinear function series'' and the operations described above
are inspired by the series expansions~\eqref{eq:Fanalytic} of analytic functions.
However, the name should not be understood literally:
when $\dim(B)\ge2$ there is no series (formal or otherwise) to sum
that contains all the information of an element of $\Mul[[B]]$.
But, if the field $K$ has characteristic zero,
then $\Mul[[B]]$ has a quotient, $\SymMul[[B]]$, whose elements are at least formally
equivallent to functions obtained by summing series.
Specifically, elements of $\SymMul[[B]]$ are the sequences
\[
\alpha=(\alpha_0,\alpha_1,\alpha_2,\ldots)
\]
where $\alpha_k$ is a symmetric multilinear function, and $\SymMul[[B]]$ has
operations of symmetric
product and symmetric composition.
In the case that $B$ is a Banach algebra over $\Cpx$,
the map $F\mapsto(F(0),F_1,F_2,\ldots)$ sending a $B$--valued analytic function $F$
defined near $0$ in $B$ to the sequence of terms in its series
expansion~\eqref{eq:Fanalytic} is an embedding of the algebra of germs of $B$--valued
analytic functions at $0$ in $B$ into $\SymMul[[B]]$, and the composition of such functions
(taking value $0$ at $0$) corresponds to symmetric composition in $\SymMul[[B]]$.
Thus, $\SymMul[[B]]$ can be viewed as the algebra of formal series obtained directly from the
theory of analytic functions by disregarding all questions of convergence (and allowing the
algebra $B$ to be over an arbitrary field), while $\Mul[[B]]$ can be thought of as
obtained from $\SymMul[[B]]$
by allowing the multilinear functions to be nonsymmetric and redefining
the operations appropriately.

\smallskip
The above algebraic constructions
are useful for describing operations on free random variables.
Freeness is a notion that applies to elements of noncommutative probability spaces,
and is analogous to independence in classical probability theory.
It was discovered by Voiculescu~\cite{V85}
and is the fundamental idea in free probability theory.
See, for example, the book~\cite{VDN} for definitions, early results and references.
We now briefly recall some important concepts, while specifying notation that we will use throughout
the paper.
Given a unital algebra $B$ (over a general field $K$),
a {\em $B$--valued noncommutative probability space}
is a pair $(A,E)$ where $A$ is a unital algebra containing $B$ as a unital subalgebra
and where $E:A\to B$ is a conditional expectation, namely a linear idempotent map restricting
to the identity on $B$ and satisfying the conditional expectation property
\[
E(b_1ab_2)=b_1E(a)b_2\qquad(b_1,b_2\in B,\,a\in A).
\]
Elements of $A$ are called random variables, and the {\em distribution series} of $a\in A$ is
$\Phit_a\in\Mul[[B]]$ given by $\Phit_{a,0}=E(a)$ and
\begin{equation}\label{eq:Phita}
\Phit_{a,n}(b_1,\ldots,b_n)=E(ab_1ab_2\cdots ab_na).
\end{equation}
Thus, the distribution series of $a$ encodes all the moments of $a$.
Recall that subalgebras $A_i\subseteq A$ for $i\in I$, each containing $B$, are said to be {\em free}
if $E(a_1a_2\cdots a_n)=0$ whenever $n\in\Nats$, $a_j\in A_{i_j}\cap\ker E$
and $i_1\ne i_2,\,\ldots,i_{n-1}\ne i_n$.
Random variables $a_i$ ($i\in I$) are said to be free if the algebras $\alg(\{a_i\}\cup B)$
generated by them together with $B$ are free.
The mixed moments
\[
E(a_{i_0}b_1a_{i_1}b_2\cdots a_{i_{n-1}}b_na_{i_n})\qquad(i_j\in I,\,b_j\in B)
\]
of free random variables $a_i$ are determined by the freeness relation in terms of the moments
of the $a_i$ individually.

Suppose $a_1$ and $a_2$ are free random variables.
It is of particular interest to describe the distribution series $\Phit_{a_1+a_2}$ and $\Phit_{a_1a_2}$
of their sum and product in terms of $\Phit_{a_1}$ and $\Phit_{a_2}$.
We will refer to these as the problems of additive and multiplicative free convolution,
respectively.

These problems were solved in the case $B=\Cpx$, (taken as a complex algebra), by Voiculescu
in~\cite{V86} and~\cite{V87}.
Note that in this case, $\Phit_a$ is just a formal power series with complex coefficients.
Voiculescu's solutions may be described as follows:
Voiculescu defined the R--transform $R_f$ and the S--transform $S_f$
of a formal power series $f$, (the latter in the case that $f$ has nonzero constant coefficient).
These may be defined as the unique formal power series
satisfying
\begin{align}
(z+z^2f(z))\cinv&=z(1+z\,R_f(z))^{-1}, \label{eq:Rf} \\
(zf(z))\cinv&=z(1+z)^{-1}\,S_f(z) \label{eq:Sf}
\end{align}
where the left--hand--sides of~\eqref{eq:Rf} and~\eqref{eq:Sf} are the inverses
with respect to composition
of the formal power series $z+z^2f(z)$ and $zf(z)$, respectively.
The series $f$ can be recovered from each of $R_f$ and $S_f$.
Voiculescu showed
\begin{align}
R_{\Phit_{x+y}}&=R_{\Phit_x}+R_{\Phit_y} \label{eq:VR} \\
S_{\Phit_{xy}}&=S_{\Phit_x}\,S_{\Phit_y} \label{eq:VS}
\end{align}
when $x$ and $y$ are free (over $B=\Cpx$), where in~\eqref{eq:VS} it is assumed $E(x)$ and $E(y)$
are nonzero.

For a general algebra $B$,
the problem of additive free convolution was shown by Voiculescu~\cite{V95} to have a solution;
he constructed certain canonical random variables and showed them to be additive
under hypothesis of freeness.
The solution to additive free convolution was given a beautiful combinatorial description by Speicher
in~\cite{Sp1} and~\cite{Sp},
using the lattice of noncrossing partitions.
Voiculescu~\cite{V95} also showed that, if $B$ is a unital Banach algebra, then for random variables
in a $B$--valued Banach noncommutative probability space,
the solution
to the symmetrized version of the problem of additive free convolution can be expressed
in terms of $B$--valued analytic functions.
Specifically, suppose $(A,E)$ is a B--valued Banach noncommutative probability space, namely
$A$ is a unital Banach algebra containing $B$ as a closed, unital subalgebra
and $E:A\to B$ is a bounded conditional expectation.
For $a\in A$, let
\[
\Phi_a(b)=\sum_{n=0}^\infty E(a(ba)^n)=E(a(1-ba)^{-1}),\qquad(b\in B,\,\|b\|<\|a\|^{-1}).
\]
Then $\Phi_a$ is a $B$--valued analytic function in the indicated neighborhood of $0$ in $B$,
and is the symmetric analogue of $\Phit_a$ described in~\eqref{eq:Phita},
because from the $n$th term of the Taylor expansion of $\Phi_a$ 
we recover the symmetric moments
\begin{equation}\label{eq:symmom}
\Phi_{a,n}(b_1,\ldots,b_n)=\frac1{n!}\sum_{\sigma\in S_n}
E(ab_{\sigma(1)}ab_{\sigma(2)}\cdots ab_{\sigma(n)}a)
\end{equation}
of $a$.
(Note that when $B=\Cpx$, then $\Phi_a$ and $\Phit_a$ coincide, where we equate a complex--valued
analytic function with its power series expansion.)
Given a $B$--valued analytic function $F$ defined in a neighborhood of $0$ in $B$,
Voiculsecu's R--transfrom of $F$ is the unique (germ of a)
$B$--valued analytic function $R_F$ defined
in some neighborhood of $0$ in $B$ and satisfying
\begin{equation}\label{eq:CR}
C_F\cinv(b)=(1+b\,R_F(b))^{-1}b
\end{equation}
where $C_F\cinv$ is the inverse with respect to composition
of the function
\begin{equation}\label{eq:CF}
C_F(b)=b+b\,F(b)\,b.
\end{equation}
Then for $x,y\in A$ free with respect to $B$, Voiculescu's result yields
\begin{equation}\label{eq:Rplus}
R_{\Phi_{x+y}}=R_{\Phi_x}+R_{\Phi_y}.
\end{equation}

The first solution of the problem of multiplicative free convolution beyond the case $B=\Cpx$
was obtained by Aagaard~\cite{Aa}.
He treated the case of $B$ a commutative Banach algebra, and showed that the $S$--transform
is multiplicative for free random variables in a $B$--valued Banach noncommutative probability space.
In~\cite{D}, we treated the case of $B$ a general Banach algebra.
For a $B$--valued analytic function $F$ defined in a neighborhood of $0$ in $B$
and satisfying that $F(0)$ is invertible in $B$,
we define the S--transform of $F$ to be the unique
(germ of a) $B$--valued analytic function $S_F$ defined in some neighborhood
of $0$ in $B$ and satisfying
\begin{equation}\label{eq:PsiS}
D_F\cinv(b)=b(1+b)^{-1}\,S_F(b),
\end{equation}
where $D_F\cinv$ is the inverse with respect to composition
of the function 
\begin{equation}\label{eq:DF}
D_F(b)=bF(b).
\end{equation}
We showed that the S--transform satifies a twisted multiplicativity property with respect to freeness,
namely
\begin{equation}\label{eq:Stmult}
S_{\Phi_{xy}}(b)=S_{\Phi_y}(b)\,S_{\Phi_x}(\bt)
\end{equation}
where
\[
\bt=(S_{\Phi_y}(b))^{-1}\,b\,S_{\Phi_y}(b)
\]
is $b$ conjugated by the inverse of $S_{\Phi_y}(b)$,
assuming $x$ and $y$ are free in $(A,E)$ and $E(x)$ and $E(y)$ are invertible.

To summarize, the R--transform and S--transform of a symmetric moment generating function $\Phi_a$
of a random variable $a$ over a Banach algebra $B$
are equivalent to $\Phi_a$ via function theoretic relations~\eqref{eq:CR} and~\eqref{eq:CF},
respectively~\eqref{eq:PsiS} and~\eqref{eq:DF};
the additivity and twisted multiplicativity results~\eqref{eq:Rplus} and~\eqref{eq:Stmult}
allow computation of the generating functions for symmetric moments of the sum $x+y$ and,
respectively, product $xy$ of free random variables $x$ and $y$
in terms of those for the original variables $x$ and $y$.

As we will see, it turns out to be more natural to consider the reciprocal
of the S--transform, which we call the T--transform.
Thus, from~\eqref{eq:PsiS} we have the following relation that defines the T--transform:
\begin{equation}\label{eq:Tdef}
D_F\cinv(b)\,T_F(b)\,(1+b)=b.
\end{equation}
The twisted multiplicativity property~\eqref{eq:Stmult} for free random variables $x$ and $y$ becomes
\begin{equation}\label{eq:Ttmult}
T_{\Phi_{xy}}(b)=T_{\Phi_x}(T_{\Phi_y}(b)\,b\,(T_{\Phi_y}(b))^{-1})\,T_{\Phi_y}(b).
\end{equation}
Substituting $D_F(b)=b\,F(b)$ for $b$ in~\eqref{eq:Tdef}, we get
\[
b\,T_F(b\,F(b))\,(1+b\,F(b))=b\,F(b)
\]
and equating Taylor expansions we deduce
\begin{equation}\label{eq:Tdefeq}
T_F(b\,F(b))\,(1+b\,F(b))=F(b),
\end{equation}
which may be taken as an equivalent definition of $T_F$.

Our main applications of the algebra of formal multilinear function series, and our motivation
for introducing it, are the {\em unsymmetrized R--transform} and {\em unsymmetrized T--transform}, denoted
$\Rt$ and $\Tt$.
These do for the distribution series $\Phit_a$ what the R--transform and T--transform
do for the generating functions $\Phi_a$, but in the context of $\Mul[[B]]$ rather than
$B$--valued analytic functions.
The unsymmetrized R-- and T--transforms generalize the usual ones in two directions.
On the one hand, all moments of a random variable $a$ are treated, not only the symmetric ones.
On the other hand, $\Mul[[B]]$ and the unsymmetrized R-- and T--transforms
are purely algebraic: neither norms nor boundedness is required and
$B$ can be an arbitrary algebra over an arbitrary field of any characteristic.
Amazingly, the equations defining the unsymmetrized R-- and T--transforms and 
giving their properties for free random variables are exactly the same as those we have seen
for the usual R-- and T--transforms, but now interpreted in $\Mul[[B]]$.
(Moreover, the properties of the usual R--and T--tranforms follow from the 
properties we will prove for the unsymmetrized versions by passing to the quotient $\SymMul[[B]]$.)
Thus, given $\beta\in\Mul[[B]]$, its unsymmetrized R--transform is the unique element $\Rt_\beta$
of $\Mul[[B]]$ satisfying
\begin{equation}\label{eq:Rt}
\alpha\cinv=(1+I\,\Rt_\beta)\,I,
\end{equation}
where $\alpha\cinv$ is the inverse with respect to formal composition
in $\Mul[[B]]$ of
\begin{equation}\label{eq:Rtalpha}
\alpha=I+I\,\beta\,I,
\end{equation}
where all products in~\eqref{eq:Rt} and~\eqref{eq:Rtalpha} are formal products in $\Mul[[B]]$
and where $I$ is the identity element with respect to formal composition as given in~\eqref{eq:I}.
We immediately see that~\eqref{eq:Rt} and~\eqref{eq:Rtalpha}
are the same as~\eqref{eq:CR} and~\eqref{eq:CF},
but in this different context.
Furthermore, we will show that if $(A,E)$ is a $B$--valued noncommutative probability space and
if $x,y\in A$ are free, then
\begin{equation}\label{eq:Rtplus}
\Rt_{\Phit_{x+y}}=\Rt_{\Phit_x}+\Rt_{\Phit_y}.
\end{equation}
If $\beta=(\beta_0,\beta_1,\ldots)\in\Mul[[B]]$ with $\beta_0$ an invertible element of the
unitization of $B$, then the unsymmetrized T--transform of $\beta$ is the 
unique element $\Tt_\beta$ of $\Mul[[B]]$ satisfying
\begin{equation}\label{eq:Tt}
(\Tt_\beta\circ(I\,\beta))\,(1+I\,\beta)=\beta,
\end{equation}
and we will show that for $x,y\in A$ free with $E(x)$ and $E(y)$ invertible, we have
\begin{equation}\label{eq:Tttmult}
\Tt_{\Phit_{xy}}=(\Tt_{\Phit_x}\circ(\Tt_{\Phit_y}\,I\,\Tt_{\Phit_y}^{-1}))\,\Tt_{\Phit_y},
\end{equation}
where in~\eqref{eq:Tt} and~\eqref{eq:Tttmult}, all products are formal products in $\Mul[[B]]$,
the inverse is with respect to formal product,
$\circ$ denotes formal composition in $\Mul[[B]]$ and $I$ is as in~\eqref{eq:I}.
Here we see that~\eqref{eq:Tt} and~\eqref{eq:Tttmult} are
the same as~\eqref{eq:Tdefeq} and~\eqref{eq:Ttmult}, respectively, in this new context.

Ultimately, formulas such as~\eqref{eq:Rtplus} and~\eqref{eq:Tttmult} are combinatorial
in nature, even if the succinct notation is reminiscent of and inspired by function theory.
It is no coincidence that our proof of~\eqref{eq:Rtplus}
looks very similar to Speicher's work in~\cite{Sp1} and~\cite{Sp}.
In fact, our results on the unsymmetrized R--transform can be viewed as duplicating some of
Speicher's results, but from a dual perspective.
Our proofs of properties of the unsymmetrized R--transform use combinatorial techniques quite similar
to Speicher's, namely the lattice of noncrossing partitions.
To prove our results on the unsymmetrized T--transform,
we introduce the combinatorial structure of
{\em noncrossing linked partitions}.
A noncrossing linked partition of $\{1,\ldots,n\}$ can be viewed as a noncrossing partition
with possibly some links of a restricted nature drawn between certain blocks of the partition.
We denote the set of all noncrossing linked partitions of $\{1,\ldots,n\}$ by $\NCL(n)$;
it has at least two natural partial orderings, and it contains the order--embedded
lattice $\NC(n)$ of all noncrossing partitions of $\{1,\ldots,n\}$.
Even in the case $B=\Cpx$, our results give a new interpretation of the $n$th coefficient
of the T--transform of $\Phi_a$ as a sum over $\NCL(n+1)$ of certain products in the moments of $a$,
analogous to Speicher's expression of the coefficients of the R--transform as
a similar sum over noncrossing partitions.
We use this to show that the cardinality of $\NCL(n)$ is the large Schr\"oder number $r_n$
and we obtain an expression for the large Schr\"oder numbers as a sum over noncrossing partitions
of a product of Catalan numbers.

The contents of the rest of this paper are outlined below.
In~\S\ref{sec:Mul}, we introduce the algebra of formal multilinear function series
over an algebra $B$ and the operation of formal composition and prove basic properties.
In~\S\ref{sec:SymMul}, when the underlying field is of characteristic zero, we consider the
quotient $\SymMul[[B]]$ obtained by symmetrizing, and describe its relationship
to polynomials and, in the case of a Banach algebra, to $B$--valued analytic functions.
In~\S\ref{sec:crv}, we consider an algebraic version of the Fock space (over $B$)
considered in~\cite{D}, and we construct random variables acting on this space and having
arbitrary distribution.
In~\S\ref{sec:NCL}, we introduce the partially ordered set of noncrossing linked partitions.
In~\S\ref{sec:Rtrans}, we introduce the unsymmetrized R--transform
and prove~\eqref{eq:Rt} and additivity~\eqref{eq:Rtplus}.
In~\S\ref{sec:Strans}, we introduce the unsymmetrized T--transform
and prove~\eqref{eq:Tt} and twisted multiplicativity~\eqref{eq:Tttmult}.
In~\S\ref{sec:scalar}, in the case that $B=\Cpx$,
we give the formula for moments in terms of the coefficients in the T--transform
and we apply this to prove that the cardinality of $\NCL(n+1)$ is the large
Schr\"oder number $r_n$.

\medskip
\noindent{\em Acknowledgements.} The author thanks
Joachim Cuntz, Siegfried Echterhoff and the Mathematics Institute of the
Westf\"alische Wilhelms--Universit\"at
M\"unster for their generous hospitality during the author's year--long visit, when
this research was conducted.
The author also thanks Catherine Yan for helpful conversations.

\section{Formal multilinear function series}
\label{sec:Mul}

\begin{defi}\label{def:fmfs}
Let $B$ be an algebra over a field $K$.
By $\Lc(B)$ we denote the set of all linear mappings from $B$ to itself,
and for $n\ge1$ we denote by $\Lc_n(B)$ the set of all $n$--multilinear mappings
\[
\alpha_n:\underset{n\text{ times}}{\underbrace{B\tdt B}}\to B,
\]
Moreover, we set $\Bt$ equal to $B$ if $B$ has an identity element and to the unitization
of $B$ otherwise, and we let $1$ denote the identity element of $\Bt$.
A {\em formal multilinear function series} over $B$ is a sequence
$\alpha=(\alpha_0,\alpha_1,\ldots)$, where $\alpha_0\in\Bt$ and $\alpha_n\in\Lc_n(B)$
for $n\ge1$.
We let $\Mul[[B]]$ denote the set of all formal multilinear function series over $B$.
\end{defi}

Given $\alpha\in\Mul[[B]]$, we imagine a function from $B$ to $\Bt$ sending $b$ to
\begin{equation}\label{eq:Bseries}
\alpha_0+\alpha_1(b)+\alpha_2(b,b)+\alpha_3(b,b,b)+\cdots,
\end{equation}
although no such function need exist.
However, the formal series~\eqref{eq:Bseries} serves to inspire the
operations of multiplication
and composition of formal multilinear function series, which we now define.

\begin{defi}\label{def:addmult}
Let $\alpha,\beta\in\Mul[[B]]$.
Then the {\em sum} $\alpha+\beta$ and {\em formal product}
$\alpha\beta$ are the elements of $\Mul[[B]]$ defined by
\[
(\alpha+\beta)_n=\alpha_n+\beta_n
\]
and $(\alpha\beta)_0=\alpha_0\beta_0$, while for $n\ge1$
\[
(\alpha\beta)_n(b_1,\ldots,b_n)=\sum_{k=0}^n\alpha_k(b_1,\ldots,b_k)\beta_{n-k}(b_{k+1},\ldots,b_n).
\]
If $\beta_0=0$, then the {\em formal composition} $\alpha\circ\beta\in\Mul[[B]]$ is defined
by $(\alpha\circ\beta)_0=\alpha_0$, while for $n\ge1$
\begin{multline}\label{eq:acircb}
(\alpha\circ\beta)_n(b_1,\ldots,b_n)= \\
=\sum_{k=1}^n\sum_{\substack{p_1,\ldots,p_k\ge1 \\ p_1+\cdots+p_k=n}}\begin{aligned}[t]
\alpha_k(\beta_{p_1}(b_1,\ldots,b_{p_1}),\beta_{p_2}(b_{q_2+1},\ldots,b_{q_2+p_2}),\ldots&,\\
\beta_{p_k}(b_{q_k+1},\ldots,b_{q_k+p_k}&)),\end{aligned}
\end{multline}
where $q_j=p_1+\cdots+p_{j-1}$
\end{defi}

With the obvious scalar multiplication and addition defined
above, $\Mul[[B]]$ is clearly a vector space over $K$.
It is, moreover, a $B,B$--bimodule in the obvious way.
The following proposition gives some of basic algebraic properties of $\Mul[[B]]$;
in particular, it is an algebra over $K$ and formal composition behaves as expected.

\begin{prop}\label{prop:Mul}
Let $\alpha,\beta,\gamma\in\Mul[[B]]$.
Then
\renewcommand{\labelenumi}{(\roman{enumi})}
\begin{enumerate}
\item the formal product is associative: $(\alpha\beta)\gamma=\alpha(\beta\gamma)$;
\item distributivity holds: $(\alpha+\beta)\gamma=\alpha\gamma+\beta\gamma$
and $\alpha(\beta+\gamma)=\alpha\beta+\alpha\gamma$;
\item the element $1=(1,0,0,\ldots)\in\Mul[[B]]$ is a multiplicative identity element;
\item an element
$\alpha=(\alpha_0,\alpha_1,\ldots)\in\Mul[[B]]$ has a multiplicative inverse
if and only if $\alpha_0$ is an invertible element of $\Bt$;
\item formal composition is associative: if $\beta_0=0$ and $\gamma_0=0$, then
$(\alpha\circ\beta)\circ\gamma=\alpha\circ(\beta\circ\gamma)$;
\item if $\gamma_0=0$, then $(\alpha+\beta)\circ\gamma=\alpha\circ\gamma+\beta\circ\gamma$
and $(\alpha\beta)\circ\gamma=(\alpha\circ\gamma)(\beta\circ\gamma)$;
\item the element $I=(0,\id_B,0,0,\ldots)\in\Mul[[B]]$
is an identity element for the operation of composition;
\item an element
$\alpha=(0,\alpha_1,\alpha_2,\ldots)\in\Mul[[B]]$ has an inverse with respect to
formal composition, denoted $\alpha\cinv$,
if and only if $\alpha_1$ is an invertible element of $\Lc(B)$.
\end{enumerate}
\end{prop}
\begin{proof}
Parts~(i)--(iii) are routine verifications.

For (iv), clearly if $\alpha$ is invertible, then so is $\alpha_0$.
If $\alpha_0$ is invertible, then
we find a left multiplicative inverse $\alpha'$ for $\alpha$
by taking $\alpha'_0=(\alpha_0)^{-1}$ and defining recursively, for $n\ge1$,
\[
\alpha'_n(b_1,\ldots,b_n)=-\sum_{k=0}^{n-1}\alpha'_k(b_1,\ldots,b_k)
\alpha_{n-k}(b_{k+1},\ldots,b_n)(\alpha_0)^{-1}.
\]
A right inverse for $\alpha$ is found similarly.

Parts~(v)--(vii) are routine (though sometimes tedious) verifications.

For~(viii), if $\alpha\cinv$ is an inverse with respect to formal composition for $\alpha$,
then $(\alpha\cinv)_1$ is an inverse for $\alpha_1$ in $\Lc(B)$.
On the other hand, if $\alpha_1$ has an inverse $(\alpha_1)\cinv$ in $\Lc(B)$, then we solve
$\alpha'\circ\alpha=I$
by setting $\alpha'_1=(\alpha_1)\cinv$ and recursively defning $\alpha'_n$ for $n\ge2$ by
\begin{align*}
\alpha_n'(&\alpha_1(b_1),\ldots,\alpha_1(b_n))= \\
&=-\sum_{k=1}^{n-1}\sum_{\substack{p_1,\ldots,p_k\ge1 \\ p_1+\cdots+p_k=n}}
\alpha'_k(\alpha_{p_1}(b_{q_1+1},\ldots,b_{q_1+p_1}),\ldots,\alpha_{p_k}(b_{q_k+1},\ldots,b_{q_k+p_k})),
\end{align*}
with $q_j=p_1+\cdots+p_{j-1}$,
and we similarly solve $\alpha\circ\alpha''=I$
by setting $\alpha''_1=(\alpha_1)\cinv$ and recursively defning
\begin{align}
\alpha&_n''(b_1,\ldots,b_n)= \label{eq:alpha''} \\
&=-(\alpha_1)\cinv\bigg(
\sum_{k=2}^n\sum_{\substack{p_1,\ldots,p_k\ge1 \\ p_1+\cdots+p_k=n}}
\alpha_k(\alpha''_{p_1}(b_{q_1+1},\ldots,b_{q_1+p_1}),\ldots,\alpha''_{p_k}(b_{q_k+1},\ldots,b_{q_k+p_k}))
\bigg). \notag
\end{align}
\end{proof}

\begin{defi}
If $\alpha\in\Mul[[B]]$ is nonzero, then
the {\em lower degree} of $\alpha$ is the least $n\in\{0,1,2,\ldots\}$
such that $\alpha_n\ne0$, and is denoted $\ldeg(\alpha)$.
\end{defi}

\begin{prop}
Let $\alpha,\beta\in\Mul[[B]]$.
Then
\[
\ldeg(\alpha\beta)\ge\ldeg(\alpha)+\ldeg(\beta).
\]
If $\beta_0=0$, then
\[
\ldeg(\alpha\circ\beta)\ge\ldeg(\alpha)\ldeg(\beta).
\]
\end{prop}

If $\alpha^{(k)}\in\Mul[[B]]$ for $k\in\Nats$
and if $\lim_{k\to\infty}\ldeg(\alpha^{(k)})=\infty$,
then $\sum_{k=1}^\infty\alpha^{(k)}$ is well defined as an element of $\Mul[[B]]$,
and we have obvious identities such as
\[
\beta(\sum_{k=1}^\infty\alpha^{(k)})=\sum_{k=1}^\infty\beta\alpha^{(k)}.
\]
Now one sees immediately that the formula for the geomtric series holds also in $\Mul[[B]]$.
\begin{prop}\label{prop:1+a}
Let $\alpha\in\Mul[[B]]$ and suppose $\alpha_0=0$.
Then
\[
(1-\alpha)^{-1}=1+\sum_{k=1}^\infty\alpha^k.
\]
\end{prop}

\section{Symmetric formal multilinear function series}
\label{sec:SymMul}

In this section, $B$ will be an algebra over a field of characteristic zero.

\begin{defi}
$\Sym$ will denote the usual {\em symmetrization operator} on multilinear functions.
Thus, for $\alpha_n\in\Lc_n(B)$, $\Sym\,\alpha_n\in\Lc_n(B)$ is given by
\[
(\Sym\,\alpha_n)(b_1,\ldots,b_n)=\frac1{n!}\sum_{\sigma\in S_n}
\alpha_n(b_{\sigma(1)},\ldots,b_{\sigma(n)}),
\]
where $S_n$ is the group of permutations of $\{1,\ldots,n\}$.
Of course, $\Sym$ is an idempotent operator.
If $\alpha=(\alpha_0,\alpha_1,\alpha_2,\ldots)\in\Mul[[B]]$, then we let
\[
\Sym\,\alpha=(\alpha_0,\alpha_1,\Sym\,\alpha_2,\Sym\,\alpha_3,\ldots),
\]
and the set of all {\em symmetric formal multilinear function series} is
\[
\SymMul[[B]]=\{\Sym\,\alpha\mid\alpha\in\Mul[[B]]\}.
\]
\end{defi}

Clearly $\Sym$ preserves the vector space operations on $\Mul[[B]]$, but
the product and composition of two symmetric elements in $\Mul[[B]]$ need not
be symmetric.

\begin{defi}
Let $\alpha,\beta\in\Mul[[B]]$.
Then the {\em symmetric product} of $\alpha$ and $\beta$ is 
\[
\alpha\cdot_s\beta=\Sym(\alpha\beta).
\]
If $\beta_0=0$, then the {\em symmetric composition} of $\alpha$ and $\beta$ is
\[
\alpha\circ_s\beta=\Sym(\alpha\circ\beta).
\]
\end{defi}

\begin{lemma}
Let $\alpha,\beta\in\Mul[[B]]$.
Then
\begin{equation}\label{eq:Symprod}
\alpha\cdot_s\beta=(\Sym\,\alpha)\cdot_s\beta=\alpha\cdot_s(\Sym\,\beta).
\end{equation}
If $\beta_0=0$, then
\begin{equation}\label{eq:Symcomp}
\alpha\circ_s\beta=\alpha\circ_s(\Sym\,\beta)=(\Sym\,\alpha)\circ_s\beta.
\end{equation}
\end{lemma}
\begin{proof}
We have
\[
(\alpha\cdot_s\beta)_n(b_1,\ldots,b_n)=\frac1{n!}\sum_{k=0}^n\sum_{\sigma\in S_n}
\alpha_k(b_{\sigma(1)},\ldots,b_{\sigma(k)})\beta_{n-k}(b_{\sigma(k+1)},\ldots,b_{\sigma(n)}),
\]
while
\begin{multline*}
((\Sym\,\alpha)\cdot_s\beta)_n(b_1,\ldots,b_n)= \\
=\frac1{n!}\sum_{k=0}^n\sum_{\sigma\in S_n}
\frac1{k!}\sum_{\rho\in S_k}
\alpha_k(b_{\sigma\circ\rho(1)},\ldots,b_{\sigma\circ\rho(k)})
\beta_{n-k}(b_{\sigma(k+1)},\ldots,b_{\sigma(n)}).
\end{multline*}
But the map $S_n\times S_k\to S_n$ given by $(\sigma,\rho)\mapsto\sigma\circ(\rho\oplus\id_{n-k})$
is $k!$--to--1.
This shows the first equality in~\eqref{eq:Symprod}.
The second follows similarly.

We have
\[
(\alpha\circ_s\beta)_n(b_1,\ldots,b_n)
=\frac1{n!}\sum_{\sigma\in S_n}\sum_{k=1}^n\sum_{\substack{p_1,\ldots,p_k\ge1 \\ p_1+\cdots+p_k=n}}
\begin{aligned}[t]
\alpha_k(&\beta_{p_1}(b_{\sigma(q_1+1)},\ldots,b_{\sigma(q_1+p_1)}), \\
&\ldots,\beta_{p_k}(b_{\sigma(q_k+1)},\ldots,b_{\sigma(q_k+p_k)})),
\end{aligned}
\]
where $q_j=p_1+\cdots+p_{j-1}$,
while
\begin{align*}
(\alpha&\circ_s(\Sym\,\beta))_n(b_1,\ldots,b_n)= \\
&=\frac1{n!}\sum_{\sigma\in S_n}\sum_{k=1}^n\sum_{\substack{p_1,\ldots,p_k\ge1 \\ p_1+\cdots+p_k=n}} \\
&\qquad\frac1{p_1!p_2!\cdots p_k!}\sum_{\tau_1\in S_{p_1},\ldots,\tau_k\in S_{p_k}}
\begin{aligned}[t]
\alpha_k(&\beta_{p_1}(b_{\sigma(q_1+\tau_1(1))},\ldots,b_{\sigma(q_1+\tau_1(p_1))}), \\
&\ldots,\beta_{p_k}(b_{\sigma(q_k+\tau_k(1))},\ldots,b_{\sigma(q_k+\tau_k(p_k))})).\end{aligned}
\end{align*}
But the map $S_n\times S_{p_1}\tdt S_{p_k}\to S_n$ given by
$(\sigma,\tau_1,\ldots,\tau_k)\mapsto\sigma\circ(\tau_1\oplus\cdots\oplus\tau_k)$ is
$p_1!p_2!\cdots p_k!$--to--1.
This shows the first equality in~\eqref{eq:Symcomp}.
On the other hand, we have
\begin{multline*}
((\Sym\,\alpha)\circ_s\beta)_n(b_1,\ldots,b_n)= \\
=\frac1{n!}\sum_{\sigma\in S_n}\sum_{k=1}^n\frac1{k!}\sum_{\tau\in S_k}
\sum_{\substack{p_1,\ldots,p_k\ge1 \\ p_1+\cdots+p_k=n}}
\begin{aligned}[t]
\alpha_k(&\beta_{p_{\tau(1)}}(b_{\sigma(q_{\tau(1)}+1)},\ldots,b_{\sigma(q_{\tau(1)}+p_{\tau(1)})}), \\
&\ldots,\beta_{p_{\tau(k)}}(b_{\sigma(q_{\tau(k)}+1)},\ldots,b_{\sigma(q_{\tau(k)}+p_{\tau(k)})}))
\end{aligned}
\end{multline*}
and the corresponding map $S_n\times S_k\to S_n$ is $k!$--to--1.
This shows the second equality in~\eqref{eq:Symcomp}.
\end{proof}

The following are pretty much analogues for the operations on
$\SymMul[[B]]$
of the properties listed in Proposition~\ref{prop:Mul}.
In particular, $\SymMul[[B]]$ is a unital algebra and symmetric composition
behaves naturally.

\begin{prop}
Let $\alpha,\beta,\gamma\in\SymMul[[B]]$.
Then
\renewcommand{\labelenumi}{(\roman{enumi})}
\begin{enumerate}
\item the symmetric product is associative:
\[
(\alpha\cdot_s\beta)\cdot_s\gamma=\alpha\cdot_s(\beta\cdot_s\gamma);
\]
\item distributivity holds:
\begin{align*}
(\alpha+\beta)\cdot_s\gamma&=\alpha\cdot_s\gamma+\beta\cdot_s\gamma \\
\alpha\cdot_s(\beta+\gamma)&=\alpha\cdot_s\beta+\alpha\cdot_s\gamma;
\end{align*}
\item the element $1=(1,0,0,\ldots)\in\SymMul[[B]]$ is an identity element
for the symmetric product operation;
\item if an element
$\alpha=(\alpha_0,\alpha_1,\ldots)\in\SymMul[[B]]$ has a multiplicative inverse $\alpha^{-1}$
in $\Mul[[B]]$, then $\Sym(\alpha^{-1})$ is the inverse of $\alpha$ for the symmetric product;
\item symmetric composition is associative: if $\beta_0=0$ and $\gamma_0=0$, then
\[
(\alpha\circ_s\beta)\circ_s\gamma=\alpha\circ_s(\beta\circ_s\gamma);
\]
\item if $\gamma_0=0$, then
\[
(\alpha+\beta)\circ_s\gamma=\alpha\circ_s\gamma+\beta\circ_s\gamma
\]
and
\[
(\alpha\cdot_s\beta)\circ_s\gamma=(\alpha\circ_s\gamma)\cdot_s(\beta\circ_s\gamma);
\]
\item the element $I=(0,\id_B,0,0,\ldots)\in\SymMul[[B]]$
is an identity element for the operation of symmetric composition;
\item if an element
$\alpha=(0,\alpha_1,\alpha_2,\ldots)\in\SymMul[[B]]$ has an inverse with respect to composition,
denoted $\alpha\cinv$,
then $\Sym(\alpha\cinv)$ is the inverse of $\alpha$ with respect to symmetric composition.
\end{enumerate}
\end{prop}
\begin{proof}
Parts (i)--(iii) follow from~\eqref{eq:Symprod} and
the analogous properties in $\Mul[[B]]$.
For (iv), using~\eqref{eq:Symprod} we have
\[
\alpha\cdot_s\Sym(\alpha^{-1})=\Sym(\alpha\alpha^{-1})=1=
\Sym(\alpha^{-1}\alpha)=\Sym(\alpha^{-1})\cdot_s\alpha.
\]
Parts (v)--(vii) follow from~\eqref{eq:Symcomp} and
the analogous properties in $\Mul[[B]]$.
For (viii), using~\eqref{eq:Symcomp} we have
\[
\alpha\circ_s\Sym(\alpha\cinv)=\Sym(\alpha\circ\alpha\cinv)=I
=\Sym(\alpha\cinv\circ\alpha)=\Sym(\alpha\cinv)\circ_s\alpha.
\]
\end{proof}

We now assume that $B$ is an algebra over the field of complex numbers.
The algebra $\SymMul[[B]]$, together with the (partially defined)
operation of symmetric composition, can be viewed as analogues of
formal power series over $B$.
When $B$ is a Banach algebra, $\SymMul[[B]]$ contains the algebra of germs
of $B$--valued analytic functions at $0$,
and symmtric composition in $\SymMul[[B]]$ corresponds to usual composition of analytic functions.
We finish this section with a brief explanation of this assertion, which is based on
Hille and Phillips~\cite{HP}, where proofs can be found; (see Chapters~III and~XXVI).

A {\em polynomial} of degree no more than $m$ over $B$ is a function $P:B\to B$ such that
for all $b,h\in B$ there are $c_0,\ldots,c_m\in B$ such that for all $\lambda\in\Cpx$,
\[
P(b+\lambda h)=\sum_{k=0}^mc_k\lambda^k.
\]
The polynomial $P$ is {\em homogeneous} of degree $m$ if
\[
P(\lambda h)=\lambda^mP(h)\qquad(h\in B,\,\lambda\in\Cpx).
\]
The {\em polar part} of such a homogeneous polynomial $P$ is the symmetric $m$--multilinear
function $\polar(P)\in\Lc_m(B)$ given by
\begin{equation}\label{eq:pol}
\begin{split}
\polar(P)&(h_1,\ldots,h_m)= \\
&=\frac1{m!}\Delta^m_{h_1,\ldots,h_m}P(0)
:=\frac1{m!}\sum_{k=1}^m(-1)^{m-k}\sum_{1\le i(1)<\cdots<i(k)\le m}
P(\sum_{j=1}^kh_{i(j)})
\end{split}
\end{equation}
and $P$ is recovered as the diagonal $b\mapsto\polar(P)(b,b,\ldots,b)$ of $\polar(P)$.
Moreover, given $\alpha_m\in\Lc_m(B)$, its diagonal 
\[
\diag(\alpha_m):b\mapsto\alpha_m(b,b,\ldots,b)
\]
is a homogeneous polynomial of degree $m$, whose polar part is equal to $\alpha_m$.
In this way, $\Lc_m(B)$ is identified with the set of all homogenous polynomials of degree
$m$ over $B$.
Any polynomial $Q$ of degree no more than $m$ over $B$ can be uniquely written
as $Q=Q_0+\cdots+Q_m$, where $Q_k$ is a homogeneous polynomial of degree $k$.
We define the {\em multilinear function series} of $Q$ to be
\[
\mfs(Q)=(Q_0,\polar(Q_1),\polar(Q_2),\ldots,\polar(Q_m),0,0,\ldots)\in\SymMul[[B]].
\]
We have that $\mfs$ is a bijection from the set of all polynomials over $B$ onto the set
of all elements of $\SymMul[[B]]$ having finite support.
The following proposition shows that this bijection preserves the algebra structure and composition.
Thus, $\SymMul[[B]]$ is related to the polynomials over $B$ as the set of formal power series $\Cpx[[z]]$
in one variable with complex coefficients is related to the polynomial algebra $\Cpx[X]$.
\begin{prop}\label{prop:PQ}
If $P$ and $Q$ are polynomials over $B$, then
\begin{align}
\mfs(P+Q)&=\mfs(P)+\mfs(Q) \label{eq:P+Q} \\
\mfs(PQ)&=\mfs(P)\cdot_s\mfs(Q). \label{eq:PQ}
\end{align}
Moreover, if $Q(0)=0$, then
\begin{equation}\label{eq:PcircQ}
\mfs(P\circ Q)=\mfs(P)\circ_s\mfs(Q).
\end{equation}
\end{prop}
\begin{proof}
Additivity~\eqref{eq:P+Q} is clear.
Using~\eqref{eq:P+Q}, in order to show~\eqref{eq:PQ} we may without loss of generality
assume that $P$ is homogeneous of degree $m$ and $Q$ is homogenous of degree $n$.
Then $PQ$ is homogeneous of degree $m+n$.
Clearly, $(\mfs(P)\cdot_s\mfs(Q))_k=0$ if $k\ne m+n$.
Let $\alpha_m=\polar(P)\in\Lc_m(B)$ and $\beta_n=\polar(Q)\in\Lc_n(B)$.
Then
\[
\begin{split}
(\mfs(P)\,\cdot_s\,\mfs(Q))_{m+n}(b,\ldots,b)&=(\mfs(P)\mfs(Q))_{m+n}(b,\ldots,b)= \\
&=\alpha_m(b,\ldots,b)\beta_n(b,\ldots,b)=P(b)Q(b),
\end{split}
\]
and we conclude that~\eqref{eq:PQ} holds.

In order to show~\eqref{eq:PcircQ}, using~\eqref{eq:P+Q}
we may without loss of generality assume that $P$ is homogenous of degree $m$.
Let $n$ be the degree of $Q$.
Let $\alpha_m=\polar(P)$, so that
\[
\alpha=(0,\ldots,0,\alpha_m,0,0,\ldots)=\mfs(P)
\]
and let
\[
\beta=(0,\beta_1,\beta_2,\ldots,\beta_n,0,0,\ldots)=\mfs(Q).
\]
Then for $k\ge1$ and $b\in B$,
\begin{equation}\label{eq:acbk}
\begin{split}
(\alpha\circ_s\beta)_k(b,\ldots,b)&=(\alpha\circ\beta)_k(b,\ldots,b)= \\[1ex]
&=\sum_{\substack{p_1,\ldots,p_m\in\{1,\ldots,n\} \\ p_1+\cdots+p_m=k}}
\alpha_m(\beta_{p_1}(b,\ldots,b),\ldots,\beta_{p_m}(b,\ldots,b)).
\end{split}
\end{equation}
On the other hand,
\[
\begin{split}
(P\circ Q)(b)
&=\alpha_m\big(\sum_{p_1=1}^n\beta_{p_1}(b,\ldots,b),\ldots,\sum_{p_m=1}^n\beta_{p_m}(b,\ldots,b)\big)= \\[1ex]
&=\sum_{p_1,\ldots,p_m\in\{1,\ldots,n\}}
\alpha_m(\beta_{p_1}(b,\ldots,b),\ldots,\beta_{p_m}(b,\ldots,b)).
\end{split}
\]
But for any given $p_1,\ldots,p_m\in\{1,\ldots,n\}$,
the map
\[
b\mapsto\alpha_m(\beta_{p_1}(b,\ldots,b),\ldots,\beta_{p_m}(b,\ldots,b))
\]
is easily seen to be a homogenous polynomial of degree $p_1+\cdots+p_m$.
So selecting the degree $k$ homogenous part $(P\circ Q)_k$ of $P\circ Q$ and comparing to~\eqref{eq:acbk},
we find $(P\circ Q)_k=\diag((\alpha\circ_s\beta)_k)$, and~\eqref{eq:PcircQ} is proved.
\end{proof}

Now suppose that $B$ is a Banach algebra over $\Cpx$.
Recall that an $m$--multilinear function $\alpha_m\in\Lc_m(B)$ is bounded if
there is a positive real number $M\ge0$ such that 
\[
\alpha_m(h_1,\ldots,h_m)\le M\|h_1\|\,\|h_2\|\,\cdots\,\|h_m\|
\]
for all $h_1,\ldots,h_m\in B$, and then the norm $\|\alpha_m\|$ is the least such constant $M$.
Let $\alpha=(\alpha_0,\alpha_1,\ldots)\in\SymMul[[B]]$ and suppose $\alpha_m$ is bounded for all $m$.
We define the {\em radius of convergence} of $\alpha$ to be 
\[
r(\alpha)=(\limsup_{n\to\infty}\|\alpha_n\|^{1/n})^{-1}\in[0,+\infty].
\]
Given any $\alpha=(\alpha_0,\alpha_1,\ldots)\in\SymMul[[B]]$,
we say $\alpha$ has {\em positive radius of convergence}
if all terms $\alpha_m$ are bounded and $r(\alpha)>0$.

If $F$ is a $B$--valued analytic function defined in a neighborhood of $0$ in $B$,
then the Taylor series expansion of $F$ is
\[
F(b)=F(0)+\sum_{n=1}^\infty F_n(b,\ldots,b),
\]
where $F_n\in\Lc_n(B)$ is the bounded symmetric $n$--multilinear function such that
$n!F_n(h_1,\ldots,h_n)$ is the $n$--fold Fr\'echet derivative of $F$ at $0$ with increments
$h_1,\ldots,h_n$.
We define the {\em multilinear function series} of $F$ to be
\[
\mfs(F)=(F(0),F_1,F_2,\ldots)\in\SymMul[[B]].
\]
Then $\mfs(F)$ has positive radius of convergence.
On the other hand,
if $\alpha\in\SymMul[[B]]$ is such that $\alpha_n$ is bounded for all $n$
and $r(\alpha)>0$, then
\begin{equation}\label{eq:diagalpha}
\diag(\alpha)(b):=\alpha_0+\sum_{n=1}^\infty\alpha_n(b,\ldots,b)
\end{equation}
is defined for all $b\in B$ with $\|b\|<r(\alpha)$ and $\diag(\alpha)$ is $B$--valued analytic
on this neighborhood of $0$ in $B$, such that $\mfs(\diag(\alpha))=\alpha$.
Therefore, $\mfs$ is a bijection from the set of germs of $B$--valued analytic functions
defined in neighborhoods of $0$ in $B$ onto the set of all $\alpha\in\SymMul[[B]]$
having positive radius of convergence.
The following proposition shows that $\mfs$ is an algebra homorphism and preserves composition.
\begin{prop}\label{prop:mfsanalytic}
Let $F$ and $G$ be $B$--valued analytic functions, each defined on a neighborhood of $0$ in $B$.
Then
\begin{align}
\mfs(F+G)&=\mfs(F)+\mfs(G) \label{eq:mfsF+G} \\
\mfs(FG)&=\mfs(F)\cdot_s\mfs(G). \label{eq:mfsFG}
\end{align}
Moreover, if $G(0)=0$, then
\begin{equation}
\mfs(F\circ G)=\mfs(F)\circ_s\mfs(G). \label{eq:mfsFcircG}
\end{equation}
Furthermore, if $\mfs(F)$ is invertible with respect to symmetric multiplication or, respectively,
symmetric composition in $\SymMul[[B]]$, then $F$ is invertible
with respect to multiplication or, respectively, composition as a $B$--valued analytic function.
\end{prop}
\begin{proof}
Equations~\eqref{eq:mfsF+G}, \eqref{eq:mfsFG} and~\eqref{eq:mfsFcircG}
follow using the convergence properties of Taylor expansions of analytic functions
and the properties of Fr\'echet derivatives.
If $\mfs(F)$ is invertible with respect to symmtric multiplication, then $F(0)$ is an invertible
element of $B$, so $F(b)$ is a invertible element of $B$ for all $b$ in some neighborhood of $0$ in $B$
and the multiplicative inverse $F^{-1}$ is analytic in this neighborhood of $B$.
If $\mfs(F)$ is invertible with respect to symmetric composition, then $F(0)=0$
and the first Fr\'echet derivative $F_1$ of $F$ is invertible as a map from $B$ to $B$.
By the open mapping theorem, $F_1$ has bounded inverse, so by the usual inverse function
theorem for maps between Banach spaces, it follows that the restriction of
$F$ to some neighborhood of $0$ in $B$ is invertible with respect to composition
and $F\cinv$ is analytic.
\end{proof}

\section{Canonical random variables}
\label{sec:crv}

In~\cite{D}, we constructed a Banach space analogue of full Fock space and of certain
annihilation and creation operators that were used to model
$B$--valued random variables, where $B$ is a Banach algebra.
In this section, we consider a purely algebraic version of this construction.
These results and proofs are similar to those of~\cite{D}.

Let $B$ be a unital algebra over the complex numbers and let $I$ be a nonempty set.
(For convenience of labeling, we usually assume that the index set $I$ contains as elements $1$ and $2$.)

Let 
\[
D=\bigoplus_{i\in I}B=\{d:I\to B\mid d(i)=0\text{ for all but finitely many }i\}
\]
be the algebraic direct sum of $|I|$ copies of $B$.
We equip $D$ with the obvious left action of $B$, given by $(bd)(i)=b(d(i))$.
Let
\[
\Qc=B\Omega\oplus\bigoplus_{k=1}^\infty(D^{\otimes k})\otimes B
\]
be the algebraic direct sum of tensor products over $\Cpx$,
where $B\Omega$ denotes a copy of $B$ whose identity element is distinguished
with the name $\Omega$.
Let $P:\Qc\to B\Omega=B$ be the idempotent that
is the identity operator on $B\Omega$ and that
sends every $D^{\otimes k}\otimes B$ to zero.
Let $\Lc(\Qc)$ denote the algebra of $\Cpx$--linear operators from $\Qc$ to $\Qc$.
We have the left and right actions of $B$ on $\Qc$,
\begin{align*}
\lambda&:B\to\Lc(\Qc) \\
\rho&:B\to\Lc(\Qc),
\end{align*}
given by
\begin{alignat*}{2}
\lambda(b):b_0\Omega&\mapsto bb_0\Omega,\qquad
&\lambda(b):d_1\otdt d_k\otimes b_0&\mapsto(bd_1)\otimes d_2\otdt d_k\otimes b_0 \\
\rho(b):b_0\Omega&\mapsto b_0b\Omega,\qquad
&\rho(b):d_1\otdt d_k\otimes b_0&\mapsto d_1\otdt d_k\otimes (b_0b).
\end{alignat*}
Note that $\lambda(b_1)$ and $\rho(b_2)$ commute for all $b_1,b_2\in B$.
Let $\Ec:\Lc(\Qc)\to B$ be $\Ec(X)=P(X\Omega)$.
Then $\Ec\circ\lambda=\id_B$.
We think of $B$ as embedded in $\Lc(\Qc)$ via $\lambda$, acting on the left,
and we will often omit to write $\lambda$.
The short proof of the next proposition is exactly like that of~\cite[Prop.\ 3.2]{D}.
\begin{prop}
The restriction of $\Ec$ to the commutant $\Lc(\Qc)\cap\rho(B)'$ of $\rho(B)$
satisfies the conditional expectation property
\[
\Ec(b_1Xb_2)=b_1\Ec(X)b_2,\qquad (X\in\Lc(\Qc)\cap\rho(B)',\,b_1,b_2\in B).
\]
\end{prop}

For $i\in I$, consider the {\em creation operator} $L_i\in\Lc(\Qc)$ defined by
\[
L_i:b_0\Omega\mapsto\delta_i\otimes b_0,\qquad
L_i:d_1\otdt d_k\otimes b_0\mapsto \delta_i\otimes d_1\otdt d_k\otimes b_0
\]
where $\delta_i\in D$ is the characteristic function of $i$.

Recall that for $n\ge1$, $\Lc_n(B)$ denotes the set of all $n$--multilinear maps $B\tdt B\to B$
and $\Lc_0$ denotes $B$.
As in~\cite{D}, for $n\ge1$, $\alpha_n\in\Lc_n(B)$ and $i\in I$, we define $V_{i,n}(\alpha_n)$
and $W_{i,n}(\alpha_n)$ in $\Lc(\Qc)$ by
by
\begin{align*}
V_{i,n}(\alpha_n)(b_0\Omega)&=0 \\[1ex]
V_{i,n}(\alpha_n)(d_1\otdt d_k\otimes b_0)&=
\begin{cases}
0,&k<n \\
\alpha_n(d_1(i),\ldots,d_n(i))b_0\Omega,&k=n \\
\alpha_n(d_1(i),\ldots,d_n(i))d_{n+1}\otdt d_k\otimes b_0,&k>n
\end{cases}
\end{align*}
and
\begin{align*}
W_{i,n}(\alpha_n)&(b_0\Omega)=0 \displaybreak[2] \\[1ex]
W_{i,n}(\alpha_n)&(d_1\otdt d_k\otimes b_0)= \\
&=\begin{cases}
0,&k<n \\
\alpha_n(d_1(i),\ldots,d_n(i))\delta_i\otimes b_0,&k=n \\
\alpha_n(d_1(i),\ldots,d_n(i))\delta_i\otimes d_{n+1}\otdt d_k\otimes b_0,&k>n.
\end{cases}
\end{align*}
For $n=0$ and $\alpha_0\in B$, we let
\[
V_{i,0}(\alpha_0)=\alpha_0\qquad W_{i,0}(\alpha_0)=\alpha_0L_i.
\]
Thus, $V_{i,n}(\alpha_n)$ is an $n$--fold annihilation operator,
while $W_{i,n}(\alpha_n)$ is an $n$--fold annihilation combined with a creation,
though because of commutativity issues, $W_{i,n}(\alpha_n)$ cannot
in general be written as a product of $L_i$ and $V_{i,n}(\alpha_n)$.

All the formulas in~\cite[Lemma 3.3]{D} continue to hold, and the
proof of~\cite[Prop.\ 3.4]{D} carries over to give the following freeness result.
\begin{lemma}\label{lem:free}
For $i\in I$ let $A_i\subseteq\Lc(\Qc)\cap\rho(B)'$ be the subalgebra
generated by
\[
\lambda(B)\cup\{L_i\}\cup\{V_{i,n}(\alpha_n)\mid n\in\Nats,\,\alpha_n\in\Lc_n(B)\}
\cup\{W_{i,n}(\alpha_n)\mid n\in\Nats,\,\alpha_n\in\Lc_n(B)\}.
\]
Then the family $(A_i)_{i\in I}$ is free with respect to $\Ec$.
\end{lemma}

Since every element of $\Qc$ belongs to 
\[
\Qc_N:=B\Omega\oplus\bigoplus_{k=1}^N(D^{\otimes k})\otimes B
\]
for some $N\in\Nats$, given $\alpha\in\Mul[[B]]$,
we may define elements of $\Lc(\Qc)$
\begin{align*}
V_i(\alpha)&=\sum_{n=0}^\infty V_{i,n}(\alpha_n) \\
W_i(\alpha)&=\sum_{n=0}^\infty W_{i,n}(\alpha_n)
\end{align*}
in the obvious way, namely, if $f\in \Qc_N$, then
\[
V_i(\alpha)f=\sum_{n=0}^NV_{i,n}(\alpha_n)f\qquad
W_i(\alpha)f=\sum_{n=0}^NW_{i,n}(\alpha_n)f.
\]
The next result is an extension of Lemma~\ref{lem:free} that follows
immediately from it.
\begin{prop}\label{prop:free}
For $i\in I$ let $\Afr_i\subseteq\Lc(\Qc)\cap\rho(B)'$ be the subalgebra
generated by
\[
\lambda(B)\cup\{L_i\}\cup\{V_i(\alpha)\mid\alpha\in\Mul[[B]]\}
\cup\{W_i(\alpha)\mid\alpha\in\Mul[[B]]\}.
\]
Then the family $(\Afr_i)_{i\in I}$ is free with respect to $\Ec$.
\end{prop}

\begin{defi}
Given a $B$--valued noncommutative probability space $(A,E)$ and $a\in A$,
the {\em distribution series} of $a$ is the element $\Phit_a\in\Mul[[B]]$ given by
\begin{align*}
\Phit_{a,0}&=E(a), \\
\Phit_{a,n}(b_1,\ldots,b_n)&=E(ab_1ab_2\cdots ab_na),\quad(n\ge1).
\end{align*}
\end{defi}

We now construct random variables in $(\Lc(\Qc),\Ec)$ of specific forms
having arbitrary distribution series.
This result is an analogue of Propositions~3.6 and~4.4 of~\cite{D}.
\begin{prop}\label{prop:can}
Let $\beta\in\Mul[[B]]$.
Then there is a unique $\alpha\in\Mul[[B]]$ such that if
\begin{equation}\label{eq:X}
X=L_i+V_i(\alpha)\in\Lc(\Qc)
\end{equation}
for some $i\in I$,
then $\Phit_X=\beta$.
If $\beta_0$ is an invertible element of $B$, then there is a unique $\alpha\in\Mul[[B]]$ such that
if
\begin{equation}\label{eq:Yf}
Y=W_i(\alpha)+V_i(\alpha)\in\Lc(\Qc)
\end{equation}
for some $i\in I$,
then $\Phit_Y=\beta$.
\end{prop}
\begin{proof}
Let $\alpha=(\alpha_0,\alpha_1,\ldots)\in\Mul[[B]]$
and let $X$ and $Y$ be as in~\eqref{eq:X} and~\eqref{eq:Yf}.
For $N\ge0$, let
\begin{align*}
X_N&=L_i+\sum_{n=0}^NV_{i,n}(\alpha_n) \\
Y_N&=\sum_{n=0}^N(W_{i,n}(\alpha_n)+V_{i,n}(\alpha_n)).
\end{align*}
Then $E(X)=\alpha_0$
and, for $b_1,\ldots,b_N\in B$,
\begin{align*}
E(Xb_1Xb_2\cdots Xb_NX)&=E(X_Nb_1X_Nb_2\cdots X_Nb_NX_N)= \\
&=\alpha_N(b_1,\ldots,b_N)+
E(X_{N-1}b_1X_{N-1}b_2\cdots X_{N-1}b_{N-1}X_{N-1}).
\end{align*}
Thus, taking $\alpha_0=\beta_0$, we may choose $\alpha_n$ recursively so that $\Phit_X=\beta$.

We also have $E(Y)=\alpha_0$
and, for $b_1,\ldots,b_N\in B$,
\begin{align*}
E(Yb_1Yb_2\cdots Yb_NY)&=E(Y_Nb_1Y_Nb_2\cdots Y_Nb_NY_N)= \\
&=\alpha_N(b_1\alpha_0,\ldots,b_N\alpha_0)+
E(Y_{N-1}b_1Y_{N-1}b_2\cdots Y_{N-1}b_{N-1}Y_{N-1}).
\end{align*}
Thus, if $\beta_0$ is invertible, then taking $\alpha_0=\beta_0$,
we may choose $\alpha_n$ recursively so that $\Phit_Y=\beta$.
\end{proof}

\begin{defi}
The element $X$ constructed in Proposition~\ref{prop:can} is called
the {\em additive canonical random variable}
with distribution series $\beta$,
and the element $Y$ is called
the {\em multiplicative canonical random variable}
with distribution series $\beta$.
\end{defi}

\section{Noncrossing linked partitions}
\label{sec:NCL}

\begin{defi}
Let $n\in\Nats$ and $E,F\subseteq\{1,\ldots,n\}$.
We say that $E$ and $F$ are {\em crossing} if there exist $i_1,i_2\in E$ and
$j_1,j_2\in F$ with $i_1<j_1<i_2<j_2$.
Otherwise, we say $E$ and $F$ are {\em noncrossing}.
We say that $E$ and $F$ are {\em nearly disjoint} if for every $i\in E\cap F$,
one of the following holds:
\renewcommand{\labelenumi}{(\alph{enumi})}
\begin{enumerate}
\item $i=\min(E)$, $|E|>1$ and $i\ne\min(F)$,
\item $i\ne\min(E)$, $i=\min(F)$ and $|F|>1$.
\end{enumerate}
\end{defi}

\begin{defi}\label{def:le}
Let $\pi$ and $\sigma$ be sets of subsets of $\{1,\ldots,n\}$.
We write $\pi\le\sigma$ if for every $E\in\pi$
there is $E'\in\sigma$ with $E\subseteq E'$.
\end{defi}

As usual, a {\em noncrossing partition} of $\{1,\ldots,n\}$ is a partition of $\{1,\ldots,n\}$
into disjoint subsets, any two distinct elements of which are noncrossing.
The set of all noncrossing partitions of $\{1,\ldots,n\}$ is a lattice under $\le$ and is denoted $\NC(n)$.
As usual, the largest element $\{\{1,\ldots,n\}\}$ of $\NC(n)$ is denoted $1_n$,
and the smallest element $\{\{1\},\ldots,\{n\}\}$ is written $0_n$.

An {\em interval partition} of $\{1,\ldots,n\}$ is a partition of $\{1,\ldots,n\}$ into disjoint subsets
that are intervals.
We write $\IP(n)$ for the set of all interval partitions of $\{1,\ldots,n\}$.

\begin{defi}\label{def:ncnp}
Let $n\in\Nats$.
A set $\pi$ of nonempty subsets of $\{1,\ldots,n\}$ is said to be a
{\em noncrossing linked partition} of $\{1,\ldots,n\}$
if
the union of $\pi$ is all of $\{1,\ldots,n\}$ and any two distinct elements of $\pi$
are noncrossing and nearly disjoint.
We let $\NCL(n)$ denote the set of all noncrossing linked partitions of $\{1,\ldots,n\}$.
\end{defi}

For $\pi\in\NCL(n)$, we will sometimes refer to an element of $\pi$ as a {\em block} of $\pi$.

\begin{remark}\label{rmk:ncnp}
Let $\pi\in\NCL(n)$.
\renewcommand{\labelenumi}{(\roman{enumi})}
\begin{enumerate}
\item  Any given element $\ell$ of $\{1,\ldots,n\}$ belongs to either
exactly one or exactly two blocks of $\pi$;
we will say $\ell$ is {\em singly} or {\em doubly covered} by $\pi$, accordingly.
\item The elements $1$ and $n$ are singly covered by $\pi$.
\item Any two distinct elements $E$ and $F$ of $\pi$ have at most one element in common.
Moreover, if $E\cap F\ne\emptyset$, then both $E$ and $F$ have at least two elements.
\item The lattice $\NC(n)$ of noncrossing partitions of $\{1,\ldots,n\}$ is a subset of $\NCL(n)$
and consists of those elements $\pi$ of $\NCL(n)$ such that all elements of $\{1,\ldots,n\}$
are singly covered by $\pi$, (i.e.\ the elements of $\NCL(n)$ that are actual partitions
of $\{1,\ldots,n\}$ into disjoint subsets).
\end{enumerate}
\end{remark}

All the elements of $\NCL(3)$ and $\NCL(4)$ are listed in Tables~\ref{tab:ncnp3} and~\ref{tab:ncnp4},
respectively, while selected elements of $\NCL(5)$ are listed in Table~\ref{tab:ncnp5}.
\begin{table}[b]
\caption{The elements $\pi$ of $\NCL(3)$ and graphical representations $\Gc_\pi$.}
\label{tab:ncnp3}
\begin{tabular}{r|l||r|l}
$\pi$\hspace*{2ex} & \hspace*{2ex}$\Gc_\pi$ & $\pi$\hspace*{2ex} & \hspace*{2ex}$\Gc_\pi$ \\ \hline\hline
$(1,2,3)$ 
&
\begin{picture}(40,18)(-3,0) 
 \multiput(0,0)(12,0){3}{\circle*{3}}
 \multiput(0,0)(12,0){3}{\line(0,1){12}}
 \put(0,12){\line(1,0){24}}
\end{picture}
& \hspace*{1.5ex} $(1,2)(3)$
& 
\begin{picture}(30,18)(-3,0) 
 \multiput(0,0)(12,0){3}{\circle*{3}}
 \multiput(0,0)(12,0){2}{\line(0,1){12}}
 \put(0,12){\line(1,0){12}}
\end{picture} \\ \hline
$(1,2)(2,3)$
& 
\begin{picture}(40,18)(-3,0) 
 \multiput(0,0)(12,0){3}{\circle*{3}}
 \multiput(0,0)(12,0){3}{\line(0,1){12}}
 \put(12,0){\line(1,1){12}}
 \put(0,12){\line(1,0){12}}
\end{picture}
& \hspace*{1.5ex} $(1,3)(2)$
&
\begin{picture}(30,18)(-3,0) 
 \multiput(0,0)(12,0){3}{\circle*{3}}
 \multiput(0,0)(24,0){2}{\line(0,1){12}}
 \put(0,12){\line(1,0){24}}
\end{picture} \\ \hline
$(1)(2,3)$
&
\begin{picture}(40,18)(-3,0) 
 \multiput(0,0)(12,0){3}{\circle*{3}}
 \multiput(12,0)(12,0){2}{\line(0,1){12}}
 \put(12,12){\line(1,0){12}}
\end{picture}
& \hspace*{1.5ex} $(1)(2)(3)$
&
\begin{picture}(30,18)(-3,0) 
 \multiput(0,0)(12,0){3}{\circle*{3}}
\end{picture}
\end{tabular}
\end{table}
\begin{table}[b]
\caption{The elements $\pi$ of $\NCL(4)$ and graphical representations $\Gc_\pi$.}
\label{tab:ncnp4}
\begin{tabular}{r|l||r|l}
$\pi$\hspace*{2ex} & \hspace*{2ex}$\Gc_\pi$ & $\pi$\hspace*{2ex} & \hspace*{2ex}$\Gc_\pi$ \\ \hline\hline
$(1,2,3,4)$ &
\begin{picture}(60,18)(-3,0) 
 \multiput(0,0)(12,0){4}{\circle*{3}}
 \multiput(0,0)(12,0){4}{\line(0,1){12}}
 \put(0,12){\line(1,0){36}}
\end{picture} 
& \hspace*{1.5ex} $(1,2,3)(4) $ &
\begin{picture}(46,18)(-3,0) 
 \multiput(0,0)(12,0){4}{\circle*{3}}
 \multiput(0,0)(12,0){3}{\line(0,1){12}}
 \put(0,12){\line(1,0){24}}
\end{picture} \\ \hline
$(1,2,3)(3,4) $ &
\begin{picture}(60,18)(-3,0) 
 \multiput(0,0)(12,0){4}{\circle*{3}}
 \multiput(0,0)(12,0){4}{\line(0,1){12}}
 \put(0,12){\line(1,0){24}}
 \put(24,0){\line(1,1){12}}
\end{picture}
& \hspace*{1.5ex}$(1,2,4)(3) $ &
\begin{picture}(46,18)(-3,0) 
 \multiput(0,0)(12,0){4}{\circle*{3}}
 \multiput(0,0)(12,0){2}{\line(0,1){12}}
 \put(36,0){\line(0,1){12}}
 \put(0,12){\line(1,0){36}}
\end{picture} \\ \hline
$(1,2,4)(2,3) $ &
\begin{picture}(60,18)(-3,0) 
 \multiput(0,0)(12,0){4}{\circle*{3}}
 \multiput(0,0)(12,0){2}{\line(0,1){12}}
 \put(36,0){\line(0,1){12}}
 \put(0,12){\line(1,0){36}}
 \put(12,0){\line(3,2){12}}
 \put(24,0){\line(0,1){8}}
\end{picture}
& \hspace*{1.5ex} $(1,3,4)(2) $ &
\begin{picture}(46,18)(-3,0) 
 \multiput(0,0)(12,0){4}{\circle*{3}}
 \put(0,0){\line(0,1){12}}
 \multiput(24,0)(12,0){2}{\line(0,1){12}}
 \put(0,12){\line(1,0){36}}
\end{picture} \\ \hline
$(1,2)(3,4)$ &
\begin{picture}(60,18)(-3,0) 
 \multiput(0,0)(12,0){4}{\circle*{3}}
 \multiput(0,0)(12,0){4}{\line(0,1){12}}
 \multiput(0,12)(24,0){2}{\line(1,0){12}}
\end{picture} 
& \hspace*{1.5ex} $(1,2)(2,3,4) $ &
\begin{picture}(46,18)(-3,0) 
 \multiput(0,0)(12,0){4}{\circle*{3}}
 \multiput(0,0)(12,0){4}{\line(0,1){12}}
 \put(12,0){\line(1,1){12}}
 \multiput(0,12)(24,0){2}{\line(1,0){12}}
\end{picture} \\ \hline
$(1,4)(2,3)$ &
\begin{picture}(60,18)(-3,0) 
 \multiput(0,0)(12,0){4}{\circle*{3}}
 \multiput(0,0)(36,0){2}{\line(0,1){12}}
 \multiput(12,0)(12,0){2}{\line(0,1){8}}
 \put(0,12){\line(1,0){36}}
 \put(12,8){\line(1,0){12}}
\end{picture} 
& \hspace*{1.5ex} $(1,2)(3)(4)$ &
\begin{picture}(46,18)(-3,0) 
 \multiput(0,0)(12,0){4}{\circle*{3}}
 \multiput(0,0)(12,0){2}{\line(0,1){12}}
 \put(0,12){\line(1,0){12}}
\end{picture} \\ \hline
$(1,2)(2,3)(4)$ &
\begin{picture}(60,18)(-3,0) 
 \multiput(0,0)(12,0){4}{\circle*{3}}
 \multiput(0,0)(12,0){3}{\line(0,1){12}}
 \put(0,12){\line(1,0){12}}
 \put(12,0){\line(1,1){12}}
\end{picture}
& \hspace*{1.5ex} $(1,2)(2,3)(3,4)$ &
\begin{picture}(46,18)(-3,0) 
 \multiput(0,0)(12,0){4}{\circle*{3}}
 \multiput(0,0)(12,0){4}{\line(0,1){12}}
 \put(0,12){\line(1,0){12}}
 \multiput(12,0)(12,0){2}{\line(1,1){12}}
\end{picture} \\ \hline
$(1,2)(2,4)(3)$ &
\begin{picture}(60,18)(-3,0) 
 \multiput(0,0)(12,0){4}{\circle*{3}}
 \multiput(0,0)(12,0){2}{\line(0,1){12}}
 \multiput(0,12)(24,0){2}{\line(1,0){12}}
 \put(36,0){\line(0,1){12}}
 \put(12,0){\line(1,1){12}}
\end{picture}
& \hspace*{1.5ex} $(1,3)(2)(4)$ &
\begin{picture}(46,18)(-3,0) 
 \multiput(0,0)(12,0){4}{\circle*{3}}
 \multiput(0,0)(24,0){2}{\line(0,1){12}}
 \put(0,12){\line(1,0){24}}
\end{picture} \\ \hline 
$(1,3)(2)(3,4) $ &
\begin{picture}(60,18)(-3,0) 
 \multiput(0,0)(12,0){4}{\circle*{3}}
 \multiput(0,0)(24,0){2}{\line(0,1){12}}
 \put(24,0){\line(1,1){12}}
 \put(36,0){\line(0,1){12}}
 \put(0,12){\line(1,0){24}}
\end{picture}
& \hspace*{1.5ex} $(1,4)(2)(3)$ &
\begin{picture}(46,18)(-3,0) 
 \multiput(0,0)(12,0){4}{\circle*{3}}
 \multiput(0,0)(36,0){2}{\line(0,1){12}}
 \put(0,12){\line(1,0){36}}
\end{picture}  \\ \hline
$(1)(2,3,4)$ &
\begin{picture}(60,18)(-3,0) 
 \multiput(0,0)(12,0){4}{\circle*{3}}
 \multiput(12,0)(12,0){3}{\line(0,1){12}}
 \put(12,12){\line(1,0){24}}
\end{picture} 
& \hspace*{1.5ex} $(1)(2,3)(4) $ &
\begin{picture}(46,18)(-3,0) 
 \multiput(0,0)(12,0){4}{\circle*{3}}
 \multiput(12,0)(12,0){2}{\line(0,1){12}}
 \put(12,12){\line(1,0){12}}
\end{picture} \\ \hline
$(1)(2,3)(3,4) $ &
\begin{picture}(60,18)(-3,0) 
 \multiput(0,0)(12,0){4}{\circle*{3}}
 \multiput(12,0)(12,0){3}{\line(0,1){12}}
 \put(12,12){\line(1,0){12}}
 \put(24,0){\line(1,1){12}}
\end{picture}
& \hspace*{1.5ex} $(1)(2,4)(3) $ &
\begin{picture}(46,18)(-3,0) 
 \multiput(0,0)(12,0){4}{\circle*{3}}
 \multiput(12,0)(24,0){2}{\line(0,1){12}}
 \put(12,12){\line(1,0){24}}
\end{picture}  \\ \hline
$(1)(2)(3,4)$ &
\begin{picture}(46,18)(-3,0) 
 \multiput(0,0)(12,0){4}{\circle*{3}}
 \multiput(24,0)(12,0){2}{\line(0,1){12}}
 \put(24,12){\line(1,0){12}}
\end{picture}
& \hspace*{1.5ex} $(1)(2)(3)(4)$ &
\begin{picture}(60,18)(-3,0) 
 \multiput(0,0)(12,0){4}{\circle*{3}}
\end{picture} 
\end{tabular}
\end{table}
\begin{table}[b]
\caption{Selected elements of $\NCL(5)$ and graphical representations.}
\label{tab:ncnp5}
\begin{tabular}{r|l||r|l}
$\pi$\hspace*{2ex} & \hspace*{2ex}$\Gc_\pi$ & $\pi$\hspace*{2ex} & \hspace*{2ex}$\Gc_\pi$ \\ \hline\hline
$(1,2,3,5)(3,4)$ &
\begin{picture}(72,18)(-3,0) 
 \multiput(0,0)(12,0){5}{\circle*{3}}
 \multiput(0,0)(12,0){3}{\line(0,1){12}}
 \put(48,0){\line(0,1){12}}
 \put(0,12){\line(1,0){48}}
 \put(24,0){\line(3,2){12}}
 \put(36,0){\line(0,1){8}}
\end{picture} 
& \hspace*{1.5ex} $(1,2,5)(2,3,4)$ &
\begin{picture}(72,18)(-3,0) 
 \multiput(0,0)(12,0){5}{\circle*{3}}
 \multiput(0,0)(12,0){2}{\line(0,1){12}}
 \put(48,0){\line(0,1){12}}
 \put(0,12){\line(1,0){48}}
 \multiput(24,0)(12,0){2}{\line(0,1){8}}
 \put(12,0){\line(3,2){12}}
 \put(24,8){\line(1,0){12}}
\end{picture} \\ \hline
$(1,2,5)(2,4)(3)$ &
\begin{picture}(72,18)(-3,0) 
 \multiput(0,0)(12,0){5}{\circle*{3}}
 \multiput(0,0)(12,0){2}{\line(0,1){12}}
 \put(48,0){\line(0,1){12}}
 \put(0,12){\line(1,0){48}}
 \put(36,0){\line(0,1){8}}
 \put(12,0){\line(3,2){12}}
 \put(24,8){\line(1,0){12}}
\end{picture}
& \hspace*{1.5ex} $(1,2)(2,5)(3,4)$ &
\begin{picture}(72,18)(-3,0) 
 \multiput(0,0)(12,0){5}{\circle*{3}}
 \multiput(0,0)(12,0){2}{\line(0,1){12}}
 \put(48,0){\line(0,1){12}}
 \put(0,12){\line(1,0){12}}
 \put(24,12){\line(1,0){24}}
 \put(12,0){\line(1,1){12}}
 \multiput(24,0)(12,0){2}{\line(0,1){8}}
 \put(24,8){\line(1,0){12}}
\end{picture}
\end{tabular}
\end{table}
In order to simplify notation, in these tables we use parentheses $(\cdots)$ rather than brackets
$\{\cdots\}$ to indicate the groupings.
We also draw graphical representations of elements of $\NCL(n)$.
These are modifications of the usual pictures of noncrossing
partitions, in the following way.
The noncrossing partitions $\pi\in\NC(n)\subseteq\NCL(n)$ are drawn in the usual way and
with all angles being right angles.
Suppose $\pi\in\NCL(n)\backslash\NC(n)$.
If $E\in\pi$ and if part~(ii--a) of Definition~\ref{def:ncnp} holds for some $F\in\pi$, $F\ne E$,
then the line connecting $\min(E)$ to the other elements of $E$ is drawn diagonally.

\begin{prop}\label{prop:po}
The relation $\le$ of Definition~\ref{def:le} is a partial order on $\NCL(n)$.
\end{prop}
\begin{proof}
It is clear that $\le$ is transitive and reflexive.
Suppose  $\pi\le\sigma$ and $\sigma\le\pi$, and let us show $\pi=\sigma$.
Let $E\in\pi$.
Then there are $E'\in\sigma$ and $E''\in\pi$ such that $E\subseteq E'\subseteq E''$.
We claim that $E=E''$.
If not, then by Remark~\ref{rmk:ncnp}, we must have $|E|\ge2$ and $|E\cap E''|=1$,
which contradicts $E\subseteq E''$.
Therefore, $E=E''$ and $\pi\subseteq\sigma$.
By symmetry $\sigma\subseteq\pi$.
\end{proof}

\begin{remark}
If $n\ge4$, then $\NCL(n)$ is not a lattice.
We will treat the case $n=4$; the other cases are similar.
Consider the elements
\[
\pi=\{\{1\},\{2,3\},\{4\}\},\qquad
\sigma=\{\{1\},\{2,4\},\{3\}\}
\]
of $\NCL(4)$.
We will show that $\pi$ and $\sigma$ have no least common upper bound.
Consider
\[
\rho=\{\{1\},\{2,3,4\}\},\qquad
\tau=\{\{1,2,4\},\{2,3\}\}
\]
in $\NCL(4)$.
Then $\rho$ and $\tau$ lie above both $\pi$ and $\sigma$.
However, suppose $\alpha\in\NCL(4)$ and $\alpha\le\rho$, $\alpha\le\tau$.
We must have $\{1\}\in\alpha$.
This precludes any other element of $\alpha$ having $1$ as an element.
If $E\in\alpha$ and $2\in E$, then $2=\min(E)$.
This implies that exactly one element of $\alpha$ contains $2$.
This element of $\alpha$ must be either $\{2\}$, $\{2,3\}$ or $\{2,4\}$.
Therefore, either $\pi\not\le\alpha$ or $\sigma\not\le\alpha$.
\end{remark}

The partial ordering of $\NCL(n)$ restricts to the usual partial ordering
on the noncrossing partitions $\NC(n)\subseteq\NCL(n)$.
As is well known, $\NC(n)$ is a lattice.
Moreover, given $\pi\in\NCL(n)$, by taking unions of blocks of $\pi$
that intersect, we find a smallest element $\pihat\in\NC(n)$
such that $\pi\le\pihat$.
\begin{defi}
We call $\pihat$ as described above the {\em noncrossing partition generated by} $\pi$.
\end{defi}

\begin{remark}\label{rem:gen}
Obviously, if $p$ and $q$ belong to
the same block of $\pihat$, then there are $k\ge1$ and $B_1,\ldots,B_k\in\pi$
with $p\in B_1$, $q\in B_k$ and $B_j\cap B_{j+1}\ne\emptyset$ whenever $1\le j\le k-1$.
We may clearly choose $B_1,\ldots,B_k$ to be distinct.
\end{remark}

\begin{defi}\label{def:contr}
Let $n\in\Nats$, let $\pi$ be a set of subsets of $\{1,\ldots,n\}$ and let $X\subseteq\{1,\ldots,n\}$
be a nonempty subset.
Let $\chi:X\to\{1,\ldots,|X|\}$ be the order preserving bijection.
The {\em renumbered restriction} of $\pi$ to $X$ is the set
\[
\pi\contr_X:=\{\chi(E\cap X)\mid E\in\pi,\,E\cap X\ne\emptyset\}.
\]
\end{defi}

Clearly, if $\pi\in\NC(n)$, then $\pi\contr_X\in\NC(|X|)$.
The analogous statement for $\NCL$ is not true.

We now discuss some notions that justify calling the elements of $\NCL(n)$ noncrossing
linked partitions.
\begin{defi}
Let $n\ge1$, $\pi\in\NCL(n)$, $B\in\pi$ and $i=\min(B)$.
If $i$ is doubly covered by $\pi$, then let $\Bchk=B\backslash\{i\}$.
If $i$ is singly covered by $\pi$, then let $\Bchk=B$.
Let
\[
\pichk=\{\Bchk\mid B\in\pi\}.
\]
We call $\pichk$ the {\em unlinking} of $\pi$.
\end{defi}

It is clear from the definitions that $\pichk\in\NC(n)$ and $\pichk\le\pi$.
Moverover, $\pichk$ is a maximal element of $\{\sigma\in\NC(n)\mid\sigma\le\pi\}$.
However, it is not necessarily the unique such element:
with $\pi=(1,2)(2,3)$ we
have $\pichk=(1,2)(3)$, but also $\sigma\le\pi$, with $\sigma=(1)(2,3)$.
It is also clear that the map $\pi\mapsto\pichk$ is order preserving.

\begin{prop}\label{prop:NCL1}
Let $n\in\Nats$, $n\ge2$ and let
\[
\NCL^{(1)}(n)=\{\pi\in\NCL(n)\mid\pichk=1_n\}.
\]
If $\pi\in\NCL^{(1)}(n)$, then $1$ and $2$ belong to the same block of $\pichk$.
Let $u(\pi)$ be the renumbered restriction of $\pichk$ to $\{2,\ldots,n\}$.
Then $u$ is an order preserving bijection from $\NCL^{(1)}(n)$ onto $\NC(n-1)$.
\end{prop}
\begin{proof}
First, we will prove that $1$ and $2$ are in the same block of $\pichk$.
By the observation made in Remark~\ref{rem:gen}, there are $k\ge1$ and distinct $B_1,\ldots,B_k\in\pi$
such that $1\in B_1$, $2\in B_k$ and $B_j\cap B_{j+1}\ne\emptyset$.
If $k=1$, then $1,2\in B_1$ and, since $1$ is singly covered by $\pi$, $\Bchk_1=B_1\in\pichk$.
If $k\ge2$, then let $i_j$ be such that $B_j\cap B_{j+1}=\{i_j\}$.
(See Remark~\ref{rmk:ncnp}.)
Then $i_1>1$, so $i_1=\min(B_2)$.
Thus, $i_2\ne\min(B_2)$ and $i_2=\min(B_3)$.
Arguing in the way and by induction, we get $1<i_1<i_2<\cdots<i_{k-1}$, and $i_{k-1}=\min(B_k)$.
Since $2\in B_k$, we have $i_{k-1}=2$ and, therefore, $k=2$ and $1,2\in B_1$.
Now, as before, $1,2\in\Bchk_1\in\pichk$.
This proves that $1$ and $2$ belong to the same block of $\pichk$.

Since $\pichk\in\NC(n)$, we immediately get $u(\pi)\in\NC(n-1)$.
Moreover, since the map $\pi\mapsto\pichk$ is order preserving, the same is true of $u$.
In order to show that $u:\NCL^{(1)}(n)\to\NC(n-1)$ is one--to--one and onto,
we will construct its inverse map.
Given $\tau\in\NC(n-1)$, let $\taut$ be the result of adding $1$
to every element of every element of $\tau$,
so that $\tau=\taut\contr_{\{2,\ldots,n\}}$.
Given $B\in\taut$, let $B'=B\cup\{\min(B)-1\}$, and let
\[
v(\tau)=\{B'\mid B\in\taut\}.
\]
We claim $v(\tau)\in\NCL^{(1)}(n)$ and that $v$ is the inverse map of $u$.
Let us begin by showing that distinct elements of $v(\tau)$ are noncrossing.
Suppose, to obtain a contradiction, there are $B_1,B_2\in\taut$ with $B_1\ne B_2$ and there are
$i_1,i_2\in B_1'$, $j_1,j_2\in B_2'$ with $i_1<j_1<i_2<j_2$.
Since $B_1$ and $B_2$ are noncrossing, we must have either $i_1=\min(B_1')$ or $j_1=\min(B_2')$
or both.
Suppose both occur.
Then $i_1+1,i_2\in B_1$ and $j_1+1,j_2\in B_2$.
Since $B_1$ and $B_2$ are disjoint, we have $i_1+1<j_1+1<i_2<j_2$, contradicting that $B_1$ and $B_2$
are noncrossing.
Similar arguments apply to obtain a contradiction assuming only one of $i_1=\min(B_1')$ and
$j_1=\min(B_2')$ holds.
Thus, distinct elements of $v(\tau)$ are noncrossing.
Let us now show that distinct elements of $v(\tau)$ are nearly disjoint.
Note that $v(\tau)$ has no elements that are singletons.
Let $B_1,B_2\in\taut$ with $B_1\ne B_2$, suppose $i\in B_1'\cap B_2'$.
Because $B_1$ and $B_2$ are disjoint, exactly one of
$i\in(B_1'\backslash B_1)\cap B_2$ and $i\in(B_2'\backslash B_2)\cap B_1$ holds.
In the first case, $i=\min(B_1')$ and $i\ne\min(B_2')$, while in the second case
$i\ne\min(B_1')$ and $i=\min(B_2')$.
This shows $v(\tau)\in\NCL(n)$.

Let us now show $v(\tau)\hat{\;}=1_n$.
If $v(\tau)\hat{\;}\ne1_n$, then there is $k\in\{2,\ldots,n\}$ such that every block
of $v(\tau)$ is a subset of either $\{1,\ldots,k-1\}$ or $\{k,\ldots,n\}$.
But there is $B\in\taut$ with $k\in B$, and $B\subseteq B'\in v(\tau)$.
By the hypothesis on $k$, we must have $k=\min(B)$.
Therefore, $k-1\in B'$, contrary to the hypothesis on $k$.
Therefore, $v(\tau)\in\NCL^{(1)}(n)$.

That we have $u\circ v(\tau)=\tau$ follows easily from the observation
that if $B\in\taut$, then
\[
(B')\check{\;}=
\begin{cases}
B\cup\{1\},&2\in B \\
B,&2\notin B.
\end{cases}
\]

Let us show $(v\circ u)(\pi)=\pi$, for $\pi\in\NCL^{(1)}(n)$.
First, we observe that either $\pi=1_n$ or every block of $\pi$ meets some other block of $\pi$.
Indeed, if $B\in\pi$, $B\ne\{1,\ldots,n\}$ and $B$ is disjoint from every other block of $\pi$,
then $B\in\pihat$, contradicting $\pihat=1_n$.
Therefore, for all $B\in\pi$, $\min(B)$ is doubly covered by $\pi$, unless $\min(B)=1$.
Hence,
\[
\Bchk=
\begin{cases}
B,&1\in B \\
B\backslash\{\min(B)\},&1\notin B.
\end{cases}
\]
Letting $\tau=u(\pi)$, this implies
$\taut=\{B\backslash\{\min(B)\}\mid B\in\pi\}$.
From this, we get $v(\tau)=\pi$.
\end{proof}

\begin{remark}
The map $u$ is an order isomorphism from $\NCL^{(1)}(n)$ onto $\NC(n-1)$ when $2\le n\le3$,
but not for $n\ge4$.
To see this when $n=4$, take $\pi=(1,2)(2,3)(3,4)$ and $\sigma=(1,2,4)(2,3)$.
Then $u(\pi)=(1)(2)(3)$ and $u(\sigma)=(1,3)(2)$.
Thus, $\pi\not\le\sigma$, but $u(\pi)\le u(\sigma)$.
\end{remark}

\begin{cor}\label{cor:NCLNC}
If $\sigma\in\NC(n)$, then the cardinality of the set
\[
\NCL^{(\sigma)}(n):=\{\pi\in\NCL(n)\mid\pihat=\sigma\}
\]
is
\[
|\NCL^{(\sigma)}(n)|=\prod_{B\in\sigma}c_{|B|-1},
\]
where for $k\in\Nats_0$,
\begin{equation}\label{eq:ck}
c_k=
\begin{cases}
\frac1k\binom{2k}{k-1},&k\ge1 \\
1,&k=0  
\end{cases}
\end{equation}
is the Catalan number (with domain extended to include $k=0$).
Consequently,
\begin{equation}\label{eq:NCLNC}
|\NCL(n)|=\sum_{\sigma\in\NC(n)}\big(\prod_{B\in\sigma}c_{|B|-1}\big).
\end{equation}
\end{cor}
\begin{proof}
Clearly, $\NCL^{(\sigma)}(n)$ can be written as the Cartesian product over $B\in\sigma$
of $\NCL^{(1)}(|B|)$.
By Proposition~\ref{prop:NCL1}, $|\NCL^{(1)}(n)|=|\NC(n-1)|$ if $n\ge2$.
It is well known that $\NC(k)$ has cardinality equal to the Catalan number $c_k$.
Finally, one sees $|\NCL^{(1)}(1)|=1$ immediately.
\end{proof}

The set $\NCL(n)$ has another natural order, defined below.
\begin{defi}
Let $\pi,\sigma\in\NCL(n)$.
We write $\pi\lenc\sigma$ if $\pihat\le\sigmahat$ and either $\pihat\ne\sigmahat$
or $\pihat=\sigmahat$ but $\pichk\le\sigmachk$.
\end{defi}
Clearly $\lenc$ is a partial order making $\NCL(n)$ into a lattice.
Both $\lenc$ and the order $\le$ from Definition~\ref{def:le} extend the usual
order on $\NC(n)$.

\begin{defi}
If $\pi$ is a set of subsets of $\{1,\ldots,n\}$
and if $J\subseteq\{1,\ldots,n\}$, we say that $J$ {\em splits} $\pi$
if for every $F\in\pi$, either $F\subseteq J$ or $F\cap J=\emptyset$.
\end{defi}

\begin{defi}
If $m,n\in\Nats$ and if $\pi\in\NCL(m)$ and $\sigma\in\NCL(n)$, then let
\[
\pi\oplus\sigma\in\NCL(m+n)
\]
be $\pi\oplus\sigma=\pi\cup\sigmat$, where $\sigmat$
is the right translation of $\sigma$ by $m$, i.e. $\sigmat$
is obtained by adding $m$ to every element of every element of $\sigma$.
\end{defi}

Take two disjoint copies of the nonnegative integers, written
$\Nats_0=\{0,1,2,3,\ldots\}$ and $\Nats_0^*=\{0^*,1^*,2^*,\ldots\}$.

\begin{defi}\label{def:Spi}
Let $n\in\Nats$ and $\pi\in\NCL(n)$.
Define the function
\[
S_\pi:\{1,\ldots,n\}\to\Nats_0\cup\Nats_0^*
\]
as follows:
\begin{itemize}
\item $S_\pi(j)=0^*$ if $j$ is not the minimal element of a block of $\pi$;
\item $S_\pi(j)=|F|-1$ if $j=\min(F)$ for some $F\in\pi$ and $j$ is singly covered by $\pi$;
\item $S_\pi(j)=(|F|-1)^*$ if $j=\min(F)$ for some $F\in\pi$ and $j$ is doubly covered by $\pi$.
\end{itemize}
\end{defi}

\begin{prop}\label{prop:spi}
Let $n\in\Nats$ and let $\pi,\sigma\in\NCL(n)$.
If $S_\pi=S_\sigma$, then $\pi=\sigma$.
\end{prop}
\begin{proof}
We use induction on $n$.
The case $n=1$ is clear, since $\NCL(1)$ has only one element.
Take $n\ge2$.
Let $s=S_\pi$.
Let $m$ be largest such that $s(m)\ne0^*$.
If $m=1$, then we must have $\pi=\sigma=\{\{1,\ldots,n\}\}$, so suppose $m\ge2$.
If $s(m)=k\in\Nats_0$ or $s(m)=k^*\in\Nats_0^*$, then the interval
$J=\{m,\ldots,m+k\}$ belongs to both $\pi$ and $\sigma$.
Indeed, we know that $m$ is the minimal element of blocks of cardinality $k+1$ 
of both $\pi$ and $\sigma$.
But this block cannot have any gaps, by the noncrossing property and by the choice of $m$.

Suppose $s(m)=k\in\Nats_0$
and let $\pit$ be the renumbered restriction of $\pi\backslash\{J\}$ to $J^c$.
and let
$\sigmat$ be the renumbered restriction of $\sigma\backslash\{J\}$ to $J^c$.
Then, $\pit,\sigmat\in\NCL(n-k-1)$.
Moreover, $S_\pit$ and $S_\sigmat$ are both equal to the composition of $s$ and the map
\begin{equation*}
i\mapsto
\begin{cases}
i,&1\le i\le m-1 \\
i-k-1,&m+k+1\le i\le n.
\end{cases}
\end{equation*}
By the induction hypothesis, $\sigmat=\pit$, and then $\sigma=\pi$ follows.

Suppose $s(m)=k^*\in\Nats_0^*$.
Then $k\ge1$.
Let $\pit$ and $\sigmat$ be the renumbered restrictions of $\pi\backslash\{J\}$
and, respectively, $\sigma\backslash\{J\}$ to $J^c\cup\{m\}$.
Then $\pit,\sigmat\in\NCL(n-k)$ and we have $S_\pit=S_\sigmat$ equal to the map
\[
i\mapsto
\begin{cases}
s(i),&1\le i\le m-1 \\
0^*,&i=m \\
s(i+k),&m+1\le i\le n-k.
\end{cases}
\]
By the induction hypothesis, $\sigmat=\pit$, and then $\sigma=\pi$ follows.
\end{proof}

\section{The unsymmetrized R--transform}
\label{sec:Rtrans}

In this section, we relate the additive canonical random variables of~\S\ref{sec:crv}
to the unsymmetrized R--transform defined at~\eqref{eq:Rt}
(and more fomally below), and use them to show additivity~\eqref{eq:Rtplus}
under hypothesis of freeness.

We now define the notation
\[
\alpha_\pi[b_1,\ldots,b_n],
\]
for $\alpha\in\Mul[[B]]$, $n\in\Nats_0$, $\pi\in\NC(n+1)$ and $b_1,\ldots,b_n\in B$.
(As usual, $\Nats_0$ will denote the set of nonnegative integers.)
This is inspired by work of Speicher~\cite{Sp1}, \cite{Sp} and is in essence equivalent to a part of
his definition.
However, we take a perspective dual to Speicher's in that we consider multilinear functions
on $B\tdt B$ defined by particular elements in $B$--valued noncommutative probability spaces,
rather than considering functions on $B,B$--bimodules.

\begin{defi}\label{def:alphapi}
If $\pi=1_{n+1}$, then let
\begin{equation}\label{eq:alpha_1}
\alpha_\pi[b_1,\ldots,b_n]=
\alpha_n(b_1,\ldots,b_n),
\end{equation}
where in the case $n=0$, the right--hand--side of~\eqref{eq:alpha_1} means $\alpha_0$.
Otherwise, let $J=\{m,m+1,\ldots,m+\ell-1\}$ be the right--most block
of $\pi$ that is an interval
and let $\pi'\in\NC(n-\ell+1)$ be the renumbered restriction of $\pi$ to
$J^c=\{1,\ldots,n\}\backslash J$.
If there is nothing to the right of the interval $J$, i.e.\ if $m+\ell-1=n+1$, then set
\begin{equation}\label{eq:atrt}
\alpha_\pi[b_1,\ldots,b_n]=
\alpha_{\pi'}[b_1,\ldots,b_{n-\ell}]b_{m-1}\alpha_{\ell-1}(b_m,\ldots,b_n).
\end{equation}
If $m+\ell-1<n+1$, then set
\begin{equation}\label{eq:alphapi}
\alpha_\pi[b_1,\ldots,b_n]=
\alpha_{\pi'}[b_1,\ldots,b_{m-2},b_{m-1}\alpha_{\ell-1}(b_m,\ldots,b_{m+\ell-2})b_{m+\ell-1},
b_{m+\ell},\ldots,b_n].
\end{equation}
By recursion,~\eqref{eq:alpha_1}, \eqref{eq:atrt} and~\eqref{eq:alphapi}
define $\alpha_\pi[b_1,\ldots,b_n]$ for arbitrary $\pi\in\NC(n+1)$.
\end{defi}

Some examples of $\alpha_\pi[b_1,\ldots,b_n]$ are given in Table~\ref{tab:alphapi}.
\begin{table}[b]
\caption{Some examples of $\alpha_\pi[b_1,b_2,b_3]$.}
\label{tab:alphapi}
\begin{tabular}{r|c|l}
$\pi\in\NC(4)$ & $\Gc_\pi$ & $\alpha_\pi[b_1,b_2,b_3]$ \\ \hline\hline
$(1,2,3,4)$ &
\begin{picture}(46,18)(-3,0) 
 \multiput(0,0)(12,0){4}{\circle*{3}}
 \multiput(0,0)(12,0){4}{\line(0,1){12}}
 \put(0,12){\line(1,0){36}}
\end{picture}
& $\alpha_3(b_1,b_2,b_3)$ \\ \hline
$(1,2,3)(4)$ &
\begin{picture}(46,18)(-3,0) 
 \multiput(0,0)(12,0){4}{\circle*{3}}
 \multiput(0,0)(12,0){3}{\line(0,1){12}}
 \put(0,12){\line(1,0){24}}
\end{picture}
& $\alpha_2(b_1,b_2)b_3\alpha_0$ \\ \hline
$(1,3,4)(2)$ &
\begin{picture}(46,18)(-3,0) 
 \multiput(0,0)(12,0){4}{\circle*{3}}
 \put(0,0){\line(0,1){12}}
 \multiput(24,0)(12,0){2}{\line(0,1){12}}
 \put(0,12){\line(1,0){36}}
\end{picture}
& $\alpha_2(b_1\alpha_0b_2,b_3)$ \\ \hline
$(1,2)(3,4)$ &
\begin{picture}(46,18)(-3,0) 
 \multiput(0,0)(12,0){4}{\circle*{3}}
 \multiput(0,0)(12,0){4}{\line(0,1){12}}
 \multiput(0,12)(24,0){2}{\line(1,0){12}}
\end{picture}
& $\alpha_1(b_1)b_2\alpha_1(b_3)$ \\ \hline
$(1,4)(2,3)$ &
\begin{picture}(46,18)(-3,0) 
 \multiput(0,0)(12,0){4}{\circle*{3}}
 \multiput(0,0)(36,0){2}{\line(0,1){12}}
 \multiput(12,0)(12,0){2}{\line(0,1){8}}
 \put(0,12){\line(1,0){36}}
 \put(12,8){\line(1,0){12}}
\end{picture}
& $\alpha_1(b_1\alpha_1(b_2)b_3)$ \\ \hline
$(1)(2)(3)(4)$ &
\begin{picture}(46,18)(-3,0) 
 \multiput(0,0)(12,0){4}{\circle*{3}}
\end{picture}
& $\alpha_0b_1\alpha_0b_2\alpha_0b_3\alpha_0$
\end{tabular}
\end{table}

\begin{defi}
Let $\Nats_{-1}=\{-1,0,1,2,3,\ldots\}$.
For $n\in\Nats$ and $\pi\in\NC(n)$, let
$K_\pi:\{1,\ldots,n\}\to\Nats_{-1}$ be
\[
K_\pi(j)=
\begin{cases}
|F|-1,&j=\min(F),\,F\in\pi \\
-1,&j\text{ not the minimal element of a block of }\pi.
\end{cases}
\]
\end{defi}

Clearly, $\pi$ is determined by $K_\pi$.

\begin{lemma}\label{lem:ncp}
Let $\alpha\in\Mul[[B]]$ and for $j\ge-1$ let
\begin{equation}\label{eq:xk}
x_j=
\begin{cases}
V_{1,j}(\alpha_j),&j\ge0 \\
L_1,&j=-1
\end{cases}
\end{equation}
in $\Lc(\Qc)$.
Let $n\in\Nats_0$, and let
$k:\{1,\ldots,n+1\}\to\Nats_{-1}$.
Let $b_1,\ldots,b_n\in B$.
If there is $\pi\in\NC(n+1)$ such that $K_\pi=k$, then
\begin{equation}\label{eq:xktermpi}
x_{k(1)}b_1x_{k(2)}b_2\cdots x_{k(n)}b_nx_{k(n+1)}=\alpha_\pi[b_1,\ldots,b_n].
\end{equation}
If there is no $\pi\in\NCL(n+1)$ such that $K_\pi=k$, then
\begin{equation}\label{eq:Ec0}
\Ec(x_{k(1)}b_1x_{k(2)}b_2\cdots x_{k(n)}b_nx_{k(n+1)})=0.
\end{equation}
\end{lemma}
\begin{proof}
One easily shows (see~\cite[Lemma 3.3]{D})
\begin{equation}\label{eq:V1n}
\begin{aligned}
V_{1,0}(\alpha_0)&=\alpha_0 \\
V_{1,n}(\alpha_n)b_1L_1b_2L_1\cdots b_nL_1&=\alpha_n(b_1,\ldots,b_n),\quad(n\ge1).
\end{aligned}
\end{equation}
Let $N=|\{j\mid k(j)\ge0\}|$.
We proceed by induction on $N$.
If $N=0$, then $k(1)=-1$.
But $PL_1=0$, where $P$ is the projection used in defining $\Ec$, so~\eqref{eq:Ec0} holds.
Moreover, there is no $\pi$ such that $K_\pi=k$ and the lemma holds when $N=0$.

Suppose $N=1$.
The only way in which we can have $k=K_\pi$ for some $\pi\in\NC(n+1)$ is if $\pi=1_{n+1}$
and $k(1)=n$, $k(j)=-1$ for all $j\in\{2,\ldots,n+1\}$.
But then, by~\eqref{eq:V1n}, 
\[
x_{k(1)}b_1x_{k(2)}b_2\cdots x_{k(n)}b_nx_{k(n+1)}=\alpha_n(b_1,\ldots,b_n)=
\alpha_{1_{n+1}}[b_1,\ldots,b_n].
\]
On the other hand, if $N=1$ and $k\ne K_{1_{n+1}}$, then we easily see that~\eqref{eq:Ec0} holds.

Suppose $N\ge2$.
Let $m$ be largest such that $k(m)\ge0$.
Then $m\ge2$.
Suppose $k(m)>n-m+1$.
Then
\[
x_{k(m)}b_m\cdots x_{k(n)}b_nx_{k(n+1)}\Omega
=V_{1,k(m)}(\alpha_{k(m)})(b_m\delta_1\otdt b_n\delta_1\otimes 1)=0,
\]
so~\eqref{eq:Ec0} holds.
But also there is no $\pi\in\NC(n+1)$ such that $K_\pi=k$.

Suppose $k(m)=n-m+1$.
Then there is $\pi\in\NC(n+1)$ with $K_\pi=k$ if and only if there is $\pi'\in\NC(m-1)$
such that $K_{\pi'}$ is the restriction of $k$ to $\{1,\ldots,m-1\}$, and then
$\pi=\pi'\cup\{\{m,\ldots,n+1\}\}$.
Using~\eqref{eq:V1n}, we get
\[
x_{k(m)}b_m\cdots x_{k(n)}b_nx_{k(n+1)}
=\alpha_{k(m)}(b_m,\ldots,b_n).
\]
By using the induction hypothesis and Definition~\ref{def:alphapi},
if there is $\pi\in\NC(n+1)$ such that $K_\pi=k$, then with $\pi'$ as above,
\begin{align*}
x_{k(1)}b_1\cdots x_{k(n)}b_nx_{k(n+1)}
&=x_{k(1)}b_1\cdots x_{k(m-1)}b_{m-1}\alpha_{k(m)}(b_m,\ldots,b_n) \\
&=\alpha_{\pi'}[b_1,\ldots,b_{m-2}]b_{m-1}\alpha_{k(m)}(b_m,\ldots,b_n)
=\alpha_\pi[b_1,\ldots,b_n].
\end{align*}
Using again the induction hypothesis, if there is no $\pi$ such that $K_\pi=k$, then
\begin{equation}\label{eq:Er}
\begin{split}
\Ec(x_{k(1)}b_1\cdots&x_{k(n)}b_nx_{k(n+1)})= \\
&=\Ec(x_{k(1)}b_1\cdots x_{k(m-2)}b_{m-2}x_{k(m-1)})b_{m-1}\alpha_{n-m+1}(b_m,\ldots,b_n)=0.
\end{split}
\end{equation}

Suppose $k(m)<n-m+1$.
Let $\ell=k(m)$, and let $k':\{1,\ldots,n-\ell\}\to\Nats_{-1}$ be
\[
k'(j)=
\begin{cases}
k(j),&1\le j\le m-1 \\
k(j+\ell+1),&m\le j\le n-\ell.
\end{cases}
\]
Then there is $\pi\in\NC(n+1)$ such that $K_\pi=k$ if and only if there is $\pi'$ in $\NC(n-\ell)$
such that $K_{\pi'}=k'$,
and then $\pi$ is obtained from $\pi'$ by applying the map
\[
j\mapsto
\begin{cases}
j,&1\le j\le m-1 \\
j+\ell+1,&m\le j\le n-\ell
\end{cases}
\]
to every element of every element of $\pi'$ and then adding the block $J=\{m,\ldots,m+\ell\}$.
Thus, $\pi'$ is the renumbered restriction of $\pi$ to $J^c$.
From~\eqref{eq:V1n}, we have
\[
x_{k(m)}b_m\cdots x_{k(m+\ell-1)}b_{m+\ell-1}x_{k(m+\ell)}
=\alpha_\ell(b_m,\ldots,b_{m+\ell-1}).
\]
Let $\bt=b_{m-1}\alpha_\ell(b_m,\ldots,b_{m+\ell-1})b_{m+\ell}$.
Suppose there is $\pi\in\NC(n+1)$ such that $K_\pi=k$.
Then letting $\pi'$ be as above, by the induction hypothesis and Definition~\ref{def:alphapi},
\begin{align*}
x_{k(1)}b_1\cdots&x_{k(n)}b_nx_{k(n+1)}= \\
&=x_{k(1)}b_1\cdots x_{k(m-2)}b_{m-2}x_{k(m-1)}\bt x_{k(m+\ell+1)}b_{m+\ell+1}\cdots x_{k(n)}b_nx_{k(n+1)} \\
&=\alpha_{\pi'}[b_1,\ldots,b_{m-1},\bt,b_{m+\ell+1},\ldots,b_n]
=\alpha_\pi[b_1,\ldots,b_n].
\end{align*}
Again using the induction hypothesis, if there is no $\pi$ such that $K_\pi=k$, then
\begin{multline*}
\Ec(x_{k(1)}b_1\cdots x_{k(n)}b_nx_{k(n+1)})= \\
=\Ec(x_{k(1)}b_1\cdots x_{k(m-2)}b_{m-2}x_{k(m-1)}\bt x_{k(m+\ell+1)}b_{m+\ell+1}\cdots x_{k(n)}b_nx_{k(n+1)})=0.
\end{multline*}
This completes the induction step.
\end{proof}

\begin{lemma}\label{lem:insertNC}
Let $n\in\Nats$, 
let $\pi\in\NC(n+1)$ and let 
\[
J=\{m,m+1,\ldots,m+\ell-1\}
\]
be a proper subinterval of $\{1,\ldots,n+1\}$.
Suppose that $J$ splits $\pi$.
Let $\pi''\in\NC(\ell)$ be the renumbered restriction of $\pi$ to $J$
and let $\pi'\in\NC(n-\ell+1)$ be the renumbered restriction of $\pi$ to $J^c$.
Let $\alpha\in\Mul[[B]]$ and $b_1,\ldots,b_n\in B$.
Then
\begin{multline}\label{eq:subinterval}
\alpha_\pi[b_1,\ldots,b_n]= \\
=\begin{cases}
\alpha_{\pi''}[b_1,\ldots,b_{\ell-1}]b_\ell\alpha_{\pi'}[b_{\ell+1},\ldots,b_n],&m=1 \\
\begin{aligned}[b]
\alpha_{\pi'}[b_1,&\ldots,b_{m-2}, \\
&b_{m-1}\alpha_{\pi''}[b_m,\ldots,b_{m+\ell-2}]b_{m+\ell-1},
b_{m+\ell},\ldots,b_n],\end{aligned}&2\le m\le n-\ell+1 \\
\alpha_{\pi'}[b_1,\ldots,b_{m-2}]b_{m-1}\alpha_{\pi''}[b_m,\ldots,b_n],&m=n-\ell+2.
\end{cases}
\end{multline}
\end{lemma}
\begin{proof}
Suppose $2\le m\le n-\ell+1$.
By Lemma~\ref{lem:ncp}, letting $k=K_\pi$, we have
\[
\alpha_\pi[b_1,\ldots,b_n]=x_{k(1)}b_1\cdots x_{k(n)}b_nx_{k(n+1)}.
\]
But, since $K_{\pi''}(j)=k(j+m-1)$, again by Lemma~\ref{lem:ncp} we have
\[
x_{k(m)}b_m\cdots x_{k(m+\ell-2)}b_{m+\ell-2}x_{k(m+\ell-1)}=\alpha_{\pi''}[b_m,\ldots,b_{m+\ell-2}].
\]
So letting $\bt=b_{m-1}\alpha_{\pi''}[b_m,\ldots,b_{m+\ell-2}]b_{m+\ell-1}$,
we have
\[
\alpha_\pi[b_1,\ldots,b_n]=x_{k(1)}b_1\cdots x_{k(m-2)}b_{m-2}x_{k(m-1)}
\bt x_{k(m+\ell)}b_{m+\ell}\cdots x_{k(n)}b_nx_{k(n+1)}.
\]
But
\[
K_{\pi'}(j)=
\begin{cases}
k(j),&1\le j\le m-1 \\
k(j+\ell),&m\le j\le n-\ell+1.
\end{cases}
\]
Therefore, applying Lemma~\ref{lem:ncp} again, we get~\eqref{eq:subinterval} in the case $2\le m\le n-\ell+1$.
The other cases are proved similarly.
\end{proof}

We'll want to use the particular case of an interval partition.
\begin{lemma}\label{lem:alphaint}
Let $\pi\in\IP(n+1)$ be an interval partition.
Let $\pi=\{F_1,\ldots,F_k\}$, where $|F_j|=p_j$ and $F_j=\{q_j+1,\ldots,q_j+p_j\}$
with $q_1=0$ and $q_j=p_1+\cdots+p_{j-1}$.
Then
\[
\begin{split}
\alpha_\pi[b_1,\ldots,b_n]=&\alpha_{p_1-1}(b_{q_1+1},\ldots,b_{q_1+p_1-1})b_{q_2}
\alpha_{p_2-1}(b_{q_2+1},\ldots,b_{q_2+p_2-1})\cdots \\
&\quad b_{q_k}
\alpha_{p_k-1}(b_{q_k+1},\ldots,b_{q_k+p_k-1}).
\end{split}
\]
\end{lemma}

\begin{prop}\label{prop:Xaddtve}
Let $\alpha\in\Mul[[B]]$ and let 
\[
X=L_1+V_1(\alpha).
\]
Then for every $n\ge0$ and $b_1,\ldots,b_n\in B$,
\begin{equation}\label{eq:ncp}
\Phit_{X,n}(b_1,\ldots,b_n)=\sum_{\pi\in\NC(n+1)}\alpha_\pi[b_1,\ldots,b_n].
\end{equation}
\end{prop}
\begin{proof}
We have
\begin{align}
\Phit_{X,n}(b_1,\ldots,b_n)&=\Ec(Xb_1Xb_2\cdots Xb_nX) \notag \\[1ex]
&=\sum_{k(1),\ldots,k(n+1)\ge-1}\Ec(x_{k(1)}b_1x_{k(2)}b_2\cdots x_{k(n)}b_nx_{k(n+1)}), \label{eq:xksum}
\end{align}
where $x_k$ is as in~\eqref{eq:xk}.
By Lemma~\ref{lem:ncp},
only finitely many terms of the sum~\eqref{eq:xksum} are nonzero, and the equality~\eqref{eq:ncp} holds.
\end{proof}

\begin{lemma}\label{lem:ncpXY}
Let $\alpha,\beta\in\Mul[[B]]$ and let 
\begin{equation}\label{eq:xij}
x_j^i=
\begin{cases}
V_{1,j}(\alpha_j),&i=1,\,j\ge0 \\
L_1,&i=1,\,j=-1 \\
V_{2,j}(\beta_j),&i=2,\,j\ge0 \\
L_2,&i=2,\,j=-1\end{cases}
\end{equation}
in $\Lc(\Qc)$.
Let $k:\{1,\ldots,n+1\}\to\Nats_{-1}$.
If there is $\pi\in\NC(n+1)$ such that $K_\pi=k$, then
\begin{equation*}
\sum_{i(1),\ldots,i(n+1)\in\{1,2\}}
x_{k(1)}^{i(1)}b_1x_{k(2)}^{i(2)}b_2\cdots x_{k(n)}^{i(n)}b_nx_{k(n+1)}^{i(n+1)}=
(\alpha+\beta)_\pi[b_1,\ldots,b_n].
\end{equation*}
If there is no $\pi$ such that $K_\pi=k$, then
\begin{equation*}
\Ec(x_{k(1)}^{i(1)}b_1x_{k(2)}^{i(2)}b_2\cdots x_{k(n)}^{i(n)}b_nx_{k(n+1)}^{i(n+1)})=0.
\end{equation*}
for any $i(1),\ldots,i(n+1)\in\{1,2\}$.
\end{lemma}
\begin{proof}
This can be proved by induction just like Lemma~\ref{lem:ncp} was proved, but using also
the fact that for any $\gamma_n\in\Lc_n(B)$,
we have
\[
V_{i,n}(\gamma_n)b_1L_{j_1}b_2L_{j_2}\cdots b_mL_{j_m}=0
\]
if $1\le m\le n$ and $i,j_1,\ldots,j_m\in\{1,2\}$ with $j_p\ne i$ for some $p$.
\end{proof}

\begin{lemma}\label{lem:X+Y}
Let $\alpha,\beta\in\Mul[[B]]$ and let 
\begin{align}
X&=L_1+V_1(\alpha) \label{eq:Xagain} \\
Y&=L_2+V_2(\beta) \label{eq:Y} \\
Z&=L_1+V_1(\alpha+\beta). \label{eq:Z}
\end{align}
Then $X+Y$ and $Z$ have the same distribution series.
\end{lemma}
\begin{proof}
By Proposition~\ref{prop:Xaddtve}, the distribution series of $Z$ is given by
\begin{equation}\label{eq:PhiZ}
\Phit_{Z,n}(b_1,\ldots,b_n)=\sum_{\pi\in\NC(n+1)}(\alpha+\beta)_\pi[b_1,\ldots,b_n].
\end{equation}
On the other hand, we have
\begin{multline}\label{eq:X+Ysum}
\Phit_{X+Y,n}(b_1,\ldots,b_n)=\Ec((X+Y)b_1(X+Y)b_2\cdots (X+Y)b_n(X+Y)) \\[1ex]
=\sum_{\substack{k(1),\ldots,k(n+1)\ge-1 \\ i(1),\ldots,i(n+1)\in\{1,2\}}}
\Ec(x_{i(1),k(1)}b_1x_{i(2),k(2)}b_2\cdots x_{i(n),k(n)}b_nx_{i(n+1),k(n+1)}),
\end{multline}
where $x_{i,j}$ is defined in~\eqref{eq:xij}.

By Lemma~\ref{lem:ncpXY}, only finitely many terms of the sum~\eqref{eq:X+Ysum} are nonzero and
we have the equality
\[
\Phit_{X+Y,n}(b_1,\ldots,b_n)=\sum_{\pi\in\NC(n+1)}(\alpha+\beta)_\pi[b_1,\ldots,b_n].
\]
Combined with~\eqref{eq:PhiZ}, this finishes the proof.
\end{proof}

\begin{lemma}\label{lem:ag}
Let $\beta\in\Mul[[B]]$ and let $\alpha=I+I\beta I$.
Then $\alpha$ is invertible with respect to formal composition, and
\begin{equation}\label{eq:ag}
\alpha\cinv=(1+I\gamma)^{-1}I=I(1+\gamma I)^{-1},
\end{equation}
where
\begin{equation}\label{eq:gam}
\gamma=((1+\beta I)^{-1}\beta)\circ\alpha\cinv=(\beta(1+I\beta)^{-1})\circ\alpha\cinv.
\end{equation}
Moreover, $\gamma$ is the unique element of $\Mul[[B]]$ making~\eqref{eq:ag} hold.
\end{lemma}
\begin{proof}
The second equality in
each of~\eqref{eq:ag} and~\eqref{eq:gam} follows from Proposition~\ref{prop:1+a}.
Since $\alpha_0=0$ and $\alpha_1=\id_B$ is invertible, $\alpha$ is invertible with respect to
composition.
Using some of the properties listed in Proposition~\ref{prop:Mul},
we calculate
\begin{align*}
((1+I\gamma)^{-1}I)\circ\alpha&=((1+I\gamma)^{-1}\circ\alpha)\alpha
=((1+I\gamma)\circ\alpha)^{-1}\alpha=(1+\alpha(\gamma\circ\alpha))^{-1}\alpha \\
&=(1+\alpha(1+\beta I)^{-1}\beta)^{-1}\alpha=(1+I\beta)^{-1}\alpha=I,
\end{align*}
where in the last line we used $\alpha=(1+I\beta)I=I(1+\beta I)$.
Thus~\eqref{eq:ag} holds.
Uniqueness of $\gamma$ follows from writing
\[
(1+I\gamma)^{-1}I=I+\sum_{k=1}^\infty(-1)^k(I\gamma)^kI
\]
and observing that this element of $\Mul[[B]]$ determines $\gamma$.
\end{proof}

\begin{defi}
Given $\beta\in\Mul[[B]]$, the {\em unsymmetrized R--transform} of $\beta$ is
\begin{equation}\label{eq:Rtdef}
\Rt_\beta=((1+\beta I)^{-1}\beta)\circ(I+I\beta I)\cinv\in\Mul[[B]].
\end{equation}
\end{defi}

The following is a direct application of Lemma~\ref{lem:ag}.
\begin{prop}\label{prop:Rtchar}
Given $\beta\in\Mul[[B]]$, $\Rt_\beta$ is the unique element of $\Mul[[B]]$ 
satisfying
\begin{equation}\label{eq:gammasat}
(I+I\beta I)\cinv=(1+I\,\Rt_\beta)^{-1}I.
\end{equation}
\end{prop}

\begin{prop}\label{prop:Rtcrv}
Let $\alpha\in\Mul[[B]]$ and consider the additive canonical random variable
\[
X=L_1+V_1(\alpha)\in\Lc(\Qc).
\]
Then $\Rt_{\Phit_X}=\alpha$.
\end{prop}
\begin{proof}
By Proposition~\ref{prop:Rtchar}, we must show
\begin{equation}\label{eq:IPhiI}
(I+I\Phit_XI)\cinv=(1+I\alpha)^{-1}I.
\end{equation}
Let $\beta=\sum_{k=1}^\infty(-1)^k(I\alpha)^kI$, so that, by Proposition~\ref{prop:1+a},
\[
I+\beta=(1+I\alpha)^{-1}I=I(1+\alpha I)^{-1}.
\]
We will show
\begin{equation}\label{eq:Phibeta}
\Phit_X\circ(I+\beta)=\alpha+\alpha I\alpha.
\end{equation}
This will prove~\eqref{eq:IPhiI}, because we will have
\[
\begin{split}
(I+I\Phit_XI)\circ(I+\beta)&=I+\beta+(I+\beta)(\Phit_X\circ(I+\beta))(I+\beta) \\
&=I+\beta+(I+\beta)\alpha(1+I\alpha)(I+\beta) \\
&=I+\beta+(I+\beta)\alpha I=I+I\alpha I+\beta(1+\alpha I) \\
&=I+I\alpha I+(I(1+\alpha I)^{-1}-I)(1+\alpha I)
=I.
\end{split}
\]
For $n\ge2$, we have
\begin{equation}\label{eq:Ibeta}
\begin{split}
(I+\beta)_n(b_1,\ldots,b_n)&=\big(\sum_{k=1}^\infty(-1)^k(I\alpha)^kI\big)_n(b_1,\ldots,b_n) \\
&=\sum_{k=1}^{n-1}(-1)^k
\sum_{\substack{p_1,\ldots,p_k\ge1 \\ p_1+\cdots+p_k=n-1}}
\begin{aligned}[t]
&b_{q_1+1}\alpha_{p_1-1}(b_{q_1+2},\ldots,b_{q_1+p_1}) \\
&\cdot b_{q_2+1}\alpha_{p_2-1}(b_{q_2+2},\ldots,b_{q_2+p_2})\cdots \\
&\cdot b_{q_k+1}\alpha_{p_k-1}(b_{q_k+2},\ldots,b_{q_k+p_k})b_n\end{aligned} \\
&=\sum_{\tau\in\IP(n-1)}(-1)^{|\tau|}b_1\alpha_\tau[b_2,\ldots,b_{n-1}]b_n,
\end{split}
\end{equation}
where $q_1=0$ and $q_j=p_1+\cdots+p_{j-1}$.
For the last equality above we used Lemma~\ref{lem:alphaint}.

For $n\ge1$, we have
\[
(\Phit_X\circ(I+\beta))_n(b_1,\ldots,b_n)=
\sum_{k=1}^n\sum_{\substack{p_1,\ldots,p_k\ge1 \\ p_1+\cdots+p_k=n}}
\begin{aligned}[t]
\Phit_{X,k}((I+\beta)_{p_1}(b_{q_1+1},\ldots,b_{q_1+p_1}),\cdots& \\
(I+\beta)_{p_k}(b_{q_k+1},\ldots,b_{q_k+p_k})&),
\end{aligned}
\]
where $q_j=p_1+\cdots+p_{j-1}$.
Let $c_j=(I+\beta)_{p_j}(b_{q_j+1},\ldots,b_{q_j+p_j})$.
Using Proposition~\ref{prop:Xaddtve}, we have
\[
(\Phit_X\circ(I+\beta))_n(b_1,\ldots,b_n)=
\sum_{k=1}^n\sum_{\substack{p_1,\ldots,p_k\ge1 \\ p_1+\cdots+p_k=n}}
\sum_{\pi\in\NC(k+1)}\alpha_\pi[c_1,\ldots,c_k].
\]
Equation~\eqref{eq:Ibeta} shows
\[
c_j=
\begin{cases}
\sum_{\tau_j\in\IP(p_j-1)}(-1)^{|\tau_j|}b_{q_j+1}
\alpha_{\tau_j}[b_{q_j+2},\ldots,b_{q_j+p_j-1}]b_{q_j+p_j},&p_j\ge2 \\
b_{q_j+1},&p_j=1.
\end{cases}
\]
For each $j$ such that $p_j\ge2$, fix some choice of $\tau_j\in\IP(p_j-1)$.
Let
\[
d_j=
\begin{cases}
b_{q_j+1}
\alpha_{\tau_j}[b_{q_j+2},\ldots,b_{q_j+p_j-1}]b_{q_j+p_j},&p_j\ge2 \\
b_{q_j+1},&p_j=1.
\end{cases}
\]
If $j\in\{1,\ldots,k\}$ with $p_j\ge2$, then using Lemma~\ref{lem:insertNC} we get
\[
\alpha_\pi[d_1,\ldots,d_k]
=\alpha_{\pit}[d_1,\ldots d_{j-1},b_{q_j+1},\ldots,b_{q_j+p_j},d_{j+1},\ldots,d_k],
\]
where $\pit$ is obtained from $\pi$, loosely speaking, by sliding a (translated) copy
of $\tau_j$ between $j$ and $j+1$.
Applying Lemma~\ref{lem:insertNC} repeatedly and using the convention $\tau_j=\emptyset\in\IP(0)$
whenever $p_j=1$, we obtain
\[
\alpha_\pi[d_1,\ldots,d_k]=\alpha_\sigma[b_1,\ldots,b_n],
\]
where $\sigma\in\NC(n+1)$ is obtained from $\pi$, loosely speaking, by inserting an appropriately
translated copy of $\tau_j$ in between $j$ and $j+1$, for every $j\in\{1,\ldots,k\}$.
More precisely, 
\[
\sigma=\pi\wedge(\tau_1,\ldots,\tau_k):=\piar\cup\bigcup_{j=1}^k\tauar_j,
\]
where $\piar$ is the result of applying the map $j\mapsto q_j+1$ to every element of every element of $\pi$,
and $\tauar_j$ is the result of adding $q_j+1$ to every element of every element of $\tau_j$.
Therefore, we have
\begin{multline*}
(\Phit_X\circ(I+\beta))_n(b_1,\ldots,b_n)=
\sum_{k=1}^n\sum_{\substack{p_1,\ldots,p_k\ge1 \\ p_1+\cdots+p_k=n}}
\sum_{\pi\in\NC(k+1)} \\[1ex]
\quad\sum_{\tau_1\in\IP(p_1-1),\ldots,\tau_k\in\IP(p_k-1)}
(-1)^{|\tau_1|+\cdots+|\tau_k|}
\alpha_{\pi\wedge(\tau_1,\ldots,\tau_k)}[b_1,\ldots,b_n].
\end{multline*}
Hence,
\begin{equation}\label{eq:Nsig}
(\Phit_X\circ(I+\beta))_n(b_1,\ldots,b_n)=\sum_{\sigma\in\NC(n+1)}
N_\sigma\;\alpha_{\sigma}[b_1,\ldots,b_n],
\end{equation}
where
\[
N_\sigma=\sum_{(\pi,\tau_1,\ldots,\tau_k)\in\Ic_\sigma}(-1)^{|\tau_1|+\cdots+|\tau_k|}
\]
and where $\Ic_\sigma$
is the set of all $(k+1)$--tuples $(\pi,\tau_1,\ldots,\tau_k)$ for $k\in\{1,\ldots,n\}$
with $\pi\in\NC(k+1)$ and $\tau_1\in\IP(p_1-1),\ldots,\tau_k\in\IP(p_k-1)$
for some  $p_1,\ldots,p_k\ge1$ such that $p_1+\cdots+p_k=n$,
such that $\sigma=\pi\wedge(\tau_1,\ldots,\tau_k)$.

We now set about calculating $N_\sigma$ for $\sigma\in\NC(n+1)$.
By {\em interval block} of $\sigma$, we will mean a block (i.e.\ element) of $\sigma$
that is an interval, i.e.\ is of the form $\{m,m+1,\ldots,m+\ell-1\}$.
We will say an interval block is {\em internal} if it has neither $1$ nor $n+1$ as element.
If we have $\sigma=\pi\wedge(\tau_1,\ldots,\tau_k)$, then all elements of $\tau_j$
become internal interval blocks of $\sigma$.
Let $\Bc_\sigma$ be the set of all internal interval blocks of $\sigma$.
We now describe a map from the power set of $\Bc_\sigma$ into $\Ic_\sigma$.
Let $\Sc\subseteq\Bc_\sigma$.
Let $A=\{1,\ldots,n+1\}\backslash(\bigcup\Sc)$ and let $k=|A|-1$.
Since $1,n+1\in A$, we have $k\ge1$.
Let us write $A=\{r_1,r_2,\ldots,r_{k+1}\}$, where
$1=r_1<r_2<\cdots<r_{k+1}=n+1$.
Let $\Sc_j$ be the set of those elements of $\Sc$ that lie between $r_j$ and $r_{j+1}$.
Then each $\Sc_j$ is either empty or is a set of intervals that can be laid end--to--end
without gaps.
For each $j\in\{1,\ldots,k\}$, let $\tau_j\in\IP(r_{j+1}-r_j-1)$ be the appropriate
translation of $\Sc_j$, i.e.\ 
the result of subtracting $r_j-1$ from every element of every element of $\Sc_j$.
Let $\pi$ be the result of applying the map $r_j\mapsto j$ to every element of every
element of $\sigma\backslash(\bigcup\Sc)$.
Then $\sigma=\pi\wedge(\tau_1,\ldots,\tau_k)$, and $(\pi,\tau_1,\ldots,\tau_k)\in\Ic_\sigma$.
The map
\begin{equation}\label{eq:Sc}
\Sc\mapsto(\pi,\tau_1,\ldots,\tau_k)
\end{equation}
that we have described has as inverse the map that sends $(\pi,\tau_1,\ldots,\tau_k)$
to the set of interval blocks of $\sigma$ that arise from blocks of the $\tau_j$ in the
realization of $\sigma$ as $\pi\wedge(\tau_1,\ldots,\tau_k)$.
Thus, the map~\eqref{eq:Sc} is a bijection from the power set of $\Bc_\sigma$ onto $\Ic_\sigma$.
Moreover, if $(\pi,\tau_1,\ldots,\tau_k)\in\Ic_\sigma$ corresponds to $\Sc\in\Bc_\sigma$,
then
\[
(-1)^{|\tau_1|+\cdots+|\tau_k|}=(-1)^{|\Sc|}.
\]
Therefore, 
\[
N_\sigma=\sum_{\Sc\subseteq\Bc_\sigma}(-1)^{|\Sc|}=
\begin{cases}
1,&\Bc_\sigma=\emptyset \\
0,&\Bc_\sigma\ne\emptyset.
\end{cases}
\]
Since every noncrossing partition of $\{1,\ldots,n+1\}$ has an interval block,
the only elements $\sigma\in\NC(n+1)$ without internal interval blocks
are $1_{n+1}$, which has exactly one element, and those of the form
\[
\sigma=\eta_p:=\{\{1,\ldots,p\},\{p+1,\ldots,n+1\}\}
\]
for $1\le p\le n$.
Since
\begin{align*}
\alpha_{1_{n+1}}[b_1,\ldots,b_n]&=\alpha_n(b_1,\ldots,b_n) \\
\alpha_{\eta_p}[b_1,\ldots,b_n]&=\alpha_{p-1}(b_1,\ldots,b_{p-1})b_p\alpha_{n-p}(b_{p+1},\ldots,b_n),
\end{align*}
from~\eqref{eq:Nsig} we get, for every $n\ge1$,
\[
\begin{split}
(\Phit_X\circ(I+\beta))_n(b_1,&\ldots,b_n)= \\
&=\alpha_n(b_1,\ldots,b_n)
+\sum_{p=1}^n\alpha_{p-1}(b_1,\ldots,b_{p-1})b_p\alpha_{n-p}(b_{p+1},\ldots,b_n) \\
&=(\alpha+\alpha I\alpha)_n(b_1,\ldots,b_n).
\end{split}
\]
Of course, we also have
\[
(\Phit_X\circ(I+\beta))_0=\Phit_{X,0}=\Ec(X)=\alpha_0=(\alpha+\alpha I\alpha)_0.
\]
We have proved~\eqref{eq:Phibeta}, and the proof of the proposition is complete.
\end{proof}

\begin{thm}\label{thm:Radditive}
Let $x$ and $y$ be free random variables in any $B$--valued noncommutative probability space.
Then
\[
\Rt_{\Phit_{x+y}}=\Rt_{\Phit_x}+\Rt_{\Phit_y}.
\]
\end{thm}
\begin{proof}
By Proposition~\ref{prop:can}, we can find additive
canonical random variables $X,Y\in\Lc(\Qc)$ as in~\eqref{eq:Xagain} and~\eqref{eq:Y}
such that $\Phit_X=\Phit_x$ and $\Phit_Y=\Phit_y$.
By Proposition~\ref{prop:free}, $X$ and $Y$ are free.
Hence, $\Phit_{x+y}=\Phit_{X+Y}$.
By Proposition~\ref{prop:Rtcrv}, $\Rt_{\Phit_X}=\alpha$ and $\Rt_{\Phit_Y}=\beta$, where
$\alpha$ and $\beta$ as are they appear in~\eqref{eq:Xagain} and~\eqref{eq:Y}.
On the other hand, by Lemma~\ref{lem:X+Y} and Proposition~\ref{prop:Rtcrv} again,
$\Rt_{\Phit_{X+Y}}=\alpha+\beta$.
\end{proof}

\section{The unsymmetrized T--transform}
\label{sec:Strans}

In this section, we relate the multiplicative canonical random variables of~\S\ref{sec:crv}
to the unsymmetrized T--transform defined at~\eqref{eq:Tt}
(and more fomally below), and use them to show twisted mulitplicativity~\eqref{eq:Tttmult}
under hypothesis of freeness.

\begin{defi}\label{def:alpha<>}
Let $\alpha\in\Mul[[B]]$ with $\alpha_0$ invertible.
Let $n\in\Nats_0$ and let $\pi\in\NCL(n+1)$.
We will define an element
\begin{equation}\label{eq:alpha<>}
\alpha_\pi\langle b_1,\ldots,b_n\rangle
\end{equation}
of $B$, for $b_1,\ldots,b_n\in B$.
Let $J=\{m,\ldots,m+\ell-1\}$ be the right--most interval block of $\pi$.
If $m=1$, which implies $\pi=1_{n+1}$, then we set
\begin{equation}\label{eq:alpha<>1}
\alpha_\pi\langle b_1,\ldots,b_n\rangle=
\alpha_n(b_1\alpha_0,b_2\alpha_0,\ldots,b_n\alpha_0).
\end{equation}
Suppose $m>1$.
Suppose $m$ is singly covered by $\pi$ and
let $\pi'=\pi\contr_{J^c}$ be the renumbered restriction of $\pi$ to the complement of $J$.
Then $\pi'\in\NCL(n-\ell+1)$.
If $J$ is all the way to the right, i.e.\ if $m+\ell-1=n+1$, then we set
\begin{equation}\label{eq:alpha<>2}
\alpha_\pi\langle b_1,\ldots,b_n\rangle=
\alpha_{\pi'}\langle b_1,\ldots,b_{m-2}\rangle b_{m-1}\alpha_{\ell-1}(b_m\alpha_0,\ldots,b_n\alpha_0),
\end{equation}
while if $m+\ell-1\le n$, then we set
\begin{equation}\label{eq:alpha<>3}
\begin{split}
\alpha_\pi\langle&b_1,\ldots,b_n\rangle= \\
&=\alpha_{\pi'}\langle b_1,\ldots,b_{m-2},
b_{m-1}\alpha_{\ell-1}(b_m\alpha_0,\ldots,b_{m+\ell-2}\alpha_0)b_{m+\ell-1},b_{m+\ell},\ldots,b_n\rangle.
\end{split}
\end{equation}
Now suppose $m$ is doubly covered by $\pi$.
Then we must have $\ell\ge2$.
Let $\pi''=(\pi\backslash\{J\})\contr_{J^c\cup\{m\}}$
be the renumbered restriction of $\pi\backslash\{J\}$ to $\{1,\ldots,m\}\cup\{m+\ell,\ldots,n+1\}$.
Then $\pi''\in\NCL(n-\ell+2)$, and we set
\begin{equation}\label{eq:alpha<>4}
\begin{aligned}
\alpha&_\pi\langle b_1,\ldots,b_n\rangle= \\
&=\alpha_{\pi''}\langle b_1,\ldots,b_{m-2},
b_{m-1}\alpha_{\ell-1}(b_m\alpha_0,\ldots,b_{m+\ell-2}\alpha_0)\alpha_0^{-1},
b_{m+\ell-1},b_{m+\ell},\ldots,b_n\rangle.
\end{aligned}
\end{equation}
Equations~\eqref{eq:alpha<>1}, \eqref{eq:alpha<>2}, \eqref{eq:alpha<>3} and~\eqref{eq:alpha<>4}
recursively define the quantity~\eqref{eq:alpha<>}.
\end{defi}
Some examples of $\alpha_\pi\langle b_1,\ldots,b_n\rangle$
are provided in Table~\ref{tab:alpha<>}.
\begin{table}[b]
\caption{Some examples of $\alpha_\pi\langle b_1,\ldots,b_n\rangle$.}
\label{tab:alpha<>}
\begin{tabular}{r|c|l}
$\pi$\hspace*{2ex} & $\Gc_\pi$ & $\alpha_\pi\langle b_1,\ldots,b_n\rangle$ \\ \hline\hline
$(1,2,3,4)$ &
\begin{picture}(46,18)(-3,0) 
 \multiput(0,0)(12,0){4}{\circle*{3}}
 \multiput(0,0)(12,0){4}{\line(0,1){12}}
 \put(0,12){\line(1,0){36}}
\end{picture} 
& $\alpha_3(b_1\alpha_0,b_2\alpha_0,b_3\alpha_0)$ \\ \hline
$(1,2,4)(3) $ &
\begin{picture}(46,18)(-3,0) 
 \multiput(0,0)(12,0){4}{\circle*{3}}
 \multiput(0,0)(12,0){2}{\line(0,1){12}}
 \put(36,0){\line(0,1){12}}
 \put(0,12){\line(1,0){36}}
\end{picture}
& $\alpha_2(b_1\alpha_0,b_2\alpha_0b_3\alpha_0)$ \\ \hline
$(1,3,4)(2) $ &
\begin{picture}(46,18)(-3,0) 
 \multiput(0,0)(12,0){4}{\circle*{3}}
 \put(0,0){\line(0,1){12}}
 \multiput(24,0)(12,0){2}{\line(0,1){12}}
 \put(0,12){\line(1,0){36}}
\end{picture}
& $\alpha_2(b_1\alpha_0b_2\alpha_0,b_3\alpha_0)$ \\ \hline
$(1,2)(3,4)$ &
\begin{picture}(46,18)(-3,0) 
 \multiput(0,0)(12,0){4}{\circle*{3}}
 \multiput(0,0)(12,0){4}{\line(0,1){12}}
 \multiput(0,12)(24,0){2}{\line(1,0){12}}
\end{picture}
& $\alpha_1(b_1\alpha_0)b_2\alpha_1(b_3\alpha_0)$ \\ \hline
$(1,4)(2,3)$ &
\begin{picture}(46,18)(-3,0) 
 \multiput(0,0)(12,0){4}{\circle*{3}}
 \multiput(0,0)(36,0){2}{\line(0,1){12}}
 \multiput(12,0)(12,0){2}{\line(0,1){8}}
 \put(0,12){\line(1,0){36}}
 \put(12,8){\line(1,0){12}}
\end{picture}
& $\alpha_1(b_1\alpha_1(b_2\alpha_0)b_3\alpha_0)$ \\ \hline
$(1)(2)(3)(4)$ &
\begin{picture}(46,18)(-3,0) 
 \multiput(0,0)(12,0){4}{\circle*{3}}
\end{picture}
& $\alpha_0b_1\alpha_0b_2\alpha_0b_3\alpha_0$ \\ \hline
$(1,2)(2,3,4) $ &
\begin{picture}(46,18)(-3,0) 
 \multiput(0,0)(12,0){4}{\circle*{3}}
 \multiput(0,0)(12,0){4}{\line(0,1){12}}
 \put(12,0){\line(1,1){12}}
 \multiput(0,12)(24,0){2}{\line(1,0){12}}
\end{picture}
& $\alpha_1(b_1\alpha_2(b_2\alpha_0,b_3\alpha_0))$ \\ \hline
$(1,2)(2,3)(3,4)$ &
\begin{picture}(46,18)(-3,0) 
 \multiput(0,0)(12,0){4}{\circle*{3}}
 \multiput(0,0)(12,0){4}{\line(0,1){12}}
 \put(0,12){\line(1,0){12}}
 \multiput(12,0)(12,0){2}{\line(1,1){12}}
\end{picture}
& $\alpha_1(b_1\alpha_1(b_2\alpha_1(b_3\alpha_0)))$ \\ \hline
$(1,2,4)(2,3) $ &
\begin{picture}(46,18)(-3,0) 
 \multiput(0,0)(12,0){4}{\circle*{3}}
 \multiput(0,0)(12,0){2}{\line(0,1){12}}
 \put(36,0){\line(0,1){12}}
 \put(0,12){\line(1,0){36}}
 \put(12,0){\line(3,2){12}}
 \put(24,0){\line(0,1){8}}
\end{picture}
& $\alpha_2(b_1\alpha_1(b_2\alpha_0),b_3\alpha_0)$
\end{tabular}
\end{table}

\begin{lemma}\label{lem:ncnp}
Let $\alpha\in\Mul[[B]]$ with $\alpha_0$ invertible and for $s\in\Nats_0\cup\Nats_0^*$ let
\begin{equation}\label{eq:ys}
y_s=
\begin{cases}
V_{1,k}(\alpha_k),&s=k\in\Nats_0 \\
W_{1,k}(\alpha_k),&s=k^*\in\Nats_0^*
\end{cases}
\end{equation}
in $\Lc(\Qc)$.
Let $n\in\Nats_0$, let $s:\{1,\ldots,n+1\}\to\Nats_0\cup\Nats_0^*$
and let $b_1,\ldots,b_n\in B$.
If there is $\pi\in\NCL(n+1)$ such that $S_\pi=s$, where $S_\pi$ is as in Definition~\ref{def:Spi}, then
\[ 
y_{s(1)}b_1y_{s(2)}b_2\cdots y_{s(n)}b_ny_{s(n+1)}=\alpha_\pi\langle b_1,\ldots,b_n\rangle.
\] 
If there is no $\pi\in\NCL(n+1)$ such that $S_\pi=s$, then
\begin{equation}\label{eq:Ec0ys}
\Ec(y_{s(1)}b_1y_{s(2)}b_2\cdots y_{s(n)}b_ny_{s(n+1)})=0.
\end{equation}
\end{lemma}
\begin{proof}
Let $N$ be the number of $j\in\{1,\ldots,n+1\}$ such that $s(j)\ne0^*$.
We will proceed by induction on $N$.
If $N=0$, then there is no $\pi\in\NCL(n+1)$ such that $S_\pi=s$ and,
since $P\alpha_0L_1=0$,~\eqref{eq:Ec0ys} holds.

Suppose $N=1$.
The only way in which we can have $s=S_\pi$ for $\pi\in\NCL(n+1)$ is if $\pi=1_{n+1}$
and $s(1)=n$, $s(j)=0^*$ for all $j\in\{2,\ldots,n+1\}$.
In this case,
\begin{multline*}
y_{s(1)}b_1\cdots y_{s(n)}b_ny_{s(n+1)}
=V_{1,n}(\alpha_n)b_1\alpha_0L_1b_2\alpha_0L_1\cdots b_n\alpha_0L_1= \\
=\alpha_n(b_1\alpha_0,b_2\alpha_0,\ldots,b_n\alpha_0)
=\alpha_{1_{n+1}}\langle b_1,\ldots,b_n\rangle.
\end{multline*}
On the other hand, if $s\ne S_{1_{n+1}}$, then we easily see that~\eqref{eq:Ec0ys} holds.

Suppose $N\ge2$.
Let $m$ be largest such that $s(m)\ne0^*$.
Then $m\ge2$.
Suppose $s(m)=k\in\Nats_0$.
If $\pi\in\NCL(n+1)$ is such that $S_\pi=s$, then we must have $k\le n-m+1$ and
$\{m,\ldots,m+k\}\in\pi$.
Therefore, if $k>n-m+1$, then there is no $\pi\in\NCL(n+1)$ such that $S_\pi=s$ and we
also have
\[
y_{s(m)}b_m\cdots y_{s(n)}b_ny_{s(n+1)}\Omega
=V_{1,k}(\alpha_k)b_m\alpha_0\delta_1\otdt b_n\alpha_0\delta_1=0,
\]
so~\eqref{eq:Ec0ys} holds.
Suppose now $k\le n-m+1$.
Let
\[
s'(j)=
\begin{cases}
s(j),&1\le j\le m-1 \\
s(j+k+1),&m\le j \le n-k.
\end{cases}
\]
Then there is $\pi\in\NC(n+1)$ such that $S_\pi=s$ if and only if there is
$\pi'\in\NC(n-k)$ such that $S_{\pi'}=s'$,
and in this case, $\pi$ is obtained from $\pi'$ by
applying the map
\[
j\mapsto
\begin{cases}
j,&1\le j\le m-1 \\
j+k+1,&m\le j \le n-k
\end{cases}
\]
to all elements of all elements of $\pi'$, and then adding the element $\{m,\ldots,m+k\}$.
Moreover, we have
\begin{equation}\label{eq:ysm1}
y_{s(m)}b_m\cdots y_{s(m+k-1)}b_{m+k-1}y_{s(m+k)}
=\alpha_k(b_m\alpha_0,\ldots,b_{m+k-1}\alpha_0).
\end{equation}
If there is $\pi\in\NCL(n+1)$ such that $S_\pi=s$, then it arises from some $\pi'$ as above, and
in the case $k<n-m+1$, by the induction hypothesis and using Definition~\ref{def:alpha<>} we have
\begin{align*}
y_{s(1)}b_1\cdots y_{s(n)}b_ny_{s(n+1)}
&=
\begin{aligned}[t]
 &y_{s(1)}b_1\cdots y_{s(m-2)}b_{m-2} \\
 &\quad\cdot y_{s(m-1)}b_{m-1}\alpha_k(b_m\alpha_0,\ldots,b_{m+k-1}\alpha_0)b_{m+k} \\
 &\quad\cdot y_{s(m+k+1)}b_{m+k+1}\cdots y_{s(n)}b_ny_{s(n+1)}
\end{aligned} \\
&=
\begin{aligned}[t]
 \alpha_{\pi'}\langle&b_1,\ldots,b_{m-2}, \\
 &b_{m-1}\alpha_k(b_m\alpha_0,\ldots,b_{m+k-1}\alpha_0)b_{m+k}, \\
 &b_{m+k+1},\ldots,b_n\rangle
\end{aligned} \\
&=\alpha_\pi\langle b_1,\ldots,b_n\rangle.
\end{align*}
A similar calculation can be made in the case $m+k=n+1$.
If there is no $\pi\in\NCL(n+1)$ such that $S_\pi=s$, then there is no $\pi'$
such that $S_{\pi'}=s'$ and by the induction hypothesis,
\[ 
\Ec(y_{s(1)}b_1\cdots y_{s(n)}b_ny_{s(n+1)})=
\begin{aligned}[t]
 \Ec(&y_{s(1)}b_1\cdots y_{s(m-2)}b_{m-2} \\
 &\quad\cdot y_{s(m-1)}b_{m-1}\alpha_k(b_m\alpha_0,\ldots,b_{m+k-1}\alpha_0)b_{m+k} \\
 &\quad\cdot y_{s(m+k+1)}b_{m+k+1}\cdots y_{s(n)}b_ny_{s(n+1)})=0.
\end{aligned}
\]
Now suppose $s(m)=k^*\in\Nats_0^*$.
Then $k\ge1$.
If $\pi\in\NCL(n+1)$ is such that $S_\pi=s$, then we must have
$\{m,\ldots,m+k\}\in\pi$ and $k\le n-m+1$.
Therefore, if $k>n-m+1$, then there is no $\pi\in\NCL(n+1)$ such that $S_\pi=s$ and we
also have
\[
y_{s(m)}b_m\cdots y_{s(n)}b_ny_{s(n+1)}\Omega
=W_{1,k}(\alpha_k)b_m\alpha_0\delta_1\otdt b_n\alpha_0\delta_1=0,
\]
so~\eqref{eq:Ec0ys} holds.
Suppose now $k\le n-m+1$.
Let
\[
s'(j)=
\begin{cases}
s(j),&1\le j\le m-1 \\
0^*,&j=m \\
s(j+k),&m+1\le j \le n-k+1.
\end{cases}
\]
Then there is $\pi\in\NC(n+1)$ such that $S_\pi=s$ if and only if there is
$\pi'\in\NC(n-k+1)$ such that $S_{\pi'}=s'$,
and in this case, $\pi$ is obtained from $\pi'$ by
applying the map
\[
j\mapsto
\begin{cases}
j,&1\le j\le m \\
j+k,&m+1\le j \le n-k+1
\end{cases}
\]
to all elements of all elements of $\pi'$, and then adding the element $\{m,\ldots,m+k\}$.
Moreover, we have
\begin{equation}\label{eq:ysm2}
y_{s(m)}b_m\cdots y_{s(m+k-1)}b_{m+k-1}y_{s(m+k)}
=\alpha_k(b_m\alpha_0,\ldots,b_{m+k-1}\alpha_0)L_1.
\end{equation}
If there is $\pi\in\NCL(n+1)$ such that $S_\pi=s$, then it arises from some $\pi'$ as above, and
we have
\[
y_{s(1)}b_1\cdots y_{s(n)}b_ny_{s(n+1)}=
\begin{aligned}[t]
 &y_{s(1)}b_1\cdots y_{s(m-2)}b_{m-2} \\
 &\quad\cdot y_{s(m-1)}b_{m-1}\alpha_k(b_m\alpha_0,\ldots,b_{m+k-1}\alpha_0)L_1b_{m+k} \\
 &\quad\cdot y_{s(m+k+1)}b_{m+k+1}\cdots y_{s(n)}b_ny_{s(n+1)}.
\end{aligned}
\]
But
\[
b_{m-1}\alpha_k(b_m\alpha_0,\ldots,b_{m+k-1}\alpha_0)L_1=
b_{m-1}\alpha_k(b_m\alpha_0,\ldots,b_{m+k-1}\alpha_0)\alpha_0^{-1}W_{1,0}(\alpha_0).
\]
Therefore, using the induction hypothesis and Definition~\ref{def:alpha<>} we have
\begin{align*}
y_{s(1)}b_1\cdots y_{s(n)}b_ny_{s(n+1)}&=
\begin{aligned}[t]
 \alpha_{\pi'}\langle&b_1,\ldots,b_{m-2}, \\
 &b_{m-1}\alpha_k(b_m\alpha_0,\ldots,b_{m+k-1}\alpha_0)\alpha_0^{-1}, \\
 &b_{m+k},\ldots,b_n\rangle
\end{aligned} \\
&=\alpha_\pi\langle b_1,\ldots,b_n\rangle.
\end{align*}
If there is no $\pi\in\NCL(n+1)$ such that $S_\pi=s$, then there is no $\pi'$
such that $S_{\pi'}=s'$ and using the induction hypothesis we get
\[ 
\Ec(y_{s(1)}b_1\cdots y_{s(n)}b_ny_{s(n+1)})=0.
\]
This completes the induction step.
\end{proof}

\begin{prop}\label{prop:Xmultve}
Let $\alpha\in\Mul[[B]]$ and let 
\[
X=V_1(\alpha)+W_1(\alpha).
\]
Then for every $n\ge0$ and $b_1,\ldots,b_n\in B$,
\begin{equation}\label{eq:ncnp}
\Phit_{X,n}(b_1,\ldots,b_n)=\sum_{\pi\in\NCL(n+1)}\alpha_\pi\langle b_1,\ldots,b_n\rangle.
\end{equation}
\end{prop}
\begin{proof}
We have
\begin{align}
\Phit_{X,n}(b_1,\ldots,b_n)&=\Ec(Xb_1Xb_2\cdots Xb_nX) \notag \\[1ex]
&=\sum_{s(1),\ldots,s(n+1)\in\Nats_0\cup\Nats_0^*}
\Ec(y_{s(1)}b_1y_{s(2)}b_2\cdots y_{s(n)}b_ny_{s(n+1)}), \label{eq:yksum}
\end{align}
where $y_{s(j)}$ is as defined in~\eqref{eq:ys};
by Lemma~\ref{lem:ncnp},
only finitely many terms of the sum~\eqref{eq:yksum} are nonzero, and the equality~\eqref{eq:ncnp} holds.
\end{proof}

\begin{lemma}\label{lem:insertNCNP}
Let $n\in\Nats$, 
let $\pi\in\NCL(n+1)$ and suppose 
\[
J=\{m,\ldots,m+\ell-1\}\subseteq\{1,\ldots,n+1\}
\]
is a proper subinterval that splits $\pi$.
Let $\pi''$ be the renumbered restriction $\pi$ to $J$, and let $\pi'$
be the renumbered restriction of $\pi$ to the complement of $J$.
Then $\pi'\in\NCL(n-\ell+1)$ and $\pi''\in\NCL(\ell)$.
Moreover, for all $\alpha\in\Mul[[B]]$ and all $b_1,\ldots,b_n\in B$, we have
\begin{multline*}
\alpha_\pi\langle b_1,\ldots,b_n\rangle= \\
=\begin{cases}
\alpha_{\pi''}\langle b_1,\ldots,b_{\ell-1}\rangle b_\ell\alpha_{\pi'}\langle b_{\ell+1},\ldots,b_n\rangle,&m=1 \\
\begin{aligned}[b]
\alpha_{\pi'}\langle b_1,&\ldots,b_{m-2}, \\
&b_{m-1}\alpha_{\pi''}\langle b_m,\ldots,b_{m+\ell-2}\rangle b_{m+\ell-1},
b_{m+\ell},\ldots,b_n\rangle,\end{aligned}&2\le m\le n-\ell+1 \\
\alpha_{\pi'}\langle b_1,\ldots,b_{m-2}\rangle b_{m-1}\alpha_{\pi''}\langle b_m,\ldots,b_n\rangle,&m=n-\ell+2.
\end{cases}
\end{multline*}
\end{lemma}
\begin{proof}
This can be proved using Lemma~\ref{lem:ncnp} just as Lemma~\ref{lem:insertNC} was proved using Lemma~\ref{lem:ncp}.
\end{proof}

\begin{defi}
For $n\in\Nats$, let $\NCLo(n)$ be the set of all $\pi\in\NCL(n)$ such that $1$ and $n$ belong
to the same block of the noncrossing partition $\pihat$ generated by $\pi$.
\end{defi}

\begin{lemma}\label{lem:sL1}
Let $n\in\Nats_0$ and let $\pi\in\NCLo(n+1)$.
Then $k:=S_\pi(1)\in\Nats$.
Let $s:\{1,\ldots,n+1\}\to\Nats_0\cup\Nats_0^*$ be
\[
s(j)=
\begin{cases}
k^*,&j=1 \\
S_\pi(j),&2\le j\le n+1.
\end{cases}
\]
Then for all $\alpha\in\Mul[[B]]$ and all $b_1,\ldots,b_n\in B$,
\begin{equation}\label{eq:yL1}
y_{s(1)}b_1\cdots y_{s(n)}b_ny_{s(n+1)}=\alpha_\pi\langle b_1,\ldots,b_n\rangle L_1,
\end{equation}
where $y_{s(j)}$ is as defined in~\eqref{eq:ys}.
\end{lemma}
\begin{proof}
We use induction on $|\pi|$.
If $|\pi|=1$, then $\pi=1_{n+1}$ and~\eqref{eq:yL1} follows directly.
Suppose $|\pi|\ge2$.
Let $J=\{m,\ldots,m+\ell-1\}$ be the right--most interval block of $\pi$.
Then $m\ge2$.
Suppose $m+\ell-1=n+1$.
Because $\pi\in\NCLo(n+1)$, $m$ is doubly covered by $\pi$.
Moreover, we have
\[
y_{s(m)}b_m\cdots y_{s(n)}b_ny_{s(n+1)}=\alpha_{\ell-1}(b_m,\ldots,b_n)L_1,
\]
so letting $\bt=b_{m-1}\alpha_{\ell-1}(b_m,\ldots,b_n)\alpha_0^{-1}$, we have
\[
y_{s(1)}b_1\cdots y_{s(n)}b_ny_{s(n+1)}=y_{s'(1)}b_1\cdots y_{s'(m-2)}b_{m-2}y_{s'(m-1)}\bt y_{s'(m)},
\]
where $s':\{1,\ldots,m\}\to\Nats_0\cup\Nats_0^*$ is
\[
s'(j)=
\begin{cases}
s(j),&1\le j\le m-1 \\
0^*,&j=m.
\end{cases}
\]
Let $\pi'$ be the renumbered restriction of $\pi\backslash\{J\}$ to $\{1,\ldots,m\}$.
Then $\pi'\in\NCLo(m)$.
Therefore,
\[
s'(j)=
\begin{cases}
k^*,&j=1 \\
S_{\pi'}(j),&2\le j\le m.
\end{cases}
\]
Using the induction hypothesis and Definition~\ref{def:alpha<>}, 
we have
\[
y_{s'(1)}b_1\cdots y_{s'(m-2)}b_{m-2}y_{s'(m-1)}\bt y_{s'(m)}
=\alpha_{\pi'}\langle b_1,\ldots,b_{m-2},\bt\rangle L_1
=\alpha_\pi\langle b_1,\ldots,b_n\rangle L_1,
\]
as required.

Now suppose $m+\ell-1\le n$.
If $m$ is singly covered by $\pi$, then let $\pi'$ be the renumbered restriction of $\pi$ to $J^c$
and let
\[
\bt=b_{m-1}\alpha_{\ell-1}(b_m\alpha_0,\ldots,b_{m+\ell-2}\alpha_0)b_{m+\ell-1}.
\]
Then $\pi'\in\NCLo(n-\ell+1)$.
Letting
\[
s'(j)=
\begin{cases}
s(1),&j=1 \\
S_{\pi'}(j),&2\le j\le n-\ell+1
\end{cases}
\]
and invoking the induction hypothesis and Definition~\ref{def:alpha<>}, we have
\begin{align*}
y_{s(1)}b_1&\cdots y_{s(n)}b_ny_{s(n+1)}= \\
&=y_{s'(1)}b_1\cdots y_{s'(m-2)}b_{m-2}y_{s'(m-1)}\bt y_{s'(m)}b_{m+\ell}\cdots y_{s'(n-\ell)}b_ny_{s'(n-\ell+1)} \\
&=\alpha_{\pi'}\langle b_1,\ldots,b_{m-2},\bt,b_{m+\ell},\ldots,b_n\rangle L_1
=\alpha_\pi\langle b_1,\ldots,b_n\rangle L_1,
\end{align*}
as required.
If $m$ is doubly covered by $\pi$, then letting
$\pi'$ be the renumbered restriction of $\pi\backslash\{J\}$ to $J^c\cup\{m\}$, we proceed in a similar fashion
to show~\eqref{eq:yL1}.
\end{proof}

\begin{defi}\label{def:Dc}
Suppose $\pi\in\NCL(n)$ and let $F\in\pi$ be the block such that $1\in F$.
We write
\[
F=\{\ell_0,\ell_1,\ldots,\ell_k\}
\]
with $k\ge0$ and $1=\ell_0<\ell_1<\cdots<\ell_k$.
Let $\pihat$ be the noncrossing partition generated by $\pi$.
Let $\Fh\in\pihat$ be such that $F\subseteq\Fh$.
Let
\[
r(j)=\max\{i\in\Fh\mid i<\ell_j\},\quad(1\le j\le k)
\]
and let $r(k+1)=\max(\Fh)$.
We must have $r(1)=1$ and $r(j)<\ell_j\le r(j+1)$ whenever $1\le j\le k$.
For $1\le j\le k$, let
\[
\pi_j=
\begin{cases}
\pi\contr_{[r(j)+1,r(j+1)]},&\ell_j\text{ singly covered by }\pi \\[0.5ex]
(\pi\backslash\{F\})\contr_{[r(j)+1,r(j+1)]},&\ell_j\text{ doubly covered by }\pi.
\end{cases}
\]
If $\pi\in\NCLo(n)$, then $k\ge1$ and $r(k+1)=n$, and we define
\[
\Dco_n(\pi)=(\pi_1,\ldots,\pi_k).
\]
If $\pi\in\NCL(n)\backslash\NCLo(n)$, then $r(k+1)<n$, we let
\[
\sigma=\pi\contr_{[r(k+1)+1,n]}
\]
and we define
\[
\Dc_n(\pi)=(\pi_1,\ldots,\pi_k,\sigma).
\]
\end{defi}

\begin{example}\label{ex:pi18}
Let
\[
\pi=(1,5,10,13)(2,4)(3)(5,7)(6)(8,9)(11,12)(13,16)(14,15)(17,18),
\]
whose graphical representation is drawn in Figure~\ref{fig:pi18}.
\begin{figure}[ht]
\caption{Graphical representation of $\pi$ from Example~\ref{ex:pi18}.}
\label{fig:pi18}
\begin{picture}(210,32)(-3,-14)
\multiput(0,0)(12,0){18}{\circle*{3}}
\put(0,0){\line(0,1){12}}
\put(48,0){\line(0,1){12}}
\put(108,0){\line(0,1){12}}
\put(144,0){\line(0,1){12}}
\multiput(180,0)(12,0){3}{\line(0,1){12}}
\put(0,12){\line(1,0){144}}
\put(144,0){\line(1,1){12}}
\put(156,12){\line(1,0){24}}
\put(192,12){\line(1,0){12}}
\multiput(12,0)(24,0){2}{\line(0,1){8}}
\put(12,8){\line(1,0){24}}
\put(48,0){\line(3,2){12}}
\put(60,8){\line(1,0){12}}
\multiput(72,0)(12,0){3}{\line(0,1){8}}
\put(84,8){\line(1,0){12}}
\multiput(120,0)(12,0){2}{\line(0,1){8}}
\put(120,8){\line(1,0){12}}
\multiput(156,0)(12,0){2}{\line(0,1){8}}
\put(156,8){\line(1,0){12}}
\put(-2,-14){$1$}
\put(46,-14){$5$}
\put(102,-14){$10$}
\put(138,-14){$13$}
\end{picture}
\end{figure}
Then $\Fh=\{1,5,7,10,13,16\}$, we have
\[
r(1)=1,\quad r(2)=7,\quad r(3)=10,\quad r(4)=16
\]
and the graphical representations of $\pi_1,\pi_2,\pi_3$ and $\sigma$ are given in Table~\ref{tab:pij}.
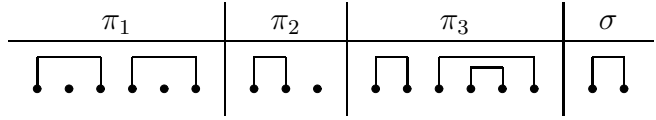
\begin{table}[ht]
\caption{Graphical representations of $\pi_1,\pi_2,\pi_3$ and $\sigma$ from Example~\ref{ex:pi18}.}
\label{tab:pij}
\begin{tabular}{c|c|c|c}
$\pi_1$ & $\pi_2$ & $\pi_3$ & $\sigma$ \\ \hline
\begin{picture}(72,24)(-6,-6)
\multiput(0,0)(12,0){6}{\circle*{3}}
\multiput(0,0)(24,0){2}{\line(0,1){12}}
\multiput(36,0)(24,0){2}{\line(0,1){12}}
\multiput(0,12)(36,0){2}{\line(1,0){24}}
\end{picture}
&
\begin{picture}(36,24)(-6,-6)
\multiput(0,0)(12,0){3}{\circle*{3}}
\multiput(0,0)(12,0){2}{\line(0,1){12}}
\put(0,12){\line(1,0){12}}
\end{picture}
&
\begin{picture}(72,24)(-6,-6)
\multiput(0,0)(12,0){6}{\circle*{3}}
\multiput(0,0)(12,0){3}{\line(0,1){12}}
\put(60,0){\line(0,1){12}}
\put(0,12){\line(1,0){12}}
\put(24,12){\line(1,0){36}}
\multiput(36,0)(12,0){2}{\line(0,1){8}}
\put(36,8){\line(1,0){12}}
\end{picture}
&
\begin{picture}(24,24)(-6,-6)
\multiput(0,0)(12,0){2}{\circle*{3}}
\multiput(0,0)(12,0){2}{\line(0,1){12}}
\put(0,12){\line(1,0){12}}
\end{picture}
\end{tabular}
\end{table}
\end{example}

\begin{lemma}\label{lem:pik}
Let $n\in\Nats$, $n\ge2$.
Then
\renewcommand{\labelenumi}{(\roman{enumi})}
\begin{enumerate}
\item
 the map $\Dco_n$ is a bijection from $\NCLo(n)$ onto
\begin{equation}\label{eq:Dcorange}
\bigcup_{k=1}^{n-1}\bigcup_{\substack{p_1,\ldots,p_k\ge1 \\ p_1+\cdots+p_k=n-1}}
\NCL(p_1)\tdt\NCL(p_k);
\end{equation}
\item
the map $\Dc_n$
is a bijection from $\NCL(n)\backslash\NCLo(n)$ onto
\begin{equation}\label{eq:Dcrange}
\bigcup_{k=0}^{n-2}\bigcup_{\substack{p_1,\ldots,p_k,p_{k+1}\ge1 \\ p_1+\cdots+p_k+p_{k+1}=n-1}}
\NCL(p_1)\tdt\NCL(p_k)\times\NCL(p_{k+1}).
\end{equation}
\end{enumerate}
\end{lemma}
\begin{proof}
It is straightforward to see that each $\pi_j$ and $\sigma$ described in Definition~\ref{def:Dc}
is a noncrossing linked partition.
Thus, $\Dco_n$ maps $\NCLo(n)$ into the set~\eqref{eq:Dcorange}
and $\Dc_n$ maps $\NCL(n)\backslash\NCLo(n)$ into the set~\eqref{eq:Dcrange}.
We will find the inverses of the maps $\Dco_n$ and $\Dc_n$.
Given $k\in\{1,\ldots,n-1\}$, $p_1,\ldots,p_k\ge1$ such that $p_1+\cdots+p_k=n-1$ and given
$\pi_j\in\NCL(p_j)$, let $\pihat_j$ be the noncrossing partition generated by $\pi_j$ and let
$E_j$ be the block of $\pihat_j$ that contains $p_j$.
Let $r_1=1$ and $r_j=1+p_1+\cdots+p_{j-1}$ ($2\le j\le k$).
Let $\ell_j=r_j+\min(E_j)$.
Let $F=\{1,\ell_1,\ell_2,\ldots,\ell_k\}$ and let
\[
\pi=\{F\}\cup\bigcup_{j=1}^k\pit_j,
\]
where $\pit_j$ is the result of adding $r_j$ to every element of every element of
\[
\begin{cases}
\pi_j,&|E_j|\ge2 \\
\pi_j\backslash\{E_j\},&|E_j|=1.
\end{cases}
\]
Then we easily see $\pi\in\NCL(n)$
and $\Dco_n(\pi)=(\pi_1,\ldots,\pi_k)$.
Also, $\pi\in\NCLo(n)$ is easily seen to be determined by $\Dco_n(\pi)$.
This proves part~(i).

For part~(ii),
let $k\ge0$, let $p_1,\ldots,p_{k+1}$ be such that $p_1+\cdots+p_{k+1}=n-1$
and let $\pi+j\in\NCL(p_j)$, ($1\le j\le k+1$).
If $k=0$ let $\tau=1_1$,
while if $k\ge1$, then let $m=n-p_{k+1}$ and $\tau=(\Dco_m)^{-1}((\pi_1,\ldots,\pi_k))$.
Let $\pi=\tau\oplus\pi_{k+1}$.
Then $\Dc_n(\pi)=(\pi_1,\ldots,\pi_{k+1})$, and we easily see that
$\pi\in\NCL(n)\backslash\NCLo(n)$ is determined by $\Dc_n(\pi)$.
This proves~(ii).
\end{proof}

\begin{lemma}\label{lem:F}
Let $n\in\Nats$ and let $\pi\in\NCL(n+1)$.
If $\pi\in\NCLo(n+1)$, then write $\Dco_n(\pi)=(\pi_1,\ldots,\pi_k)$, with $\pi_j\in\NCL(p_j)$.
If $\pi\notin\NCLo(n+1)$, then write $\Dc_n(\pi)=(\pi_1,\ldots,\pi_k,\sigma)$, with $\pi_j\in\NCL(p_j)$.
Let $r(1)=1$ and for $2\le j\le k+1$, let $r(j)=1+p_1+\cdots+p_{j-1}$.
Let $\alpha\in\Mul[[B]]$ with $\alpha_0$ invertible
and let $b_1,\ldots,b_n\in B$.
\renewcommand{\labelenumi}{(\roman{enumi})}
\begin{enumerate}
\item
If $\pi\in\NCLo(n+1)$, then 
\begin{equation}\label{eq:Fo}
\alpha_\pi\langle b_1,\ldots,b_n\rangle=
\begin{aligned}[t]
\alpha_k\big(&b_{r(1)}\alpha_{\pi_1}\langle b_{r(1)+1},\ldots,b_{r(2)-1}\rangle, \\
&b_{r(2)}\alpha_{\pi_2}\langle b_{r(2)+1},\ldots,b_{r(3)-1}\rangle, \\
&\ldots,b_{r(k)}\alpha_{\pi_k}\langle b_{r(k)+1},\ldots,b_{r(k+1)-1}\rangle\big).
\end{aligned}
\end{equation}
\item
If $\pi\notin\NCLo(n+1)$, then
\begin{equation}\label{eq:F}
\alpha_\pi\langle b_1,\ldots,b_n\rangle=
\begin{aligned}[t]
\alpha_k\big(&b_{r(1)}\alpha_{\pi_1}\langle b_{r(1)+1},\ldots,b_{r(2)-1}\rangle, \\
&b_{r(2)}\alpha_{\pi_2}\langle b_{r(2)+1},\ldots,b_{r(3)-1}\rangle, \\
&\ldots,b_{r(k)}\alpha_{\pi_k}\langle b_{r(k)+1},\ldots,b_{r(k+1)-1}\rangle\big) \\[0.3ex]
&\qquad\qquad\cdot b_{r(k+1)}\alpha_\sigma\langle b_{r(k+1)+1},\ldots,b_n\rangle.
\end{aligned}
\end{equation}
\end{enumerate}
\end{lemma}
\begin{proof}
If $\pi\notin\NCLo(n+1)$, then the interval $[r(k+1)+1,n+1]$ splits $\pi$, and by Lemma~\ref{lem:insertNCNP},
\[
\alpha_\pi\langle b_1,\ldots,b_n\rangle=\alpha_{\pi'}\langle b_1,\ldots,b_{r(k+1)-1}\rangle b_{r(k+1)}
\alpha_\sigma\langle b_{r(k+1)+1},\ldots,b_n\rangle,
\]
where $\pi'=\pi\contr_{[1,r(k+1)]}$.
If $k=0$, then we are finished.
If $k\ge1$, then~\eqref{eq:F} will follow from~(i).
Hence, it will suffice to prove~(i).

Assume $\pi\in\NCLo(n+1)$.
Let $F=\{1,\ell_1,\ldots,\ell_k\}\in\pi$ be the block containing $1$, with $1<\ell_1<\cdots<\ell_k$.
From the proof of Lemma~\ref{lem:pik}, $\ell_j=r(j)+\min(E_j)$, where $E_j$ is the block of
the noncrossing partition $\pihat_j$ generated by $\pi_j$ that contains $p_j$.
For every $j\in\{1,\ldots,k\}$,
$\ell_j\in[r(j)+1,r(j+1)]$.
If $\ell_j$ is an element of a block of $\pi$ other than $F$, then it is the minimal element of that block.
Therefore, the interval $[\ell_j-r(j),r(j+1)-r(j)]$ splits $\pi_j$.
Let $\rho_j=\pi_j\contr_{[\ell_j-r(j),r(j+1)-r(j)]}$ and, if $r(j)+1<\ell_j$,
let $\tau_j=\pi_j\contr_{[1,\ell_j-r(j)-1]}$.
Let $s=S_\pi$.
Note that $\ell_j$ is singly covered by $\pi_j$, so
\[
k_j:=S_{\pi_j}(\ell_j-r(j))\in\Nats_0.
\]
We have
\[
s(i)=
\begin{cases}
S_{\tau_j}(i-r(j)),&r(j)+1\le i\le\ell_j-1 \\
k_j^*,&i=\ell_j \\
S_{\rho_j}(i-\ell_j+1),&\ell_j+1\le i\le r(j+1).
\end{cases}
\]
By Lemma~\ref{lem:sL1}, if $r(j)+1=\ell_j$, then
\begin{align*}
b_{r(j)}y_{s(r(j)+1)}b_{r(j)+1}\cdots y_{s(r(j+1)-1)}&b_{r(j+1)-1}y_{s(r(j+1))}= \\
&=b_{r(j)}\alpha_{\rho_j}\langle b_{r(j)+1},\ldots,b_{r(j+1)-1}\rangle L_1 \\
&=b_{r(j)}\alpha_{\pi_j}\langle b_{r(j)+1},\ldots,b_{r(j+1)-1}\rangle L_1,
\end{align*}
while if $r(j)+1<\ell_j$, then by Lemmas~\ref{lem:ncnp} and~\ref{lem:sL1},
\begin{align*}
b_{r(j)}y_{s(r(j)+1)}b_{r(j)+1}\cdots&y_{s(r(j+1)-1)}b_{r(j+1)-1}y_{s(r(j+1))}= \\
&=b_{r(j)}\alpha_{\tau_j}\langle b_{r(j)+1},\ldots,b_{\ell_j-1}\rangle
b_{\ell_j}\alpha_{\rho_j}\langle b_{\ell_j+1},\ldots,b_{r(j+1)-1}\rangle L_1 \\
&=b_{r(j)}\alpha_{\pi_j}\langle b_{r(j)+1},\ldots,b_{r(j+1)-1}\rangle L_1.
\end{align*}
Therefore,
letting
\[
\bt_j=b_{r(j)}\alpha_{\pi_j}\langle b_{r(j)+1},\ldots,b_{r(j+1)-1}\rangle,
\]
we have
\[
y_{s(1)}b_1\cdots y_{s(n)}b_ny_{s(n+1)}
=V_{1,k}(\alpha_k)\bt_1L_1\bt_2L_1\cdots\bt_kL_1
=\alpha_k(\bt_1,\ldots,\bt_k).
\]
Using Lemma~\ref{lem:sL1} yields~\eqref{eq:Fo} and finishes the proof.
\end{proof}

\begin{defi}
Let $\beta\in\Mul[[B]]$ with $\beta_0$ invertible.
The {\em unsymmetrized T--transform} of $\beta$ is
\begin{equation}\label{eq:T}
\Tt_\beta=(\beta\circ(I\beta)\cinv)(1+I)^{-1}.
\end{equation}
The {\em unsymmetrized S--transform} of $\beta$ is 
$\St_\beta=\Tt_\beta^{-1}$.
\end{defi}

\begin{prop}\label{prop:Teqiv}
$\Tt_\beta$ is the unique element of $\Mul[[B]]$ satisfying
\begin{equation}\label{eq:Tchar}
(\Tt_\beta\circ(I\beta))(1+I\beta)=\beta.
\end{equation}
\end{prop}
\begin{proof}
The definition~\eqref{eq:T} is clearly equivalent to each of the following equations:
\begin{gather*}
\Tt_\beta\,(1+I)=\beta\circ(I\beta)\cinv \\
(\Tt_\beta\,(1+I))\circ(I\beta)=\beta \\
(\Tt_\beta\circ(I\beta))(1+I\beta)=\beta.
\end{gather*}
\end{proof}

\begin{prop}\label{prop:Scan}
Let $\alpha\in\Mul[[B]]$ with $\alpha_0$ invertible
and let
\[
X=V_1(\alpha)+W_1(\alpha)\in\Lc(\Qc)
\]
be the corresponding multiplicative canonical random variable.
Then $\Tt_{\Phit_X}=\alpha$.
\end{prop}
\begin{proof}
By Proposition~\ref{prop:Teqiv}, it will suffice to show
\begin{equation}\label{eq:alphaPhi}
\alpha\circ(I\Phit_X)+(\alpha\circ(I\Phit_X))I\Phit_X=\Phit_X.
\end{equation}
First, observe that the $0$th parts of both sides of~\eqref{eq:alphaPhi} are $\alpha_0$.
Let $n\in\Nats$ and $b_1,\ldots,b_n\in B$.
Then
\[
(\alpha\circ(I\Phit_X))(b_1,\ldots,b_n)
=\sum_{k=1}^n\sum_{\substack{p_1,\ldots,p_k\ge1 \\ p_1+\cdots+p_k=n}}
\begin{aligned}[t]
\alpha_k(&b_{q_1+1}\Phit_{X,p_1-1}(b_{q_1+2},\ldots,b_{q_1+p_1}), \\
&b_{q_2+1}\Phit_{X,p_2-1}(b_{q_2+2},\ldots,b_{q_2+p_2}), \\
&\ldots,b_{q_k+1}\Phit_{X,p_k-1}(b_{q_k+2},\ldots,b_{q_k+p_k})),
\end{aligned}
\]
where $q_1=0$ and $q_j=p_1+\cdots+p_{j-1}$ ($2\le j\le k$).
Using Proposition~\ref{prop:Xmultve}, we get
\begin{equation}\label{eq:pkterm}
\begin{split}
(\alpha\circ&(I\Phit_X))(b_1,\ldots,b_n)= \\
&=\sum_{k=1}^n\sum_{\substack{p_1,\ldots,p_k\ge1 \\ p_1+\cdots+p_k=n}}
\sum_{\substack{\pi_1\in\NCL(p_1) \\ \pi_2\in\NCL(p_2) \\ \stackdots \\ \pi_k\in\NCL(p_k)}}
\begin{aligned}[t]
\alpha_k(&b_{q_1+1}\alpha_{\pi_1}\langle b_{q_1+2},\ldots,b_{q_1+p_1}\rangle, \\
&b_{q_2+1}\alpha_{\pi_2}\langle b_{q_2+2},\ldots,b_{q_2+p_2}\rangle, \\
&\ldots,b_{q_k+1}\alpha_{\pi_k}\langle b_{q_k+2},\ldots,b_{q_k+p_k}\rangle).
\end{aligned} 
\end{split}
\end{equation}
By Lemma~\ref{lem:F}, the term of the sum~\eqref{eq:pkterm} is equal to $\alpha_\pi\langle b_1,\ldots,b_n\rangle$,
where
\[
\pi=(\Dco_{n+1})^{-1}((\pi_1,\ldots,\pi_k))\in\NCLo(n+1).
\]
Consequently, by part~(i) of Lemma~\ref{lem:pik},
\begin{equation}\label{eq:NCNP*}
(\alpha\circ(I\Phit_X))(b_1,\ldots,b_n)=\sum_{\pi\in\NCLo(n+1)}\alpha_\pi\langle b_1,\ldots,b_n\rangle.
\end{equation}
Similarly, we get
\begin{align}
((\alpha&\circ(I\Phit_X))I\Phit_X)(b_1,\ldots,b_n)= \notag \\
&=\alpha_0b_1\Phit_{X,n-1}(b_2,\ldots,b_n) \notag \\
&\quad+\sum_{m=1}^{n-1}(\alpha\circ(I\Phit_X))_m(b_1,\ldots,b_m)b_{m+1}\Phit_{X,n-m-1}(b_{m+2},\ldots,b_n) \notag \\[1ex]
&=\sum_{\sigma\in\NCL(n)}\alpha_0b_1\alpha_\sigma\langle b_2,\ldots,b_n\rangle \label{eq:sigsum} \\
&\quad+\sum_{m=1}^{n-1}
\sum_{k=1}^m\sum_{\substack{p_1,\ldots,p_k\ge1 \\ p_1+\cdots+p_k=m}}
\sum_{\substack{\pi_1\in\NCL(p_1) \\ \pi_2\in\NCL(p_2) \\ \stackdots \\ \pi_k\in\NCL(p_k)}}
\sum_{\sigma\in\NCL(n-m)} \label{eq:pisigsum} \\
&\qquad\qquad\qquad\qquad\qquad\qquad\qquad\qquad
\begin{aligned}[t]
\alpha_k(&b_{q_1+1}\alpha_{\pi_1}\langle b_{q_1+2},\ldots,b_{q_1+p_1}\rangle, \\
&b_{q_2+1}\alpha_{\pi_2}\langle b_{q_2+2},\ldots,b_{q_2+p_2}\rangle, \\
&\ldots,b_{q_k+1}\alpha_{\pi_k}\langle b_{q_k+2},\ldots,b_{q_k+p_k}\rangle) \\
&\qquad\cdot b_{m+1}\alpha_\sigma\langle b_{m+2},\ldots,b_n\rangle.
\end{aligned} \label{eq:pisigterm}
\end{align}
By Lemma~\ref{lem:F}, the term~\eqref{eq:pisigterm} is equal to $\alpha_\pi\langle b_1,\ldots,b_n\rangle$,
where
\[
\pi=(\Dc_{n+1})^{-1}((\pi_1,\ldots,\pi_k,\sigma))\in\NCL(n+1)\backslash\NCLo(n+1).
\]
Also, the term 
\[
\alpha_0b_1\alpha_\sigma\langle b_2,\ldots,b_n\rangle
\]
of the sum~\eqref{eq:sigsum}
is equal to $\alpha_\pi\langle b_1,\ldots,b_n\rangle$, where $\pi=(\Dc_{n+1})^{-1}(\sigma)$.
Therefore, by part~(ii) of Lemma~\ref{lem:pik},
\begin{equation}\label{eq:NCNP*c}
((\alpha\circ(I\Phit_X))I\Phit_X)(b_1,\ldots,b_n)
=\sum_{\pi\in\NCL(n+1)\backslash\NCLo(n+1)}\alpha_\pi\langle b_1,\ldots,b_n\rangle.
\end{equation}
Adding~\eqref{eq:NCNP*} and~\eqref{eq:NCNP*c} and
then using Proposition~\ref{prop:Xmultve}
finishes the proof of~\eqref{eq:alphaPhi}.
\end{proof}

In order to handle the product of two random variables, we now consider multivariable
generalizations some constructions
and results considered previously in this section.

\begin{defi}\label{def:alphai<>}
For every $i\in I$, let $\alpha^{(i)}\in\Mul[[B]]$ with $\alpha^{(i)}_0$ invertible.
Let $n\in\Nats_0$, let $\iota:\{1,\ldots,n+1\}\to I$ be a map and let $\pi\in\NCL(n+1)$.
We will define an element
\begin{equation}\label{eq:alphai<>}
\alpha^\iota_\pi\langle b_1,\ldots,b_n\rangle
\end{equation}
of $B$
for $b_1,\ldots,b_n\in B$.
If $\iota$ is not constant on all blocks of $\pi$, then we set
\[
\alpha^\iota_\pi\langle b_1,\ldots,b_n\rangle=0.
\]
Suppose $\iota$ is constant on all blocks of $\pi$.
Let $J=\{m,\ldots,m+\ell-1\}$ be the right--most interval block of $\pi$.
If $m=1$, which implies $\pi=1_{n+1}$, then we set
\begin{equation}\label{eq:alphai<>1}
\alpha^\iota_\pi\langle b_1,\ldots,b_n\rangle=
\alpha^{(\iota(1))}_n(b_1\alpha^{(\iota(1))}_0,b_2\alpha^{(\iota(1))}_0,\ldots,b_n\alpha^{(\iota(1))}_0).
\end{equation}
Suppose $m>1$.
Suppose $m$ is not contained in any other element of $\pi$.
Let $\iota':\{1,\ldots,n-\ell+1\}\to I$ be
\[
\iota'(j)=
\begin{cases}
\iota(j)&1\le j\le m-1 \\
\iota(j+\ell)&m\le j\le n-\ell+1
\end{cases}
\]
and let $\pi'=\pi\contr_{J^c}$ be the renumbered restriction of $\pi$ to the complement of $J$.
Then $\pi'\in\NCL(n-\ell+1)$ and $\iota'$ is constant on all blocks of $\pi'$.
If $J$ is all the way to the right, i.e.\ if $m+\ell-1=n+1$, then we set
\begin{equation}\label{eq:alphai<>2}
\alpha^\iota_\pi\langle b_1,\ldots,b_n\rangle=
\alpha^{\iota'}_{\pi'}\langle b_1,\ldots,b_{m-2}\rangle b_{m-1}
\alpha^{(\iota(m))}_{\ell-1}(b_m\alpha^{(\iota(m))}_0,\ldots,b_n\alpha^{(\iota(m))}_0),
\end{equation}
while if $m+\ell-1\le n$, then we set
\begin{equation}\label{eq:alphai<>3}
\alpha^\iota_\pi\langle b_1,\ldots,b_n\rangle
=\begin{aligned}[t]
\alpha^{\iota'}_{\pi'}\langle&b_1,\ldots,b_{m-2}, \\
&b_{m-1}\alpha^{(\iota(m))}_{\ell-1}(b_m\alpha^{(\iota(m))}_0,\ldots,b_{m+\ell-2}\alpha^{(\iota(m))}_0)b_{m+\ell-1}, \\
&b_{m+\ell},\ldots,b_n\rangle.
\end{aligned}
\end{equation}
Now suppose $m$ is contained in another block of $\pi$.  Then we must have $\ell\ge2$.
Let $\iota'':\{1,\ldots,n-\ell+2\}\to I$ be
\[
\iota''(j)=
\begin{cases}
\iota(j)&1\le j\le m \\
\iota(j+\ell-1)&m+1\le j\le n-\ell+2
\end{cases}
\]
and let $\pi''=(\pi\backslash\{J\})\contr_{J^c\cup\{m\}}$
be the renumbered restriction of $\pi\backslash\{J\}$ to $\{1,\ldots,m\}\cup\{m+\ell,\ldots,n+1\}$.
Then $\pi''\in\NCL(n-\ell+2)$, $\iota''$ is constant on all blocks of $\pi''$ and we set
\begin{equation}\label{eq:alphai<>4}
\alpha^\iota_\pi\langle b_1,\ldots,b_n\rangle
=\begin{aligned}[t]
\alpha^{\iota''}_{\pi''}\langle&b_1,\ldots,b_{m-2}, \\
&b_{m-1}\alpha^{(\iota(m))}_{\ell-1}(b_m\alpha^{(\iota(m))}_0,\ldots,b_{m+\ell-2}\alpha^{(\iota(m))}_0)
(\alpha^{(\iota(m))}_0)^{-1}, \\
&b_{m+\ell-1},b_{m+\ell},\ldots,b_n\rangle.
\end{aligned}
\end{equation}
Equations~\eqref{eq:alphai<>1}, \eqref{eq:alphai<>2}, \eqref{eq:alphai<>3} and~\eqref{eq:alphai<>4}
recursively define the quantity~\eqref{eq:alphai<>}
when $\iota$ is constant on all blocks of $\pi$.
\end{defi}

Here is the multivariable generalization of Lemma~\ref{lem:ncnp}.
\begin{lemma}\label{lem:incnp}
For every $i\in I$ take $\alpha^{(i)}\in\Mul[[B]]$ with $\alpha^{(i)}_0$ invertible
and for $s\in\Nats_0\cup\Nats_0^*$ and $i\in I$ let
\begin{equation}\label{eq:yis}
y_s^i=
\begin{cases}
V_{i,k}(\alpha^{(i)}_k),&s=k\in\Nats_0 \\
W_{i,k}(\alpha^{(i)}_k),&s=k^*\in\Nats_0^*
\end{cases}
\end{equation}
in $\Lc(\Qc)$.
Let $n\in\Nats_0$, let $\iota:\{1,\ldots,n+1\}\to I$ and $s:\{1,\ldots,n+1\}\to\Nats_0\cup\Nats_0^*$
be maps and let $b_1,\ldots,b_n\in B$.
If there is $\pi\in\NCL(n+1)$ such that $S_\pi=s$, where $S_\pi$ is as in Definition~\ref{def:Spi}, then
\[ 
y^{\iota(1)}_{s(1)}b_1y^{\iota(2)}_{s(2)}b_2\cdots y^{\iota(n)}_{s(n)}b_ny^{\iota(n+1)}_{s(n+1)}
=\alpha^\iota_\pi\langle b_1,\ldots,b_n\rangle.
\] 
If there is no $\pi\in\NCL(n+1)$ such that $S_\pi=s$, then
\begin{equation}\label{eq:iEc0ys}
\Ec(y^{\iota(1)}_{s(1)}b_1y^{\iota(2)}_{s(2)}b_2\cdots y^{\iota(n)}_{s(n)}b_ny^{\iota(n+1)}_{s(n+1)})=0.
\end{equation}
\end{lemma}
\begin{proof}
A key observation, which takes the place of equations~\eqref{eq:ysm1} and~\eqref{eq:ysm2},
is that if $s=S_\pi$,
if the interval $\{m,\ldots,m+k\}$ is a block of $\pi$ and if $\iota$ is not constant on this interval,
then
\[
y^{\iota(m)}_{s(m)}b_m\cdots y^{\iota(m+k-1)}_{s(m+k-1)}b_{m+k-1}y^{\iota(m+k)}_{s(m+k)}=0,
\]
while if it is constant in this interval, then
\begin{multline*}
y^{\iota(m)}_{s(m)}b_m\cdots y^{\iota(m+k-1)}_{s(m+k-1)}b_{m+k-1}y^{\iota(m+k)}_{s(m+k)}= \\[1ex]
=
\begin{cases}
\alpha^{(\iota(m))}_k(b_m\alpha^{(\iota(m))}_0,\ldots,b_{m+k-1}\alpha^{(\iota(m))}_0),
&m\text{ singly covered by }\pi \\
\alpha^{(\iota(m))}_k(b_m\alpha^{(\iota(m))}_0,\ldots,b_{m+k-1}\alpha^{(\iota(m))}_0)L_i,
&m\text{ doubly covered by }\pi.
\end{cases}
\end{multline*}
Using these additional facts, the proof proceeds similarly to the proof of Lemma~\ref{lem:ncnp}.
\end{proof}

Here is the multivariate analogue of Lemma~\ref{lem:insertNCNP}, which is proved similarly.
\begin{lemma}\label{lem:insertiNCNP}
Let $n\in\Nats$,
let $\pi\in\NCL(n+1)$ and suppose 
\[
J=\{m,\ldots,m+\ell-1\}\subseteq\{1,\ldots,n+1\}
\]
is a proper subinterval that splits $\pi$.
Let $\pi''$ be the renumbered restriction $\pi$ to $J$, and let $\pi'$
be the renumbered restriction of $\pi$ to the complement of $J$.
Let 
\[
\iota:\{1,\ldots,n+1\}\to I,
\]
let $\iota'':\{1,\ldots,\ell\}\to I$ be $\iota''(j)=\iota(j+m-1)$ and let
$\iota':\{1,\ldots,n-\ell+1\}\to I$ be
\[
\iota'(j)=
\begin{cases}
\iota(j),&1\le j\le m-1 \\
\iota(j+\ell),&m\le j\le n-\ell+1.
\end{cases}
\]
Then for all $\alpha\in\Mul[[B]]$ and all $b_1,\ldots,b_n\in B$, we have
\begin{multline*}
\alpha^\iota_\pi\langle b_1,\ldots,b_n\rangle= \\
=\begin{cases}
\alpha^{\iota''}_{\pi''}\langle b_1,\ldots,b_{\ell-1}\rangle b_\ell
\alpha^{\iota'}_{\pi'}\langle b_{\ell+1},\ldots,b_n\rangle,&m=1 \\
\begin{aligned}[b]
\alpha^{\iota'}_{\pi'}\langle b_1,&\ldots,b_{m-2}, \\
&b_{m-1}\alpha^{\iota''}_{\pi''}\langle b_m,\ldots,b_{m+\ell-2}\rangle b_{m+\ell-1},
b_{m+\ell},\ldots,b_n\rangle,\end{aligned}&2\le m\le n-\ell+1 \\
\alpha^{\iota'}_{\pi'}\langle b_1,\ldots,b_{m-2}\rangle b_{m-1}
\alpha^{\iota''}_{\pi''}\langle b_m,\ldots,b_n\rangle,&m=n-\ell+2.
\end{cases}
\end{multline*}
\end{lemma}

Now we have the multivariate analogue of Lemma~\ref{lem:F}, which is proved similarly.
\begin{lemma}\label{lem:iF}
Let $n\in\Nats$, let $\pi\in\NCL(n+1)$, let $\iota:\{1,\ldots,n+1\}\to I$ and assume
that $\iota$ is constant on all blocks of $\pi$.
If $\pi\in\NCLo(n+1)$, then write $\Dco_n(\pi)=(\pi_1,\ldots,\pi_k)$, with $\pi_j\in\NCL(p_j)$.
If $\pi\notin\NCLo(n+1)$, then write $\Dc_n(\pi)=(\pi_1,\ldots,\pi_k,\sigma)$, with $\pi_j\in\NCL(p_j)$.
Let $r(1)=1$ and for $2\le j\le k+1$, let $r(j)=1+p_1+\cdots+p_{j-1}$.
For $1\le j\le k$ (and for $j=k+1$ when $\pi\notin\NCLo(n+1)$),
let $\iota_j:\{1,\ldots,p_j\}\to I$ be $\iota_j(x)=\iota(r(j)+x)$.
For every $i\in I$, take $\alpha^{(i)}\in\Mul[[B]]$ with $\alpha^{(i)}_0$ invertible
and let $b_1,\ldots,b_n\in B$.
\renewcommand{\labelenumi}{(\roman{enumi})}
\begin{enumerate}
\item
If $\pi\in\NCLo(n+1)$, then 
\[
\alpha^\iota_\pi\langle b_1,\ldots,b_n\rangle=
\begin{aligned}[t]
\alpha^{(\iota(1))}_k\big(&b_{r(1)}\alpha^{\iota_1}_{\pi_1}\langle b_{r(1)+1},\ldots,b_{r(2)-1}\rangle, \\
&b_{r(2)}\alpha^{\iota_2}_{\pi_2}\langle b_{r(2)+1},\ldots,b_{r(3)-1}\rangle, \\
&\ldots,b_{r(k)}\alpha^{\iota_k}_{\pi_k}\langle b_{r(k)+1},\ldots,b_{r(k+1)-1}\rangle\big).
\end{aligned}
\]
\item
If $\pi\notin\NCLo(n+1)$, then
\[
\alpha^\iota_\pi\langle b_1,\ldots,b_n\rangle=
\begin{aligned}[t]
\alpha^{(\iota(1))}_k\big(&b_{r(1)}\alpha^{\iota_1}_{\pi_1}\langle b_{r(1)+1},\ldots,b_{r(2)-1}\rangle, \\
&b_{r(2)}\alpha^{\iota_2}_{\pi_2}\langle b_{r(2)+1},\ldots,b_{r(3)-1}\rangle, \\
&\ldots,b_{r(k)}\alpha^{\iota_k}_{\pi_k}\langle b_{r(k)+1},\ldots,b_{r(k+1)-1}\rangle\big) \\[0.3ex]
&\qquad\qquad\cdot b_{r(k+1)}\alpha^{\iota_{k+1}}_\sigma\langle b_{r(k+1)+1},\ldots,b_n\rangle.
\end{aligned}
\]
\end{enumerate}
\end{lemma}

\begin{thm}\label{thm:Ttrans}
Let $x$ and $y$ be free random variables in any $B$--valued noncommutative probability space
$(A,E)$, and assume $E(x)$ and $E(y)$ are invertible elements of $B$.
Then
\begin{equation}\label{eq:Txy}
\Tt_{\Phit_{xy}}=(\Tt_{\Phit_x}\circ(\Tt_{\Phit_y}\,I\,\Tt_{\Phit_y}^{-1}))\,\Tt_{\Phit_y}.
\end{equation}
\end{thm}
\begin{proof}
Let $\alpha^{(1)}=\Tt_{\Phit_x}$
and $\alpha^{(2)}=\Tt_{\Phit_y}$.
Consider multiplicative canonical random variables
\[
X_i=V_i(\alpha^{(i)})+W_i(\alpha^{(i)})\in\Lc(\Qc),\qquad(i=1,2)
\]
as constructed in Proposition~\ref{prop:can}, where we take index set $I=\{1,2\}$
in the construction of $\Qc$.
By Proposition~\ref{prop:Scan}, $\Tt_{\Phit_{X_i}}=\alpha^{(i)}$, and by Proposition~\ref{prop:free},
$X_1$ and $X_2$ are free.
Let
\[
\gamma=(\alpha^{(1)}\circ(\alpha^{(2)}\,I\,(\alpha^{(2)})^{-1}))\,\alpha^{(2)}.
\]
It will suffice to show $\Tt_{\Phit_{X_1X_2}}=\gamma$.
By Proposition~\ref{prop:Teqiv}, it will suffice to show
\begin{equation}\label{eq:s0}
(\gamma\circ(I\Phit_{X_1X_2}))(1+I\Phit_{X_1X_2})=\Phit_{X_1X_2}.
\end{equation}
Let $f_n,g_n:\{1,\ldots,n\}\to\{1,2\}$ be
\[
f_n(j)=
\begin{cases}
1,&j\text{ odd} \\
2,&j\text{ even,}
\end{cases}
\qquad\qquad
g_n(j)=
\begin{cases}
2,&j\text{ odd} \\
1,&j\text{ even.}
\end{cases}
\]
The domain of these functions will be obvious from the context, and we will drop the subscripts
and refer to these functions as $f$ and $g$, respectively.
Using the notation~\eqref{eq:yis} and Lemma~\ref{lem:incnp}, for $n\ge1$ and $b_1,\ldots,b_n\in B$,
we get
\begin{align}
\Phit_{X_1X_2,n}(b_1,\ldots,b_n)&=\Ec(X_1X_2b_1X_1X_2b_2\cdots X_1X_2b_nX_1X_2) \notag \\
&=\sum_{\substack{s_1(1),\ldots,s_1(n+1)\in\Nats_0\cup\Nats_0^* \\
 s_2(1),\ldots,s_2(n+1)\in\Nats_0\cup\Nats_0^*}}
\begin{aligned}[t]
\Ec(&y^1_{s_1(1)}y^2_{s_2(1)}b_1y^1_{s_1(2)}y^2_{s_2(2)}b_2\cdots \\
 &\cdots y^1_{s_1(n)}y^2_{s_2(n)}b_ny^1_{s_1(n+1)}y^2_{s_2(n+1)})
\end{aligned} \notag \\[1ex]
&=\sum_{\pi\in\NCL(2n+2)}\alpha^f_\pi\langle1,b_1,1,b_2,\ldots,1,b_n,1\rangle. \label{eq:s1}
\end{align}
Moreover,
\[
\Phit_{X_1X_2,0}=\Ec(X_1X_2)=\alpha^{(1)}_0\alpha^{(2)}_0=\alpha^f_{0_2}\langle1\rangle,
\]
and~\eqref{eq:s1} holds also for $n=0$.
Let
\[
\beta=\alpha^{(2)}\circ(I\Phit_{X_1X_2}).
\]
Then $\beta_0=\alpha^{(2)}_0$, while for $n\ge1$,
\begin{equation}\label{eq:s2}
\beta_n(b_1,\ldots,b_n)=\sum_{k=1}^n\sum_{\substack{p_1,\ldots,p_k\ge1 \\ p_1+\cdots+p_k=n}}
\begin{aligned}[t]
\alpha^{(2)}_k(&b_{q_1+1}\Phit_{X_1X_2,p_1-1}(b_{q_1+2},\ldots,b_{q_1+p_1}), \\
&b_{q_2+1}\Phit_{X_1X_2,p_2-1}(b_{q_2+2},\ldots,b_{q_2+p_2}), \\
&\cdots,b_{q_k+1}\Phit_{X_1X_2,p_k-1}(b_{q_k+2},\ldots,b_{q_k+p_k})),
\end{aligned}
\end{equation}
where as usual, we take $q_j=p_1+\cdots+p_{j-1}$.
Using~\eqref{eq:s1}, we have
\begin{align}
\beta_n(b_1,\ldots,b_n)=&\sum_{k=1}^n\sum_{\substack{p_1,\ldots,p_k\ge1 \\ p_1+\cdots+p_k=n}}
\sum_{\substack{\pi_1\in\NCL(2p_1) \\ \pi_2\in\NCL(2p_2) \\ \stackdots \\ \pi_k\in\NCL(2p_k)}} \notag \\
&\qquad\begin{aligned}[t]
\alpha_k^{(2)}(&b_{q_1+1}\alpha^f_{\pi_1}\langle1,b_{q_1+2},1,b_{q_1+3},\ldots,1,b_{q_1+p_1},1\rangle, \\
&b_{q_2+1}\alpha^f_{\pi_2}\langle1,b_{q_2+2},1,b_{q_2+3},\ldots,1,b_{q_2+p_2},1\rangle, \\
&\ldots,b_{q_k+1}\alpha^f_{\pi_k}\langle1,b_{q_k+2},1,b_{q_k+3},\ldots,1,b_{q_k+p_k},1\rangle)
\end{aligned}
\label{eq:s3}
\end{align}
By Lemma~\ref{lem:iF}, the term~\eqref{eq:s3} is equal to
$\alpha_\pi^g\langle b_1,1,b_2,1,\ldots,b_n,1\rangle$,
where
\[
\pi=(\Dco_{2n+1})^{-1}((\pi_1,\ldots,\pi_k)).
\]
By Lemma~\ref{lem:pik} and the requirement that the blocks of the noncrossing linked partitions
must consist of exclusively odd or exclusively even numbers in order for the terms of our series
to be nonzero, for $n\ge1$ we have
\begin{equation}\label{eq:s4}
\beta_n(b_1,\ldots,b_n)=\sum_{\pi\in\NCLo(2n+1)}
\alpha_\pi^g\langle b_1,1,b_2,1,\ldots,b_n,1\rangle.
\end{equation}
Moreover, $\beta_0=\alpha^{(2)}_0=\alpha^g_{1_1}\langle\;\rangle$,
so~\eqref{eq:s4} holds also when $n=0$.

We have $(\beta^{-1})_0=(\beta_0)^{-1}$, and we will write simply $\beta_0^{-1}$ for this element of $B$.
Using Proposition~\ref{prop:1+a}, we have
\[
\begin{split}
\beta^{-1}&=(\beta_0+(\beta-\beta_0))^{-1}=(1+\beta_0^{-1}(\beta-\beta_0))^{-1}\beta_0^{-1} \\
&=\beta_0^{-1}+\sum_{k=1}^\infty(-1)^k(\beta_0^{-1}(\beta-\beta_0))^k\beta_0^{-1}.
\end{split}
\]
Therefore, for $n\ge1$ and $b_1,\ldots,b_n\in B$,
\begin{equation}\label{eq:betainv}
\begin{split}
(\beta&^{-1})_n(b_1,\ldots,b_n)= \\
&=\sum_{k=1}^n(-1)^k\sum_{\substack{p_1,\ldots,p_k\ge1 \\ p_1+\cdots+p_k=n}}
\begin{aligned}[t]
\beta_0^{-1}\beta_{p_1}(b_{q_1+1},\ldots,b_{q_1+p_1})
&\beta_0^{-1}\beta_{p_2}(b_{q_2+1},\ldots,b_{q_2+p_2}) \\
\cdots&\beta_0^{-1}\beta_{p_k}(b_{q_k+1},\ldots,b_{q_k+p_k})\beta_0^{-1}.
\end{aligned}
\end{split}
\end{equation}
We have
\begin{equation}\label{eq:s5}
(\alpha^{(2)}I(\alpha^{(2)})^{-1})\circ(I\Phit_{X_1X_2})=\beta I\Phit_{X_1X_2}\beta^{-1}.
\end{equation}
Using~\eqref{eq:s1}, \eqref{eq:s4} and~\eqref{eq:betainv}, for $n\ge1$ and $b_1,\ldots,b_n\in B$ we get
\begin{align}
(\Phit_{X_1X_2}&\beta^{-1})_n(b_1,\ldots,b_n)= \notag \\
&=\Phit_{X_1X_2,n}(b_1,\ldots,b_n)\beta_0^{-1} \notag \\
&\quad+\sum_{m=1}^n\Phit_{X_1X_2,n-m}(b_1,\ldots,b_{n-m})
 \sum_{k=1}^m(-1)^k \notag \\
&\quad\qquad\sum_{\substack{p_1,\ldots,p_k\ge1 \\ p_1+\cdots+p_k=m}}
 \begin{aligned}[t]
 &\beta_0^{-1}\beta_{p_1}(b_{n-m+q_1+1},\ldots,b_{n-m+q_1+p_1}) \\
 &\beta_0^{-1}\beta_{p_2}(b_{n-m+q_2+1},\ldots,b_{n-m+q_2+p_2}) \\
 &\cdots\beta_0^{-1}\beta_{p_k}(b_{n-m+q_k+1},\ldots,b_{n-m+q_k+p_k})\beta_0^{-1}
\end{aligned} \notag \displaybreak[2] \\[1ex]
&=\sum_{\rho\in\NCL(2n+2)}\alpha^f_\rho\langle1,b_1,\ldots,1,b_n,1\rangle\beta_0^{-1} \notag \\
&\quad+\sum_{m=1}^n\sum_{k=1}^m(-1)^k\sum_{\substack{p_1,\ldots,p_k\ge1 \\ p_1+\cdots+p_k=m}}
\quad\sum_{\rho\in\NCL(2n-2m+2)}
\quad\sum_{\substack{\pi_1\in\NCLo(2p_1+1) \\ \pi_2\in\NCLo(2p_2+1) \\ \stackdots \\ \pi_k\in\NCLo(2p_k+1)}} \notag \\
&\quad\qquad\begin{aligned}[t]
\alpha&^f_\rho\langle1,b_1,\ldots,1,b_{n-m},1\rangle
\beta_0^{-1}\alpha^g_{\pi_1}\langle b_{n-m+q_1+1},1,\ldots,b_{n-m+q_1+p_1},1\rangle \\
&\beta_0^{-1}\alpha^g_{\pi_2}\langle b_{n-m+q_2+1},1,\ldots,b_{n-m+q_2+p_2},1\rangle \\
&\cdots\beta_0^{-1}\alpha^g_{\pi_k}\langle b_{n-m+q_k+1},1,\ldots,b_{n-m+q_k+p_k},1\rangle\beta_0^{-1}
\end{aligned} \notag \displaybreak[3] \\[1ex]
&=\sum_{k=0}^n(-1)^k\sum_{\substack{p_0\ge0 \\ p_1,\ldots,p_k\ge1 \\ p_0+p_1+\cdots+p_k=n}}
\quad\sum_{\rho\in\NCL(2p_0+2)}
\quad\sum_{\substack{\pi_1\in\NCLo(2p_1+1) \\ \pi_2\in\NCLo(2p_2+1) \\ \stackdots \\ \pi_k\in\NCLo(2p_k+1)}}
  \label{eq:s7sum} \\
&\quad\qquad\begin{aligned}[t]
\alpha&^f_\rho\langle1,b_1,\ldots,1,b_{p_0},1\rangle
\beta_0^{-1}\alpha^g_{\pi_1}\langle b_{\qt_1+1},1,\ldots,b_{\qt_1+p_1},1\rangle \\
&\beta_0^{-1}\alpha^g_{\pi_2}\langle b_{\qt_2+1},1,\ldots,b_{\qt_2+p_2},1\rangle \\
&\cdots\beta_0^{-1}\alpha^g_{\pi_k}\langle b_{\qt_k+1},1,\ldots,b_{\qt_k+p_k},1\rangle\beta_0^{-1}, \label{eq:s7term}
\end{aligned}
\end{align}
where $\qt_j=p_0+p_1+\cdots+p_{j-1}$.
We claim that the nonzero terms of the sum~\eqref{eq:s7sum} with $\{2p_0+2\}\in\rho$ and $k\ge1$ cancel exactly
the nonzero terms with $\{2p_0+2\}\notin\rho$.
Indeed, if $k\ge1$ and $\{2p_0+2\}\in\rho$, then
\[
\alpha^f_\rho\langle1,b_1,\ldots,1,b_{p_0},1\rangle
=\alpha^f_{\rho'}\langle1,b_1,\ldots,1,b_{p_0}\rangle\alpha^{(2)}_0,
\]
where 
\begin{equation}\label{eq:rho'}
\rho'=\rho\contr_{[1,2p_0+1]}=\rho\backslash\{\{2p_0+2\}\}\in\NCL(2p_0+1).
\end{equation}
But then using Lemma~\ref{lem:insertiNCNP} we get
\begin{align*}
\alpha^f_\rho\langle1,b_1,\ldots,1,b_{p_0},&1\rangle
\beta_0^{-1}\alpha^g_{\pi_1}\langle b_{\qt_1+1},1,\ldots,b_{\qt_1+p_1},1\rangle= \\
&=\alpha^f_{\rho'}\langle1,b_1,\ldots,1,b_{p_0}\rangle
\alpha^g_{\pi_1}\langle b_{\qt_1+1},1,\ldots,b_{\qt_1+p_1},1\rangle \\
&=\alpha^f_{\rho'\oplus\pi_1}\langle1,b_1,\ldots,1,b_{p_0},1,b_{\qt_1+1},1,\ldots,b_{\qt_1+p_1},1\rangle.
\end{align*}
On the other hand, if $p_0\ge0$, $\sigma\in\NCL(2p_0+2)$ and if $\{2p_0+2\}\notin\sigma$ and
the term~\eqref{eq:s7term} is nonzero, then $p_0\ge1$ and
each block of $\sigma$ consists of either all even or all odd numbers.
Letting $\sigmahat$ be the noncrossing partition generated by $\sigma$, it follows that the minimal element
of the block of $\sigmahat$ that contains $2p_0+2$ is $2\ell$ for some $\ell\in\{1,2,\ldots,p_0\}$,
and the interval $[1,2\ell-1]$ splits $\sigma$.
Therefore, we have $\sigma=\rho'\oplus\pi_0$, where
\[
\rho'=\sigma\contr_{[1,2\ell-1]},\qquad
\pi_0=\sigma\contr_{[2\ell,2p_0+2]}\in\NCLo(2p_0-2\ell+3)
\]
and
\[
\alpha^f_\sigma\langle1,b_1,\ldots,1,b_{p_0},1\rangle
=\alpha^f_\rho\langle1,b_1,\ldots,b_{\ell-1},1\rangle\beta_0^{-1}
\alpha^g_{\pi_0}\langle b_\ell,1,\ldots,b_{p_0},1\rangle,
\]
where $\rho=\rho'\cup\{\{2\ell\}\}$.
We have exhibited a bijective correspondence between
the set of nonzero terms of the sum~\eqref{eq:s7sum} having $\{2p_0+2\}\in\rho$ and $k\ge1$
and the set of nonzero terms with $\{2p_0+2\}\notin\rho$, and the corresponding terms appear in
the sum~\eqref{eq:s7sum} with opposite signs.
After cancelling, the only survivors are those terms with $k=0$, $p_0=n$, and $\{2n+2\}\in\rho$,
in which case the term~\eqref{eq:s7term} becomes
\[
\alpha^f_\rho\langle1,b_1,\ldots,1,b_n,1\rangle\beta_0^{-1}
=\alpha^f_{\rho'}\langle1,b_1,\ldots,1,b_n\rangle,
\]
where $\rho'$ is as in~\eqref{eq:rho'}.
We find
\begin{equation}\label{eq:PhiXX}
(\Phit_{X_1X_2}\beta^{-1})_n(b_1,\ldots,b_n)=
\sum_{\sigma\in\NCL(2n+1)}\alpha^f_\sigma\langle1,b_1,\ldots,1,b_n\rangle.
\end{equation}
Moreover, 
\[
(\Phit_{X_1X_2}\beta^{-1})_0=\alpha^{(1)}_0\alpha^{(2)}_0\beta_0^{-1}=\alpha^{(1)}_0
=\sum_{\sigma\in\NCL(1)}\alpha^f_\sigma\langle\;\rangle,
\]
so~\eqref{eq:PhiXX} holds also when $n=0$.
Combining~\eqref{eq:s4} with~\eqref{eq:PhiXX}, we have, for $n\ge1$ and $b_1,\ldots,b_n\in B$,
\begin{align}
(\beta I\Phit_{X_1X_2}\beta&^{-1})_n(b_1,\ldots,b_n)= \notag \\
&=\sum_{k=1}^n\sum_{\substack{\pi\in\NCLo(2k-1) \\ \sigma\in\NCL(2n-2k+1)}}
\alpha^g_\pi\langle b_1,1,\ldots,b_{k-1},1\rangle b_k
\alpha^f_\sigma\langle1,b_{k+1},\ldots,1,b_n\rangle \notag \\
&=\sum_{k=1}^n\sum_{\substack{\pi\in\NCLo(2k-1) \\ \sigma\in\NCL(2n-2k+1)}}
\alpha^g_{\pi\oplus\sigma}\langle b_1,1,\ldots,b_{n-1},1,b_n\rangle, \label{eq:s8}
\end{align}
where for the last equality we used Lemma~\ref{lem:insertiNCNP}.
However, every $\tau\in\NCL(2n)$, all of whose blocks consist of exclusively odd or
exclusively even numbers, can be uniquely written as $\tau=\pi\oplus\sigma$ for some $\pi$, $\sigma$
and $k$ as in the sum~\eqref{eq:s8}.
Therefore,
\begin{equation}\label{eq:bIPb}
(\beta I\Phit_{X_1X_2}\beta^{-1})_n(b_1,\ldots,b_n)=
\sum_{\tau\in\NCL(2n)}
\alpha^g_\tau\langle b_1,1,\ldots,b_{n-1},1,b_n\rangle.
\end{equation}
Let $\eps=\beta I\Phit_{X_1X_2}\beta^{-1}$.
Then
\begin{equation}\label{eq:s8.5}
\gamma\circ(I\Phit_{X_1X_2})=(\alpha^{(1)}\circ\eps)(\alpha^{(2)}\circ(I\Phit_{X_1X_2}))
=(\alpha^{(1)}\circ\eps)\beta.
\end{equation}
For $n\ge1$ and $b_1,\ldots,b_n\in B$, using~\eqref{eq:bIPb} we have
\begin{align}
(\alpha^{(1)}&\circ\eps)_n(b_1,\ldots,b_n)= \notag \\
&=\sum_{k=1}^n\sum_{\substack{p_1,\ldots,p_k\ge1 \\ p_1+\cdots+p_k=n}}
\begin{aligned}[t]
\alpha_k^{(1)}(&\eps_{p_1}(b_{q_1+1},\ldots,b_{q_1+p_1}), \\
&\eps_{p_2}(b_{q_2+1},\ldots,b_{q_2+p_2}), \\
&\ldots,\eps_{p_k}(b_{q_k+1},\ldots,b_{q_k+p_k}))
\end{aligned} \notag \displaybreak[2] \\
&=\sum_{k=1}^n\sum_{\substack{p_1,\ldots,p_k\ge1 \\ p_1+\cdots+p_k=n}}
\sum_{\substack{\pi_1\in\NCLo(2p_1) \\ \pi_2\in\NCLo(2p_2) \\ \stackdots \\ \pi_k\in\NCLo(2p_k)}}
\begin{aligned}[t]
\alpha_k^{(1)}(&\alpha^g_{\pi_1}\langle b_{q_1+1},\ldots,b_{q_1+p_1}\rangle, \\
&\alpha^g_{\pi_2}\langle b_{q_2+1},\ldots,b_{q_2+p_2}\rangle, \\
&\ldots,\alpha^g_{\pi_k}\langle b_{q_k+1},\ldots,b_{q_k+p_k}\rangle).
\end{aligned} \label{eq:s9}
\end{align}
By Lemma~\ref{lem:iF}, the term~\eqref{eq:s9} is equal to
$\alpha_\pi^f\langle1,b_1,1,b_2,1,\ldots,b_n\rangle$,
where
\[
\pi=(\Dco_{2n+1})^{-1}((\pi_1,\ldots,\pi_k))\in\NCLo(2n+1).
\]
By Lemma~\ref{lem:pik} and the requirement that the blocks of the noncrossing linked partitions
must consist of exclusively odd or exclusively even numbers in order for the terms of our series
to be nonzero, for $n\ge1$ we have
\begin{equation}\label{eq:s10}
(\alpha^{(1)}\circ\eps)_n(b_1,\ldots,b_n)
=\sum_{\pi\in\NCLo(2n_1)}\alpha_\pi^f\langle1,b_1,1,b_2,1,\ldots,b_n\rangle.
\end{equation}
Moreover,
\[
(\alpha^{(1)}\circ\eps)_0=\alpha^{(1)}_0=\alpha^f_{1_1}\langle\;\rangle,
\]
so~\eqref{eq:s10} holds also for $n=0$.
Using~\eqref{eq:s4} and~\eqref{eq:s8.5}, we have, for $n\ge1$,
\begin{align}
(\gamma\circ(I\Phit_{X_1X_2}&))_n(b_1,\ldots,b_n)= \notag \\
&=\sum_{m=0}^n\sum_{\substack{\pi\in\NCLo(2m+1) \\ \sigma\in\NCLo(2n-2m+1)}}
\alpha^f_\pi\langle1,b_1,\ldots,1,b_m\rangle\alpha^g_\sigma\langle b_{m+1},1,\ldots,b_n,1\rangle \notag \\
&=\sum_{m=0}^n\sum_{\substack{\pi\in\NCLo(2m+1) \\ \sigma\in\NCLo(2n-2m+1)}}
\alpha^f_{\pi\oplus\sigma}\langle1,b_1,\ldots,1,b_n,1\rangle, \label{eq:s11}
\end{align}
where we used Lemma~\ref{lem:insertiNCNP} for the last equality.
Moreover,
\[
(\gamma\circ(I\Phit_{X_1X_2}))_0=\gamma_0=\alpha^{(1)}_0\alpha^{(2)}_0
=\alpha^f_{1_1\oplus1_1}\langle1\rangle,
\]
so~\eqref{eq:s11} holds also for $n=0$.
Using~\eqref{eq:s1} and~\eqref{eq:s11}, we have for $n\ge1$,
\begin{align}
((\gamma&\circ(I\Phit_{X_1X_2}))(I\Phit_{X_1X_2}))_n(b_1,\ldots,b_n)= \notag \\
&=\sum_{\ell=1}^n\sum_{m=0}^{\ell-1}
\sum_{\substack{\pi\in\NCLo(2m+1) \\ \sigma\in\NCLo(2\ell-2m-1) \\ \rho\in\NCL(2n-2\ell+2)}}
\alpha^f_{\pi\oplus\sigma}\langle1,b_1,\ldots,1,b_{\ell-1},1\rangle
b_\ell\alpha^f_\rho\langle1,b_{\ell+1},\ldots,1,b_n,1\rangle. \notag \\
&=\sum_{\ell=1}^n\sum_{m=0}^{\ell-1}
\sum_{\substack{\pi\in\NCLo(2m+1) \\ \sigma\in\NCLo(2\ell-2m-1) \\ \rho\in\NCL(2n-2\ell+2)}}
\alpha^f_{\pi\oplus\sigma\oplus\rho}\langle1,b_1,\ldots,1,b_n,1\rangle, \label{eq:s12}
\end{align}
where again we used Lemma~\ref{lem:insertiNCNP} for the last equality.

Let us now see that~\eqref{eq:s0} holds.
For $n\ge1$,
\[
((\gamma\circ(I\Phit_{X_1X_2}))(1+I\Phit_{X_1X_2}))_n(b_1,\ldots,b_n)
\]
is equal to the sum of~\eqref{eq:s11} and~\eqref{eq:s12},
and we want to see that this equals
\[
\Phit_{X_1X_2,n}(b_1,\ldots,b_n).
\]
Using~\eqref{eq:s1}, it will suffice to see that every element $\tau$ of $\NCL(2n+2)$,
all of whose blocks consist of either all even or all odd numbers,
can be uniquely written as either $\tau=\pi\oplus\sigma$
for $\pi$ and $\sigma$ as in the summaton~\eqref{eq:s11}
or $\tau=\pi\oplus\sigma\oplus\rho$ for $\pi$, $\sigma$ and $\rho$
as in the summation~\eqref{eq:s12}.
But this is easy to see.
Indeed, given $\tau$, let $\tauhat$ be the noncrossing partition generated by $\tau$.
Let $E\in\tauhat$ be such that $1\in E$.
Then $\max(E)=2m+1$ for some $m\in\{0,\ldots,n\}$.
Letting $\pi=\tau\contr_{[1,2m+1]}$, we have $\pi\in\NCLo(2m+1)$.
Let $F\in\tauhat$ be such that $2m+2\in F$.
Then $\max(F)=2\ell$ for some $\ell\in\{m+1,\ldots,n+1\}$.
Letting $\sigma=\tau\contr_{[2m+2,2\ell]}$, we have $\sigma\in\NCLo(2\ell-2m-1)$.
If $\ell=n+1$, then $\tau=\pi\oplus\sigma$.
If $\ell\le n$, then letting $\rho=\tau\contr_{[2\ell+1,2n+2]}$, we have $\rho\in\NCL(2n-2\ell+2)$
and $\tau=\pi\oplus\sigma\oplus\rho$.
In either case, the indicated decomposition for $\tau$ is clearly the only possible decomposition of
either sort.
Of course, we also have
\[
((\gamma\circ(I\Phit_{X_1X_2}))(1+I\Phit_{X_1X_2}))_0=\gamma_0
=\alpha^{(1)}_0\alpha^{(2)}_0=\Phit_{X_1X_2,0}.
\]
This finishes the proof of~\eqref{eq:s0}, and of the theorem.
\end{proof}

\begin{remark}
The definition and properties of the usual T--transform,
given as equations\eqref{eq:Tdefeq} and ~\eqref{eq:Ttmult}
in the introduction,
follow from Proposition~\ref{prop:Teqiv} and Theorem~\ref{thm:Ttrans}
by taking the quotient from $\Mul[[B]]$ into $\SymMul[[B]]$ and 
using Proposition~\ref{prop:mfsanalytic}.
\end{remark}

\section{The scalar case and an application to the cardinality of $\NCL$}
\label{sec:scalar}

In the case $B=\Cpx$, considered as an algebra over $\Cpx$, the unsymmetrized
T--transform equals the usual T--transform.
The formula from the last section for the moments of a random variable in terms of the T--transform
of its distribution 
is especially nice and may be of interest.
We record it here, and obtain as an application the cardinality of $\NCL(n)$.

Since every multilinear function $\beta_n\in\Lc_n(\Cpx)$ is 
of the form
\[
\beta_n(z_1,\ldots,z_n)=\lambda z_1\cdots z_n
\]
for some $\lambda\in\Cpx$, we write
\[
\Mul[[\Cpx]]=\{\alpha=(\alpha_0,\alpha_1,\ldots)\mid\alpha_n\in\Cpx\}.
\]
We identify $\Mul[[\Cpx]$ with the algebra
of formal power series $\Cpx[[z]]$ in one variable with complex coefficients, and thereby identify
the unsymmetrized T--transform $\Tt_{\Phit_a}$
of the distribution series of a random variable $a$
in a $\Cpx$--valued noncommutative probability space with the
usual T--transfrom $T_{\Phit_a}$.

\begin{prop}\label{prop:CpxT}
Let $(A,\phi)$ be a $\Cpx$--valued noncommutative probability space and let $a\in A$
with $\phi(a)\ne0$.
Let
$T_{\Phit_a}=\alpha=(\alpha_0,\alpha_1,\ldots)\in\Mul[[\Cpx]]$
be the T--transform of the moment generating function of $a$.
Then for every positive integer $n$,
\[
\phi(a^n)=\sum_{\pi\in\NCL(n)}\alpha_0^{n-|\pi|}\prod_{B\in\pi}\alpha_{|B|-1}.
\]
\end{prop}
\begin{proof}
Combining Proposition~\ref{prop:Xmultve} and Proposition~\ref{prop:Scan}, we get
\[
\phi(a^n)=\sum_{\pi\in\NCL(n)}\alpha_\pi\langle1,\ldots,1\rangle.
\]
However, an examination of Definition~\ref{def:alpha<>} shows
\[
\alpha_\pi\langle1,\ldots,1\rangle=\alpha_0^{n-|\pi|}\prod_{B\in\pi}\alpha_{|B|-1}.
\]
\end{proof}

\begin{example}
We get
\begin{align*}
\phi(a)&=\alpha_0 \\
\phi(a^2)&=\alpha_0^2+\alpha_0\alpha_1 \\
\phi(a^3)&=\alpha_0^3+3\alpha_0^2\alpha_1+\alpha_0\alpha_1^2+\alpha_0^2\alpha_2 \\
\phi(a^4)&=\alpha_0^4+6\alpha_0^3\alpha_1+6\alpha_0^2\alpha_1^2+4\alpha_0^3\alpha_2
+\alpha_0\alpha_1^3+3\alpha_0^2\alpha_1\alpha_2+\alpha_0^3\alpha_3
\end{align*}
and writing $m_k=\phi(a^k)$ we obtain
\begin{align*}
\alpha_0&=m_1 \\
\alpha_1&=m_1^{-1}m_2-m_1 \\
\alpha_2&=m_1^{-2}m_3-m_1^{-3}m_2^2-m_1^{-1}m_2+m_1.
\end{align*}
\end{example}

The large Schr\"oder numbers that appear below 
are well known 
(see~\S6.2 and Exercises~6.39 and~6.40 of~\cite{St})
and have several combinatorial
characterizations.

\begin{thm}\label{thm:cardNCL}
For every $n\ge0$, the cardinality of $\NCL(n+1)$ is the large Schr\"oder number $r_n$.
\end{thm}
\begin{proof}
From Proposition~\ref{prop:CpxT}, if $a$ is a random variable whose T--transform is
$T_{\Phi_a}=\alpha=(1,1,1,\ldots)$, then $\phi(a^n)=|\NCL(n)|$ for every $n\ge1$.
Treating $F=\Phi_a$ and $T_{\Phi_a}$ as analytic functions,
we have
\[
F(z)=\sum_{n=0}^\infty|\NCL(n+1)|z^n
\]
and thus by convention~\eqref{eq:DF},
\[
D_F(w)=\sum_{n=1}^\infty|\NCL(n)|w^n.
\]
On the other hand
\[
T_F(z)=\sum_{n=0}^\infty z^n=\frac1{1-z}
\]
and from the relation~\eqref{eq:Tdef} we get
\[
D_F\cinv(z)=\frac{z(1-z)}{1+z}.
\]
Solving for $D_F$ yields
\[
F(w)=\frac1wD_F(w)=\frac{1-w-\sqrt{1-6w+w^2}}{2w},
\]
which is the generating function for the large Schr\"oder numbers.
\end{proof}

The first several values of $|\NCL(n)|$ are listed in Table~\ref{tab:cardNCNP}.
\begin{table}[hb]
\caption{Cardinalities of $\NCL(n)$.}
\label{tab:cardNCNP}
\begin{tabular}{r|c|c|c|c|c|c|c|c}
$n$\hspace{0.8em} & $1$ & $2$ & $3$ & $4$ & $5$ & $6$ & $7$ & $8$ \\ \hline
$\card(\NCL(n))$ & $1$ & $2$ & $6$ & $22$ & $90$ & $394$ & $1806$ & $8558$
\end{tabular}
\end{table}

Combining Theorem~\ref{thm:cardNCL} and equation~\eqref{eq:NCLNC}
of Corollary~\ref{cor:NCLNC}, we obtain the following expression for
the large Schr\"oder numbers.

\begin{cor}
For every $n\ge1$,
\[
r_n=\sum_{\pi\in\NC(n+1)}\;\prod_{B\in\pi}c_{|B|-1},
\]
where as in~\eqref{eq:ck}, $c_k$ is the Catalan number if $k\ge1$ and $1$ if $k=0$.
\end{cor}

\bibliographystyle{plain}

\end{document}